\documentclass[12pt,a4paper,reqno, UTF-8]{amsart}
\usepackage{amscd,amssymb,amsfonts,amsbsy,amsmath,verbatim,color,amsaddr,mathrsfs}
\usepackage{graphicx,cite,array}
\usepackage{float}
\usepackage{tikz,ifthen,bm}
\usetikzlibrary{intersections}
\usetikzlibrary{calc}

\newtheorem{thm}{Theorem}[section]
\newtheorem{theo}[thm]{Theorem}

\newtheorem{conj}[thm]{Conjecture}
\newtheorem{lem}[thm]{Lemma}
\newtheorem{ques}[thm]{Question}
\newtheorem{prop}[thm]{Proposition}
\newtheorem{exam}[thm]{Example}

\theoremstyle{definition}

\newtheorem{rema}[thm]{Remark}
\numberwithin{equation}{section}

\newcommand{\myint}[1][v]{\textup{int}_{#1}\,}

\newcommand{\D}{{\mathcal D}}

\newcommand{\cL}{{\mathcal L}}

\def\cP{{\mathcal P}}

\newcommand{\Real}{\mathbb{R}}

\newcommand{\R}{\Real}
\newcommand{\Z}{{\mathbb Z}}

\newcommand{\N}{{\mathbb N}}

\def\sign{\textup{sign}\,}

\def\cP{\mathcal{P}}

\def\bb{{\boldsymbol b}}
\def\bi{{\boldsymbol i}}
\def\bj{{\boldsymbol j}}
\def\bk{{\boldsymbol k}}

\def\ba{{\boldsymbol a}}

\def\bx{{\boldsymbol x}}
\def\by{{\boldsymbol y}}
\def\bz{{\boldsymbol z}}
\def\bu{{\boldsymbol u}}
\def\bv{{\boldsymbol v}}
\def\bw{{\boldsymbol w}}
\def\bs{{\boldsymbol s}}
\def\bt{{\boldsymbol t}}
\def\brho{{\bm \rho}}
\def\bone{\textbf{1}}

\def\biota{\boldsymbol{\iota}}
\def\ep{\varepsilon}
\def\eps{\ep}

\def\dis{\displaystyle}
\def\pt{\partial}

\def\diag{\textup{diag}}

\newcommand{\setcolsep}[1][2]
    {\setlength{\arraycolsep}{#1 pt}}
\def\newdot{\boxdot}

\makeatletter
\tikzset{%
    declare function={%
        ix(\t)={(2-\t)*\hx+(\t-1)*\rx}; iy(\t)={(2-\t)*\hy+(\t-1)*\ry};
        fx(\t)=(1-\hx)*\t/\h+\hx;  fy(\t)=(1-\hy)*\t/\h+\hy;
        delt(\t)=sqrt(1-\t)*\del;  ct(\t)=(delt(\t)+.5);
        utx(\t)=(1-\ux+\ux*\t);    uty(\t)=(1-\uy+\uy*\t);
        uutx(\t)=(1-\uux+\uux*\t); uuty(\t)=(1-\uuy+\uuy*\t);
        ax(\t)=(\axtwo/2-\axone)*\t*\t+(2*\axone-\axtwo/2)*\t;
        ay(\t)=(\aytwo/2-\ayone)*\t*\t+(2*\ayone-\aytwo/2)*\t;
        }}
\tikzset{dot2/.style={dash pattern=on 2pt off 1.2pt, very thin}}
\tikzset{red2/.style={top color=yellow!20!red!30,bottom color=red}}
\tikzset{red1/.style={left color=yellow!20!red!30,right color=red}}
\tikzset{blue2/.style={top color=yellow!20!blue!30,bottom color=blue}}
\tikzset{blue1/.style={left color=yellow!20!blue!30,right color=blue}}
\tikzset{%
    axistile/.style={>=stealth,xscale=2.4,yscale=2.5,%
    x={(.25cm,-.2cm)},y={(.1*\ct cm,.1cm)},z={(0cm,1cm)}
    }}
\tikzset{%
    iter1/.style={>=stealth,xscale=2.8,yscale=3,
        x={(-.28cm,-.1cm)}, y={(1cm,0cm)},z={(0cm,1cm)}}
}
\def\mystyle#1
    {\ifnum#1=1
    \def\sty@b{dot2}
    \def\sty@t{dot2}
    \def\sty@l{very thin}
    \def\sty@m{dot2}
    \def\sty@r{very thin}
    \else
    \def\sty@b{very thin}
    \def\sty@t{very thin}
    \def\sty@l{very thin}
    \def\sty@m{very thin}
    \def\sty@r{very thin}
    \fi
    }
\def\Para#1{
    \ifnum#1=1
    \def\parax{\ux}\def\paray{\uy}\def\parhx{\hx}
    \def\parhy{\hy}
    \coordinate(py) at(-.5*\ux,-.5*\uy,0);
    \coordinate(rb)at($(o3)+.5*(\hx,\hy,-2/3)$);
    \coordinate(lt)at($(o3)-.5*(\hx,\hy,2/3)$);
    \else
    \def\parax{\uux}\def\paray{\uuy}\def\parhx{\rx}
    \def\parhy{\ry}
    \coordinate(py) at(-.5*\uux-\ux,-.5*\uuy-\uy,0);
    \coordinate(rb)at($(o6)+.5*(\rx,\ry,-4/3)$);
    \coordinate(lt)at($(o6)-.5*(\rx,\ry,4/3)$); 
    \fi}
\pgfmathsetmacro{\ct}{9}
\def\mystart{\@ifnextchar[{\mystarti}{\mystarti[2]}}
\def\mystarti[#1]{\@ifnextchar[{\mystartii[#1]}{\mystartii[#1][.08]}}
\def\mystartii[#1][#2]{
    \pgfmathsetmacro{\r}{3}
    \pgfmathsetmacro{\v}{1.6}
    \pgfmathsetmacro{\h}{#1}
    \pgfmathsetmacro{\eps}{#2}
    \pgfmathsetmacro{\axone}{-.5}
    \pgfmathsetmacro{\ayone}{-.5}
    \pgfmathsetmacro{\axtwo}{-1.15}
    \pgfmathsetmacro{\aytwo}{-1.15}
    \pgfmathsetmacro{\sone}{2}
    \pgfmathsetmacro{\stwo}{1.8}
    \pgfmathsetmacro{\ux}{abs(\sone/6-\axone/3+\axtwo/6)}  
    \pgfmathsetmacro{\uy}{abs(\stwo/6-\ayone/3+\aytwo/6)}  
    \pgfmathsetmacro{\uux}{abs(\sone/6+\axone/3-\axtwo/6)}  
    \pgfmathsetmacro{\uuy}{abs(\stwo/6+\ayone/3-\aytwo/6)}  
    \pgfmathsetmacro{\hx}{(1-\ux)}      
    \pgfmathsetmacro{\hy}{(1-\uy)}      
    \pgfmathsetmacro{\rx}{(1-\uux)}      
    \pgfmathsetmacro{\ry}{(1-\uuy)}      
    \coordinate (c1) at (0,-.5*\h,2*\v);
    \coordinate (c2) at (0,.5*\h-.1,2*\v);
    \coordinate (c3) at (0,1.5*\h-.2,2*\v);
    \coordinate (c5) at (0,-.2,\v);
    \coordinate (c4) at (0,\h+.2,\v);
    \coordinate (c6) at (0,0,0);
    \coordinate (c7) at (0,\h,0);
    \coordinate (c8) at (0,2*\h,0);
    \coordinate (eps) at (0,0,\eps);
    \coordinate (reps) at (0,0,1/3-\eps);
    \coordinate (reps2) at (0,0,1/3-2*\eps);
    \coordinate (o) at (0,0,0);
    \coordinate (o1) at (0,0,0);
    \coordinate (o3) at ($-.5*(\ux,\uy,-2/3)$);
    \coordinate (o2) at ($(o3)+.5*({1-ix(1)},{1-iy(1)},0)-(eps)$);
    \coordinate (o4) at ($(o3)-.5*({1-ix(1)},{1-iy(1)},0)+(eps)$);
    \coordinate (o6) at (-\ux-.5*\uux,-\uy-.5*\uuy,2/3);
    \coordinate (o5) at ($(o6)+.5*({1-ix(2)},{1-iy(2)},0)-(eps)$);
    \coordinate (o7) at ($(o6)-.5*({1-ix(2)},{1-iy(2)},0)+(eps)$);
    \coordinate (o8) at ($(o7)+(0,0,1/3-\eps)$);
    \coordinate (oo1) at (-\ux/2,-\uy/2,1/3);
    \coordinate (oo2) at (-\ux-\uux/2,-\uy-\uuy/2,2/3);
    \coordinate (t1) at (\axtwo/2+\sone/2,\aytwo/2+\stwo/2,1);
    \coordinate (t2) at (\axtwo/2+5*\sone/6,\aytwo/2+5*\stwo/6,1);
    \coordinate (t3) at ($(t1)+2/3*(\sone,\stwo,0)$);
    \coordinate (t4) at ($(t1)+(\sone,\stwo,0)$);
    \coordinate (t5) at ($(t1)+(4/3*\sone,4/3*\stwo,0)$);
    }

\def\rectreg#1#2#3#4{%
    \path[#4]($(#1)+.5*(#2,-#3,0)$)--++(-#2,0,0)--++(0,#3,0)--++(#2,0,0)--cycle;
    }
\def\rectline#1#2#3#4#5{%
    \draw[#4]($(#1)+.5*(#2,-#3,0)$)--++(-#2,0,0)--++(0,#3,0);
    \draw[#5]($(#1)+.5*(#2,-#3,0)$)--++(0,#3,0)--++(-#2,0,0);
    }
\def\rect#1#2#3#4{%
    \@ifnextchar[%
        {\rect@i{#1}{#2}{#3}{#4}}
        {\rect@i{#1}{#2}{#3}{#4}[dot2]}
    }
\def\rect@i#1#2#3#4[#5]{%
    \@ifnextchar[%
        {\rect@ii{#1}{#2}{#3}{#4}[#5]}
        {\rect@ii{#1}{#2}{#3}{#4}[#5][very thin]
        }
    }
\def\rect@ii#1#2#3#4[#5][#6]{%
    \rectreg{#1}{#2}{#3}{#4}
    \rectline{#1}{#2}{#3}{#5}{#6}
    }
\def\myrect#1#2{\draw[fill=#2,opacity=.5](#1)--++(1,0,0)--++(0,1,0) --++(-1,0,0)--cycle}

\def\cub[#1]#2{%
    \draw[fill=yellow](#1)--++(1,0,0) --++(#2)--++(-1,0,0)--cycle;
    \draw[fill=red!40]($(#1)+(1,0,0)$)--++(0,1,0) --++(#2)--++(0,-1,0)--cycle;
    \rect{$(#1)+(#2)+(.5,.5,0)$}11{fill=cyan!70}
    }

\def\myycub#1[#2][#3]{
    \foreach \cub@i/\cub@sty in {0/dot2,1/very thin}
        {\path($(#1)+(.5,-.5,\cub@i)$)--++(-1,0,0)--++(0,1,0) --++(1,0,0)--cycle;
        \draw[\cub@sty]($(#1)+(.5,-.5,\cub@i)$)--++(-1,0,0)--++(0,1,0);
        \draw[very thin]($(#1)+(.5,-.5,\cub@i)$)--++(0,1,0)--++(-1,0,0);
        }
    \path[fill=#2,opacity=.5]($(#1)+#3*.5*(1,-1,0)$) --++({-#3},0,0)--++(0,0,1) --++({#3},0,0)--cycle;
    \path[fill=#2,opacity=.5]($(#1)+#3*.5*(1,-1,0)$) --++(0,{#3},0)--++(0,0,1) --++(0,{-#3},0)--cycle;
    \foreach \cub@i/\cub@sty in {{-1,-1,0}/dot2,{1,-1,0}/very thin, {1,1,0}/very thin, {-1,1,0}/very thin}
        \draw[\cub@sty]($(#1)+.5*(\cub@i)$)--++(0,0,1);
    }
\def\lateral{\@ifstar{\lateraldown}{\lateralup}}
\def\lateraldown#1#2{%
    \@ifnextchar[%
        {\lateraldowni{#1}{#2}}
        {\lateraldowni{#1}{#2}[\prm@i/3-\eps]}
    }
\def\lateraldowni#1#2[#3]{
    \@ifnextchar[%
        {\lateraldownii{#1}{#2}[#3]}
        {\lateraldownii{#1}{#2}[#3][opacity=1]}
    }

\def\lateraldownii#1#2[#3][#4]{
    \coordinate (r@tmp) at ($(#2)+(\parhx/2,\parhy/2,-{(#3)})+{(#3)}*(#1)$);
    \rect{$(r@tmp)-(\parhx/2,\parhy/2,0)$}\parhx\parhy{#4,fill=blue!40} [dot2][dot2]  ;
    \draw[#4,dot2,blue2](r@tmp)--++(-\parhx,0,0)--++
        ($(0,\paray,0)+\eps*(#1)$)--++(1,0,0)--cycle;
    \draw[#4,dot2,blue2](r@tmp)--++(0,-\parhy,0)--++
        ($(\parax,0,0)+\eps*(#1)$)--++(0,1,0)--cycle;
    \@ifundefined{prm@uncover@line}{\def\prm@uncover@line{1}}{}
    \ifnum\prm@uncover@line=1
        \draw[dot2](r@tmp)--++(-\parhx,0,0)--++(0,-\parhy,0)--++(\parhx,0);
    \fi
    }

\def\lateralup#1#2{%
    \@ifnextchar[%
        {\lateralupi{#1}{#2}}
        {\lateralupi{#1}{#2}[\prm@i/3+\eps]}
    }
\def\lateralupi#1#2[#3]{
    \@ifnextchar[%
        {\lateralupii{#1}{#2}[#3]}
        {\lateralupii{#1}{#2}[#3][opacity=1]}
    }
\def\lateralupii#1#2[#3][#4]{
    \coordinate (l@tmp) at ($(#2)-(\parhx/2,\parhy/2,#3)+{(#3)}*(#1)$);
    \draw[#4,dot2,red2](l@tmp)--++(0,\parhy,0)--++
        ($(-\parax,0,0)-\eps*(#1)$)--++(0,-1,0)--cycle;
    \draw[#4,dot2,red1](l@tmp)--++(\parhx,0,0)--++
        ($(0,-\paray,0)-\eps*(#1)$)--++(-1,0,0)--cycle;
    \rect{$(l@tmp)+(\parhx/2,\parhy/2,0)$}\parhx\parhy{#4,fill=red!50}[dot2]
    \draw[very thin]($(l@tmp)+(\parhx,0,0)$)--++($(0,-\paray,0)-\eps*(#1)$);
    \draw[very thin]($(l@tmp)+(0,\parhy,0)$)--++($(-\parax,0,0)-\eps*(#1)$);
    \@ifundefined{prm@uncover@line}{\def\prm@uncover@line{1}}{}
    \ifnum\prm@uncover@line=1
        \draw[dot2]($(l@tmp)+(-\parax,\parhy-1,0)-\eps*(#1)$)--++(0,1,0);
    \fi
    }

\def\tile#1#2{\def\para{#2}
    \ifcase\para 
        \coordinate (b) at (#1);
        \cub[b]{t1};
        \or  
        \foreach \i in {0,1,2}
            {
            \coordinate (b) at
                ($(#1)+({ax(\i)/3+\sone*(\i/9)},
                {ay(\i)/3+\stwo*(\i/9)},
                \i/3)$);
              \cub[b]{$1/3*(t2)$}
              }
            \or 
        \foreach \i in {0,1,2}
          \foreach \j in {0,1,2}
              {
              \coordinate (b) at
                ($(#1)+({ax(\i)/3+ax(\j)/9+\sone*(\i/9+2*\j/27)},
                {ay(\i)/3+ay(\j)/9+\stwo*(\i/9+2*\j/27)},
                \i/3+\j/9)$);
              \cub[b]{$1/9*(t3)$}
              }
              \or 
        \foreach \i in {0,1,2}
          \foreach \j in {0,1,2}
            \foreach \k in {0,1,2}
              {
              \coordinate (b) at
                ($(#1)+({ax(\i)/3+ax(\j)/9+ax(\k)/27+\sone*(\i/9+2*\j/27+\k/27)},
                {ay(\i)/3+ay(\j)/9+ay(\k)/27+\stwo*(\i/9+2*\j/27+\k/27)},
                \i/3+\j/9+\k/27)$);
              \cub[b]{$1/27*(t4)$}
              }
            \else 
        \foreach \i in {0,1,2}
          \foreach \j in {0,1,2}
            \foreach \k in {0,1,2}
              \foreach \l in {0,1,2}
              {
              \coordinate (b) at
                ($(#1)+({ax(\i)/3+ax(\j)/9+ax(\k)/27+ax(\l)/81 +\sone*(\i/9+2*\j/27+\k/27+4*\l/243)},
                {ay(\i)/3+ay(\j)/9+ay(\k)/27+ax(\l)/81 +\stwo*(\i/9+2*\j/27+\k/27+4*\l/243)},
                \i/3+\j/9+\k/27+\l/81)$);
              \cub[b]{$1/18*(t5)$}
              }
        \fi}

\def\prism{%
    \@ifstar{\@ifnextchar[{\prismm{t1}}{\prismm{t1}[{0,0,0}]}}
        {\@ifnextchar[{\prismm{0,0,1}}{\prismm{0,0,1}[{0,0,0}]}}}
\def\prismm#1[#2]{ 
    \rect{#2}11{fill=gray!70}
    \foreach \prm@i in {1,2}
        \foreach \prm@j/\prm@col in {-1/blue!70,0/yellow!50,1/red!50}
            \rect{$(#2)+{(\prm@i/3+\prm@j*\eps)}*(#1)$}11{fill=\prm@col}
    \rect{$(#2)+(#1)$}11{fill=cyan!50}[very thin]
    \foreach \prm@i in {%
        {(1/3-\eps)},{1/3},{(2/3-\eps)}, {2/3}}
        \draw[dot2]($(.5,-.5,0)+(#2)+\prm@i*(#1)$)--++(-1,0,0)--++(0,1,0);
    \foreach \prm@i in {%
        0,{1/3},{2/3},{1/3-\eps},{1/3+\eps}, {2/3-\eps},{2/3+\eps},1}
        \draw[fill=black]($(#2)+\prm@i*(#1)$)circle(.02);
    \draw[dashed](#2)--++(#1);
    \foreach\prm@i in {{.5,-.5,0}, {.5,.5,0}, {-.5,.5,0}}
        \draw[thin]($(\prm@i)+(#2)$)--++(#1);
    \draw[dot2]($(#2)-.5*(1,1,0)$)--++(#1);
    }

\def\shrinkprism{%
    \@ifstar{\@ifnextchar[{\shrinkprismm{1}}{\shrinkprismm{1}[{0,0,0}]}}
        {\@ifnextchar[{\shrinkprismm{0}}{\shrinkprismm{0}[{0,0,0}]}}}
\def\shrinkprismcenter#1#2{
    \foreach \prm@i/\prm@j in {
        0/o1,{(1/3-\eps)}/o2,{1/3}/o3,{(1/3+\eps)}/o4, {(2/3-\eps)}/o5,{2/3}/o6,{(2/3+\eps)}/o7,1/o8}
        \draw[fill=black]($(#2)+(0,0,{\prm@i-#1*\prm@i})+#1*(\prm@j)$) circle(.02);
    \draw[dashed](#2)--++(reps)--++($2*(eps)-#1*(\ux,\uy,0)$) --++(reps2)--++($2*(eps)-#1*(\uux,\uuy,0)$)--++(reps);}
\def\shrinkprismm#1[#2]{
    \rect{#2}11{fill=gray!50}
    \foreach \p@i in {1,2}{%
        \Para{\p@i}
        \foreach \p@j/\p@x/\p@y/\p@col in {%
            -1/1/1/blue!70,0/\parhx/\parhy/yellow!50, 1/1/1/red!50}{%
            \coordinate (prm@tmp) at (-#1*\p@j*\parax/2,-#1*\p@j*\paray/2,\p@i/3+\p@j*\eps);
            \rect{$(#2)+#1*(py)+(prm@tmp)$}{\p@x}{\p@y}{fill=\p@col}
            }
        }
    \rect{$(#2)+#1*(-\ux-\uux,-\uy-\uuy,0)+(0,0,1)$}11{fill=cyan!50}[very thin]
    \foreach \prm@i in {1,2}{%
        \Para{\prm@i}
        \foreach \prm@j/\prm@ix/\prm@iy in {-1/1/1,0/\parhx/\parhy}
            \draw[dot2] ($(#2)+#1*(py)+(.5*\prm@ix-\prm@j*\parax/2*#1, -.5*\prm@iy-\prm@j*\paray/2*#1, \prm@i/3+\prm@j*\eps)$) --++(-\prm@ix,0,0)--++(0,\prm@iy,0);
        }
    \foreach \prm@i/\prm@j in {
        0/o1,{(1/3-\eps)}/o2,{1/3}/o3,{(1/3+\eps)}/o4, {(2/3-\eps)}/o5,{2/3}/o6,{(2/3+\eps)}/o7,1/o8}
        \draw[fill=black]($(#2)+(0,0,{\prm@i-#1*\prm@i})+#1*(\prm@j)$) circle(.02);
    \shrinkprismcenter{#1}{#2};
    \foreach \prm@i/\prm@j/\prm@sty in {1/-1/thin,1/1/thin,-1/1/thin,-1/-1/dot2}
        \draw[\prm@sty]($(#2)+.5*(\prm@i,\prm@j,0)$)--++(reps)-- ++($.5*(-\prm@i*\ux,-\prm@j*\uy,2*\eps)+#1*.5*(-\ux,-\uy,0)$) -- ++($.5*(\prm@i*\ux,\prm@j*\uy,2*\eps)+#1*.5*(-\ux,-\uy,0)$) -- ++(reps2)--++ ($.5*(-\prm@i*\uux,-\prm@j*\uuy,2*\eps)+#1*.5*(-\uux,-\uuy,0)$) --++($.5*(\prm@i*\uux,\prm@j*\uuy,2*\eps)+#1*.5*(-\uux,-\uuy,0)$) --++(reps);
        }

\def\varprism{%
    \@ifstar{\varprismi{t2}}{\varprismi{0,0,1}}}
\def\varprismi#1{
    \@ifnextchar[
        {\varprismii{#1}}
        {\varprismii{#1}[{3,0,0}]}}
\def\varprismii#1[#2]{%
    \@ifnextchar[
        {\varprismiii{#1}[#2]}
        {\varprismiii{#1}[#2][1]}}

\def\varprismiii#1[#2][#3]{
    \ifnum#3=1\def\prm@uncover@line{1}\else\def\prm@uncover@line{0}\fi
    \rect{#2}11{fill=gray!50}
    \foreach \prm@i/\prm@j in {1/o3,2/o6}{
        \Para{\prm@i}
        \coordinate (r) at ($(rb)+{(\prm@i/3-\eps)}*(#1)$);
        \coordinate (l) at ($(lt)+{(\prm@i/3+\eps)}*(#1)$);
        \lateral*{#1}{$(#2)+(\prm@j)+(0,0,-\eps)$}
        \rect{$(#2)-(0,0,\prm@i/3)+(\prm@j)+\prm@i/3*(#1)$} \parhx\parhy{fill=yellow!50}
        \lateral{#1}{$(#2)+(\prm@j)+(0,0,\eps)$}
        }
    \rect{$(#2)+(o8)+(0,0,-1)+(#1)$}11{fill=cyan!50}[very thin]
    \foreach \prm@i/\prm@k in {{0,0,0}/0,{-\ux,-\uy,0}/1, {-\ux-\uux,-\uy-\uuy,0}/2}
        {\foreach \prm@j in {{.5,-.5,0},{.5,.5,0},{-.5,.5,0}}
            \draw[thin]($(\prm@i)+(\prm@j)+(#2)+{\prm@k/3}*(#1)$) --++($1/3*(#1)$);
            \ifnum#3=1
                \draw[dot2]
                ($(\prm@i)+(-.5,-.5,0)+(#2)+{\prm@k/3}*(#1)$)--++($1/3*(#1)$);
                \fi}
    \foreach \prm@i/\prm@j/\prm@ii/\prm@ix/\prm@iy in {o4/o3/1/\hx/\hy,o7/o6/2/\rx/\ry}{
        \draw[thin]($(#2)+(\prm@i)+(.5,-.5,-\prm@ii/3-\eps)+\prm@ii/3*(#1)$) --++(0,1-\prm@iy,0) --++(1-\prm@ix,0)--++(0,1,0) --++(-1,0,0)--++(0,\prm@iy-1,0)--++(\prm@ix-1,0,0);
            }
        \ifx#31
        \foreach \prm@i/\prm@ii/\prm@j in {%
            o1/{0}/{0,0,0},
            o2/{(1/3-\eps)}/{0,0,-1/3+\eps},
            o3/{1/3}/{0,0,-1/3},
            o4/{(1/3+\eps)}/{0,0,-1/3-\eps},
            o5/{(2/3-\eps)}/{0,0,-2/3+\eps},
            o6/{2/3}/{0,0,-2/3},
            o7/{(2/3+\eps)}/{0,0,-2/3-\eps},
            o8/1/{0,0,-1}}
            \draw[fill=black]($(#2)+(\prm@i)+(\prm@j)+\prm@ii*(#1)$) circle(.02);

        \draw[dashed,thin]($(#2)+(o1)$) --++($(o2)-(o1)-(reps)+{(1/3-\eps)}*(#1)$) --++($(o3)-(o2)-(eps)+\eps*(#1)$)--++($(o4)-(o3)-(eps)+\eps*(#1)$) --++($(o5)-(o4)-(reps2)+{(1/3-2*\eps)}*(#1)$)
        --++($(o6)-(o5)-(eps)+\eps*(#1)$)
        --++($(o7)-(o6)-(eps)+\eps*(#1)$)
        --++($(o8)-(o7)-(reps)+{(1/3-\eps)}*(#1)$);
        \fi
    }

\def\mysectreg{\@ifstar{\mysectregi[\uux][\uuy]}{\mysectregi[\ux][\uy]}}
\def\mysectregi[#1][#2]#3#4#5#6#7{
    \coordinate(Ga)at({(#1-1)/2-(#4)*(#1)},{(1-(#2))/2},{#4*\eps});
    \path[#7,fill=#5](${(#6)}*(Ga)+(#3)$)--++(0,{-(#6)*(1-(#2)+(#2)*(#4))},0)
        --++({#6*(1-#1+#1*#4)},0,0) --++($#6*({(ct(#4)-.5+#1-#4*#1)}, {(-ct(#4)+.5-#2+#4*#2)},0)$)
        --++($#6*({-2*ct(#4)-#1+#4*#1},0,0)$) --++($#6*(0,{2*ct(#4)+#2-#4*#2},0)$)
        --cycle;
        }
\def\colorgray{gray!80}
\def\mysectreggray{%
    \@ifstar{\mysectreggrayi[\uux][\uuy]}
        {\mysectreggrayi[\ux][\uy]}}
\def\mysectreggrayi[#1][#2]#3#4#5#6{
    \path[fill=\colorgray,opacity=.7]%
        ($({-.5-(#1)/2-\del},{.5+(#2)/2+\del},0)+(#3)$)--++ (0,{-1-(#2)-2*\del},0)--++({1+#1+2*\del},0,0) --($(.5*#4,-.5*#5,0)+(#3)$)--++(-#4,0,0) --++(0,#5,0)-- ($({-.5-(#1)/2-\del},{.5+(#2)/2+\del},0)+(#3)$)--++
        ({1+#1+2*\del},0,0)--++(0,{-1-(#2)-2*\del},0) --($(.5*#4,-.5*#5,0)+(#3)$)--++(0,#5,0)--++(-#4,0,0)--cycle;
    }

\def\mysectline{%
    \@ifstar{\mysectlinei[\uux][\uuy]}
        {\mysectlinei[\ux][\uy]}}
\def\mysectlinesty{0}
\def\mysectlinei[#1][#2]#3#4#5#6#7{%
    \coordinate (myo) at (#3,#4,#5);
    \ifnum\mysectlinesty=1
        \draw[opacity=.5,fill=#6,domain=0:1]%
        plot({#7*(ct(\x)+#1/2-\x*#1)+#3}, {#7*(-ct(\x)-#2/2)+#4}, {#7*\x*\eps+#5})--
        plot({#7*(-ct(1-\x)-#1/2)+#3},{#7*(-ct(1-\x)-#2/2)+#4},
        {#7*(1-\x)*\eps+#5})--cycle;
    \draw[opacity=.5,fill=#6,domain=0:1]
        plot({#7*(-ct(\x)-#1/2)+#3},{#7*(-ct(\x)-#2/2)+#4},
        {#7*(\x)*\eps+#5})--
        plot({#7*(-ct(1-\x)-#1/2)+#3},{#7*(ct(1-\x)+#2/2-#2*\x)+#4}, {#7*(1-\x)*\eps+#5})--cycle;
    \else
    \draw[opacity=.5,#6,domain=0:1,thick]%
        plot({#7*(ct(\x)+#1/2-\x*#1)+#3}, {#7*(-ct(\x)-#2/2)+#4}, {#7*\x*\eps+#5})
        plot({#7*(-ct(\x)-#1/2)+#3},{#7*(-ct(\x)-#2/2)+#4},
        {#7*\x*\eps+#5})
        plot({#7*(-ct(\x)-#1/2)+#3},{#7*(ct(\x)+#2/2-#2*\x)+#4}, {#7*\x*\eps+#5});
    \foreach \myrect@j in {.1,.2,...,1}
        \draw[help lines,#6,very thin] ($#7*({(ct(\myrect@j)+#1/2-\myrect@j*#1)},
            {(-ct(\myrect@j)-#2/2)},{\myrect@j*\eps})+(myo)$) --++($#7*({-2*ct(\myrect@j)-#1+\myrect@j*#1},0,0)$) --++($#7*(0,{2*ct(\myrect@j)+#2-\myrect@j*#2},0)$);
    \fi
    }
\def\mySSplus{
    \@ifstar{\mySSplusi[\uux][\uuy]}{\mySSplusi[\ux][\uy]}
    }

\def\mySSplusi[#1][#2]#3#4#5{%
    \path[fill=#4,opacity=.5]($#5/2*(#1-1,#2-1,0)+(#3)$) --++(#5-#5*#1,0,0)--++(0,-#5*#2,#5*\eps) --++(-#5,0,0)--cycle;
    \path[fill=#4,opacity=.5]($#5*.5*(#1-1,#2-1,0)+(#3)$) --++(0,#5-#5*#2,0)--++(-#5*#1,0,#5*\eps)--++(0,-#5,0)--cycle;}
\def\SSplusline#1#2{
    \mystyle{#2}
    \coordinate (SStmp0)at($#2/2*(1-\parax,\paray-1,0)+(#1)$);
    \draw[\sty@b](SStmp0)--++(#2*\parax-#2,0,0)--++(0,#2-#2*\paray,0);
    \draw[\sty@t]($(SStmp0)+#2*(0,-\paray,\eps)$) --++(-#2,0,0)--++(0,#2,0);
    \draw[\sty@l](SStmp0)--++($#2*(0,-\paray,\eps)$);
    \draw[\sty@m]($(SStmp0)+#2*(\parax-1,0,0)$) --++($-#2*(\parax,\paray,-\eps)$);
    \draw[\sty@r]($(SStmp0)+#2*(\parax-1,1-\paray,0)$) --++($#2*(-\parax,0,\eps)$);
    }
\def\sectline{\mysectlinei[\parax][\paray]}
\def\sectreggray{\mysectreggrayi[\parax][\paray]}
\def\sectreg{\mysectregi[\parax][\paray]}
\def\SSplus{\mySSplusi[\parax][\paray]}
\def\mydashline#1{
    \@ifundefined{colorgray}{\def\colorgray{gray!70}}{}
    \@ifundefined{parax}{\def\parax{\ux}}{}
    \@ifundefined{paray}{\def\paray{\uy}}{}
    \draw[\colorgray,dot2] ($(#1)+(\del+1/2+\parax/2,-\del-1/2-\paray/2,0)$) --++(-2*\del-1-\parax,0,0)--++(0,2*\del+1+\paray,0);
    }
\def\Sextension#1{\@ifnextchar[{\Sextensioni#1}{\Sextensioni#1[0]}}

\def\Sextensioni#1[#2]#3#4{
    \def\colorgray{#1}

    \rect{o1}11{fill=gray}
    \foreach \myi/\sto/\ax/\ay/\az in
        {1/oo1/{-\ux/2}/{-\uy/2}/{1/3},
            2/oo2/{-\ux-\uux/2}/{-\uy-\uuy/2}/{2/3}}
        {
        \Para{\myi}
        \sectreg{\sto}{0}{red!30}{1}{opacity=.5} 
        \sectreg{\sto}{0}{blue!40}{-1}{opacity=.5}
        \ifx#21
            \coordinate (oo4) at(\ux-\myi*\ux,\uy-\myi*\uy,\myi/3-\eps);%
            \rectreg{oo4}11{fill=blue!80}
        \fi
        \SSplus{\sto}{blue!40}{-1}
        \sectline{\ax}{\ay}{\az}{blue!90}{-1} 
        \sectreggray{\sto}{#3}{#4}1  
        \@ifundefined{cutreg}{\def\cutreg{1}}{}
        \ifnum\cutreg=1
            \path[fill=white]($(\sto)+(-#3/2,#4/9,0)$)--++(0,#4/2-#4/9,0) --++(#3/5*2,0,0)--cycle;
            \path[fill=white]($(\sto)+(#3/2,-#4/9,0)$)--++(0,-#4/2+#4/9,0) --++(-#3/5*2,0,0)--cycle;
        \fi
        \ifx#21
            \coordinate (oo3) at
                (\uux-\myi*\uux-\ux,\uuy-\myi*\uuy-\uy,\myi/3+\eps);
            \rectreg{oo3}11{fill=red!50}
        \fi
        \sectline{\ax}{\ay}{\az}{red}{1} 
        \SSplus{\sto}{red!40}{1}
        \SSplusline{\sto}{1} 
        \SSplusline{\sto}{-1}
        }
    \foreach \i/\j/\sty in {1/-1/very thin,1/1/very thin,-1/1/very thin,-1/-1/dot2}
        \draw[\sty]($.5*(\i,\j,0)$)--++(reps) ++($.5*(-\i*\ux,-\j*\uy,2*\eps)+.5*(-\ux,-\uy,0)$) ++($.5*(\i*\ux,\j*\uy,2*\eps)+.5*(-\ux,-\uy,0)$) -- ++(reps2)++ ($.5*(-\i*\uux,-\j*\uuy,2*\eps)+.5*(-\uux,-\uuy,0)$) ++($.5*(\i*\uux,\j*\uuy,2*\eps)+.5*(-\uux,-\uuy,0)$) --++(reps);
    \foreach \i/\j in {{-\ux,-\uy,1/3}/{0,0,1/3}, {-\ux-\uux,-\uy-\uuy,2/3}/{-\ux,-\uy,2/3}}
        {\draw[very thin]($.5*(1,1,0)+(\i)$)--++(0,0,\eps);
        \draw[dot2]($-.5*(1,1,0)+(\j)$)--++(0,0,-\eps);}
        \rect{o8}11{fill=cyan}
    }
\def\mybraceht{.05}\def\col{black}
\newcommand{\mybrace}[4][-1]
    {\ifthenelse{#1=1}{\def\loc{above}}{\def\loc{below}}
    \draw[\col](#2)sin++({#3*.1},{#1*\mybraceht})--++({#3*.35},0) --++({#3*.05},{#1*\mybraceht*.5})node[yshift=-#1*3pt,\loc,font=\tiny]{#4} --++({#3*.05},{#1*(-\mybraceht*.5)}) --++({#3*.35},0) cos++({#3*.1},{#1*(-\mybraceht)})--
    ++({-#3},0);}

\makeatother
\begin{document}
\title[Self-affine tiles that are tame balls]
{A class of self-affine tiles in $\R^d$ that are $d$-dimensional tame balls \myvp}
\date{}\def\myvp{\vspace*{-15pt}}
\author[G.T. Deng]{Guotai Deng\myvp}
\address{School of Mathematics and Statistics \& Hubei Key Laboratory of Mathematical Sciences, Central China Normal University, Wuhan 430079, P.R. China \myvp}
\email{hilltower@163.com}

\author[C.T. Liu]{Chuntai Liu \myvp}
\address{School of Mathematics and Computer Science, Wuhan Polytechnic University, Wuhan 430023, P. R. China \myvp}
\email{lct984@163.com}

\author[S.-M. Ngai]{Sze-Man Ngai \myvp}
    \address{Key Laboratory of High Performance Computing and Stochastic Information Processing (HPCSIP) (Ministry of Education of China), College of Mathematics and Statistics, Hunan Normal University, Changsha, Hunan 410081, P. R. China,\\ Department of Mathematical Sciences,
	Georgia Southern University Statesboro, GA 30460-8093, USA.}
\email{smngai@georgiasouthern.edu \myvp}

\thanks{%
    The first author is supported by the National Natural Science Foundation of China grant (12071171), Hubei Provincial Natural Science Foundation of China (2020CFB833), and Science Foundation of Jiangxi Education Department (GJJ202302). The second author is supported in part by the National Natural Science Foundation of China grant (11601403), and Hubei Provincial Natural Science Foundation of China (2019CFB602). The third author is supported in part by the National Natural Science Foundation of China grant (11771136), the Hunan Province' Hundred Talents Program, Construct Program of the Key Discipline in Hunan Province, and a Faculty Research Scholarly Pursuit Award from Georgia Southern University.}



\begin{abstract}
    We study a family of self-affine tiles in $\R^d$ ($d\ge\text{2}$) with noncollinear digit sets, which naturally generalizes a class studied originally by Deng and Lau in $\R^{\text{2}}$ and its extension to $\R^{\text{3}}$ by the authors. By using Brouwer's invariance of domain theorem, along with a tool which we call horizontal distance, we obtain necessary and sufficient conditions for the tiles to be $d$-dimensional tame balls. This answers positively the conjecture in an earlier paper by the authors stating that a member in a certain class of self-affine tiles is homeomorphic to a $d$-dimensional ball if and only if its interior is connected.
\end{abstract}
\maketitle
\setcounter{tocdepth}{1}
\tableofcontents

\section{Introduction}\label{S:intro}
    Let $d\ge1$ be an integer and $A$ be a $d\times d$ expanding matrix (i.e., all of its eigenvalues are greater than 1 in modulus). It is well known that (see, e.g., \cite{Hutchinson, Lagarias-Wang_1996}) for any finite set $\D\subset \R^{d}$  there exists a unique nonempty compact set $T=T(A,\D)$ such that
        \begin{equation*}\label{attractor}
          T=\bigcup\nolimits_{d_i\in\D}A^{-1}(T+d_i).
        \end{equation*}
    The above set equation can be rewritten as $AT=T+\D$ and $T$ can be expressed as
        \begin{equation}\label{texpression}
        T=\Bigl\{\sum\nolimits_{k\ge1}A^{-k}d_{i_k}:\,d_{i_k}\in \D\Bigr\}.
        \end{equation}
    We call $\D$ a {\it digit set}, $(A,\D)$ a {\it self-affine pair}, and $T$ a {\it self-affine set}. If $\#\D=|\det(A)|$ is an integer and the interior of $T$ (denoted by $T^\circ$) is nonempty, the set $T$ actually tiles $\R^d$ in the sense that there exists a discrete set $\cL\subset\R^d$ which satisfies (i) $T+\cL=\R^{d}$ and (ii) $(T^\circ+\iota_1)\cap(T^\circ+\iota_2)=\emptyset$ for all distinct $\iota_1,\iota_2\in\cL$. Such a $T$ is called a {\it self-affine tile}.

    The theory of self-affine tiles was established by Kenyon\cite{Kenyon}, Lagarias and Wang \cite{Lagarias-Wang_1996, Lagarias19962, Lagarias1997, wang} in the 1990s. In recent years, self-affine tiles play an important role in fractal geometry, wavelet theory and other fields. A basic central topic about self-affine sets is to study their topological properties and structure, including connectedness (see, e.g., \cite{Kirat2000, Aki2004, Kirat2004, Leung-Lau_2007, he, Leung-Luo_2012} and the references therein), and ball-likeness (see, e.g. \cite{Bandt-Wang_2001, Luo-Rao-Tan_2002, Leung-Lau_2007, Kamae, Conner, R3, Rn} and the references therein). Here a set $E\subset\R^d$ is said to be \textit{ball-like} if it is homeomorphic to a $d$-dimensional unit ball. If $d=2$, a ball-like set is also called a \textit{disk-like set}. All the balls mentioned in this paper are closed.

    When the integral expanding matrix $A\in\R^{d\times d}$ is given, the structure of the digit set strongly affects the topological properties of the self-affine set $T$. Hata \cite{Hata} obtained an important tool to study the connectedness of more general self-similar sets. By using a graph structure on $\D$, Kirat and Lau \cite{Kirat2000} developed the following criterion for the connectedness of self-affine tiles generated by collinear digit sets: If $\D=\{0,v,\ldots,|(\det A|-1)v\}$ and the characteristic polynomial of $A$ is height reducing, then $T=T(A,\D)$ is connected (see also \cite{Kirat2004}). As it is easy to describe all quadratic polynomials by this method, Kirat and Lau proved the connectedness of a class of tiles in $\R^2$.

    Akiyama and Gjini \cite{Aki2004} extended Kirat-Lau's conclusion to the case of $d=3,4$ by using the Schur-Cohn criterion. They confirmed that all characteristic polynomials are height reducing when $d\le4$. For $d>4$, He \textit{et al.} \cite{he} described an algorithm to determine whether a polynomial is height reducing. This method is not applicable if the digit set is nonconsecutive or noncollinear. A lot of research has focused on collinear digit sets or more special consecutive digit sets (see \cite{he, Kirat2000, LiuJC2014}). One of the reasons is that this case involves canonical number systems. It is a widely studied basic case (for canonical number systems and tiles see \cite{Gilbert, Katai, Kovacs}; for some progress on the topology of the tiles, see \cite{Aki2000, AkiThus2005}). Another reason is that in most cases, the connectedness can be revealed only by consecutive digits. If, on the other hand, a ``scattered" digit set is used, it is easy to find a disconnected sub-block, making it more difficult to study the topological structure of the tile.

    In 2011, Deng and Lau \cite{DengLau2011} initiated the study of self-affine tiles with noncollinear digit sets, and obtained the following result concerning connectedness. 
    \begin{thm}{\rm(Deng-Lau~\cite{DengLau2011})}
        Let $p,q\in\mathbb{Z}$ with  $|p|, |q|\geq 2$, $a\in \mathbb{R}$, and $b_i\in\R$, $i\in\{0,1,\dots,|p|-1\}$. Let
        \begin{equation}
            A\!=\!\begin{pmatrix}
                p& 0 \\[-3pt]
                -a&q
            \end{pmatrix},
            \D\!=\!\big\{(i,j+b_i):\, 0\le i\le |p|-1,\ 0\le j\le |q|-1\big\}.
        \end{equation}
    Then the self-affine set $T=T(A,\D)$ is a tile, and it is connected if and only if
        \begin{equation}
            \Big|\frac{b_{i+1}-b_i}{q}+\frac{\sign(p) (b_0-b_{|p|-1})-a}{q(q-\sign(p))}\Big|\le1, \text{ for all }i,
        \end{equation}
        where $\sign(x)=-1, 0$ or 1 according to $x<0$, $x=0$ or $x>0$.
    \end{thm}

    Leung and Luo \cite{Leung2015} showed the connectedness of the self-affine set $T$ in the case $\D=\{0,\bv,\pm A\bv\}$ ($\bv,A\bv$ are linearly independent), under the condition that $|\det A|=3$. The result in \cite{DengLau2011} has been extended by Liu \textit{et al.}~\cite{LiuJC2017} to the digit set $\D=\Z_m\times\Z_n$. Moreover, various other extensions in $\R^2$ have been obtained by a number of authors (see \cite{Ma-Dong-Deng_2014, Liu-Ngai-Tao_2016}).

    In the two-dimensional case, the topological structure of tiles has been studied extensively, and disk-like tiles have been identified in many cases. Due to the existence of Alexander's horned balls in $\R^3$, it is very difficult to verify whether a specific three-dimensional self-affine tile is ball-like. So, when $d\ge3$, the topology of self-affine tiles is very complex and little is known.
    In studying higher-dimensional tiles, Gelbrich \cite{Gelbrich} first proposed a part of the following question, which appeared recently in \cite{Bandtpre,DengJiang2012}.
    \begin{ques}\label{ques1}
        Can a self-affine tile \textup{(}or its boundary\textup{)} in $\R^d$ become a manifold in a nontrivial way, and is there a method to identify whether a self-affine tile (or its boundary) can be a manifold?
    \end{ques}

    We recall the definitions of tame ball and wild ball. Suppose a set $E\subset\R^d$ is homeomorphic to a $d$-dimensional unit ball. If the homeomorphism can be extended to $\R^d$, $E$ is called a \textit{tame ball}; otherwise $E$ is called a \textit{wild ball}.
    In 2015, Kamae, Luo and Tan \cite{Kamae} used Barnsley's theory of attractors with condensation \cite{Barnsley} to construct a class of $d$-dimensional self-affine tiles $T=T(A,\D)$ that are tame balls; here $A:=3\diag(1,\ldots,1)\in\R^{d\times d}$,
        $$\D:=\bigg\{\bigg(i_1,\ldots,i_n +t\Big(1-\sum_{j=1}^{n-1}i_j^2\Big)\bigg):\, i_j\in\{-1,0,1\},1\le j\le n\bigg\},$$
    and $t\in(-3,3)$. In 2016, Conner and Thuswaldner \cite{Conner} proposed an algebraic method to determine whether a given three dimensional self-affine tile is a tame ball, by using geometric topology tools developed by Cannon and Edwards. They confirmed that several known three-dimensional self-affine tiles \cite{BandtMessing, Gelbrich} are ball-like. They also constructed a self-affine tile in $\R^3$ that is not a topological manifold (its boundary is a two-dimensional wild sphere). Their method makes use of the representation of a class of algebraic geometric complexes with recursive structure (which allowed arbitrary precision calculations). These complexes were used to approximate self-affine tiles, and then determine the homeomorphism type of the boundary with the help of computer. In 2020, Thuswaldner and Zhang \cite{Thu} considered the boundary problem of self-affine tiles. Under the assumption that the digit set is consecutive, they constructed a class of self-affine tiles each of whose boundary is a two-dimensional sphere.

    In  \cite{R3}, authors of the present paper studied Question \ref{ques1} and extended some results in \cite{DengLau2011} to $\R^3$. By assuming that $A$ is a special lower triangle matrix, they constructed a family of ball-like self-affine tiles. The main technique is to restructure a self-affine tile by a cut-and-paste process, which allows one to construct a homeomorphism  from $T$ to a ball in $\R^3$. The following is a part of the main result in \cite{R3}.
    \begin{theo}\label{theomain2} {\rm (Deng-Liu-Ngai~\cite{R3})}
    Let $(A,\D)$ be given as follows:
        \begin{equation}\label{eqmatirxa'}
            A=
            \begin{pmatrix}
            p&0&0\\[-3pt]
            0&q&0\\[-3pt]
            -t&-s&r
        \end{pmatrix}\quad\text{ and }\quad
        \begin{array}{rcl}
        \D&=&\big\{(i,j,k):\,0\le i<|p|,\\
        &&\,\,0\le j<|q|,0\le k<|r|\big\},
        \end{array}
    \end{equation}
    where $p,q,r$ are integers not less than 2 in absolute value and $s,t$ are real numbers. If $st\ge0$, then $T=T(A,\D)$ is homeomorphic to a ball if and only if $|s+t|<|r(r-1)|$.
    \end{theo}

By assuming that $A$ is an upper triangular matrix,  authors of the present paper further extended some of the results in \cite{R3} to $\R^d$  for $d\ge 2$ \cite{Rn}. By constructing homotopy mappings, they obtained a necessary and sufficient condition for $T$ to be contractible (to a single point) and hence $T^\circ$ is simply connected.
    \begin{theo}\label{theoTcontractible}{\rm (Deng-Liu-Ngai~\cite{Rn})}
        For $d\ge 2$, let $p_1,\ldots,p_d$ be integers with $|p_i|\ge 2$ for $1\le i\le d$. Let  $A_1$ be the diagonal matrix $\diag(p_1,\ldots,p_{d-1})$ and $\bs=(s_1,\ldots,s_{d-1})\in\R^{d-1}$. Suppose the self-affine pair $(A,\D)$ satisfies
        \begin{equation}\label{eqexprad}
            A\!=\!
            \begin{pmatrix}
                A_1& -\bs\\[-3pt]
                0 & p_{d}
            \end{pmatrix},\quad
            \D\!=\!\big\{(i_1,\ldots,i_d):\,0\le i_j<|p_j|,\ 1\le j\le d\big\}.
        \end{equation}
    Then $T$ is a self-affine tile and the following statements are equivalent.\\
	\indent$(1)$ $T$ is connected if and only if
        \begin{equation}\label{eqcondition1}
            \max_{1\le j\le d-1}
            \left\{
                \left|
                \dfrac{s_j} {p_j(p_j-\sign(p_d))}
                \right|
            \right\}
            \le1.
        \end{equation}
    \indent$(2)$ $T^\circ$ is connected if and only if the inequality in \eqref{eqcondition1} is strict. Moreover, if $T^\circ$ is connected then $T^\circ$ is contractible to a single point.
    \end{theo}

    Although $T$ is contractible, one cannot deduce that $T$ is homeomorphic to a $d$-dimensional ball. The following is conjectured in \cite{Rn}.
    \begin{conj}\label{conj1}
       Assume the same hypotheses of Theorem \ref{theoTcontractible}. Then $T$ is homeomorphic to a $d$-dimension ball if and only if $T^\circ$ is connected.
    \end{conj}

    The main purpose of this paper is to give a positive answer to Conjecture \ref{conj1}. We use a purely analytic approach to show that the connectedness of $T^\circ$ implies that $T$ is a $d$-dimensional tame ball. The following is our main result.
    \begin{theo}\label{theomain1}
        For $d\ge 2$, let  $(A,\D)$ be given as in \eqref{eqexprad}, and let $T=T(A,\D)$ be the corresponding self-affine tile. Then the following statements are equivalent.
    \begin{enumerate}
        \item $T^\circ$ is connected;
        \item $T$ is a $d$-dimensional tame ball;
        \item the following inequality holds
            \begin{equation*}
                \max_{1\le j\le d-1}
                \left\{
                    \left|
                    \dfrac{s_j} {p_j(p_j-\sign(p_d))}
                \right|
                \right\}
                <1.
            \end{equation*}
    \end{enumerate}
    \end{theo}
    In fact, in Section~\ref{S:sec4}, we will prove a generalized version of Theorem~\ref{theomain1} by using more general digit sets (see Theorem~\ref{theo1}).

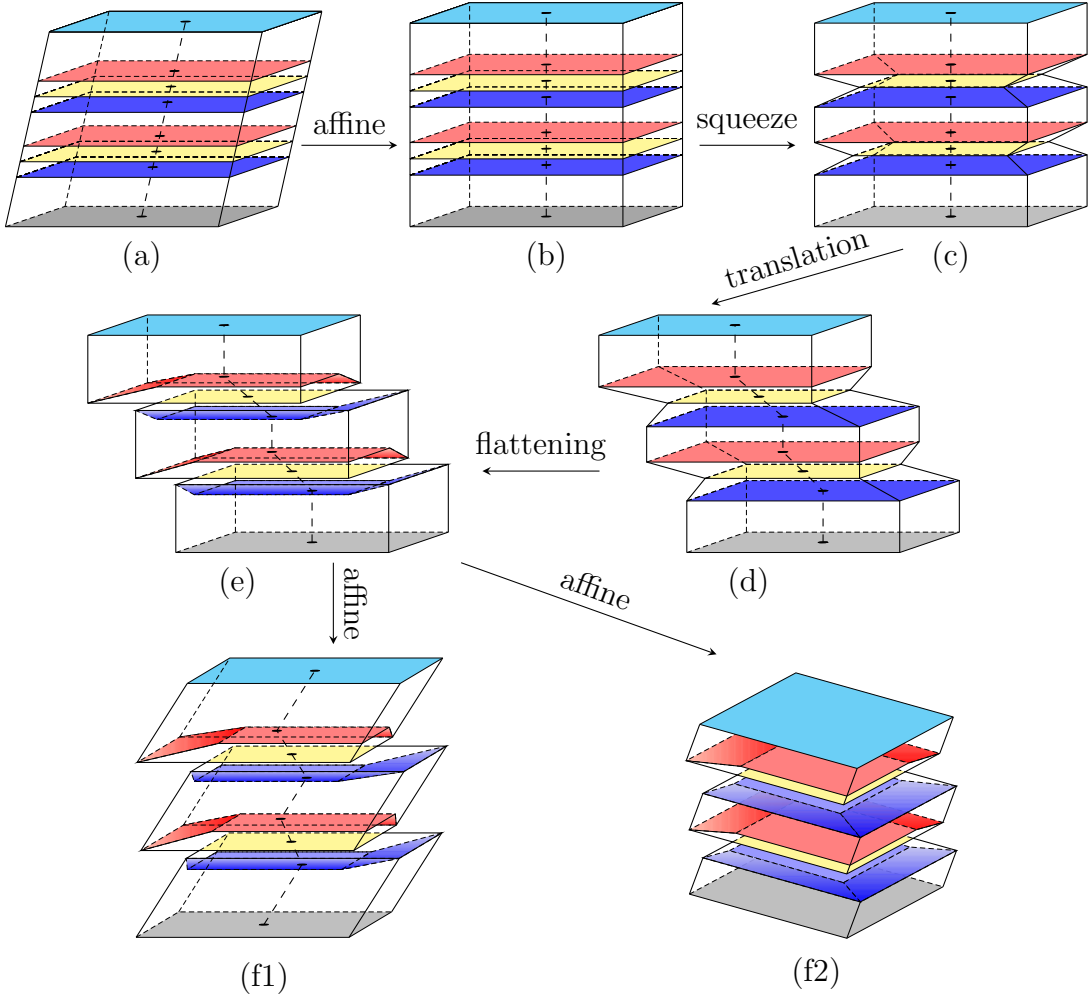
\begin{figure}
    \begin{tikzpicture}[{>=stealth,xscale=2.8,yscale=2.7,
        x={(-.28cm,-.1cm)}, y={(1cm,0cm)},z={(0cm,1cm)}}]
        \mystart[2]%
        \prism*[c1]
        \prism[c2]
        \shrinkprism[c3]
        \shrinkprism*[c4]
        \varprism[c5]
        \draw[->](0,-.25,2*\v+.35)--node[above]{affine}++(0,.45,0);
        \draw[->]($.5*(c2)+.5*(c3)+(0,-.23,.35)$)--node[above]{squeeze} ++(0,.45,0);
        \draw[->]($(c3)+(.8,0,-.08)$)--node[sloped,above]{translation\,\,} ++(0,-.9,-.27);     \draw[->]($.5*(c4)+.5*(c5)+(0,.15,.35)$)--node[above]{flattening} ++(0,-.55,0);
        \draw[->](0,-.1,\v-.1)--node[sloped,above]{affine}++(0,0,-.4);
        \draw[->](0,.5,\v-.1)--node[above,sloped]{affine}++(0,1.2,-.45);
        \foreach \i/\tex in {c1/(a),c2/(b),c3/(c)} \draw($(\i)+(0,0,-.2)$)node{\tex};
        \foreach \i/\tex in {c4/(d),c5/(e)}
            \draw($(\i)+(0,-.35,-.2)$)node{\tex};
    \end{tikzpicture}

    \begin{tikzpicture}[>=stealth,xscale=2.8,yscale=2.7, x={(-.2cm,-.125cm)}, y={(1cm,0cm)},z={(0cm,1.3cm)}]
        \mystart
        \varprism*[{0,0,0}]
        \draw(0,0,-.2)node{(f1)};
    \end{tikzpicture}
    \hskip 3cm
    \begin{tikzpicture}[{>=stealth,xscale=2.3,yscale=1.95, x={(-.55cm,-.35cm)}, y={(.9cm,-0.3cm)},z={(0cm,1.4cm)}}]
        \mystart
        \varprism*[{0,0,0}][0]
        \draw(0,0,-.4)node{(f2)};
        \draw(0,0,1)node{};\draw(0,0,-.5)node{};
    \end{tikzpicture}
    \caption{\label{figillus} Five homeomorphisms on a slant prism, from (a) to (f1) and (f2).  Figures (f1) and (f2) show the same set, with different viewpoints and with auxiliary dashed line segments in (f1). This figure is drawn with
    $p_1=p_2=p_3=3$, $\bs=(2,1.8)$ and $\D=\{(i_1+a_k^1,i_2+a_k^2,k):\,i_1,i_2,k=0,1,2\}$, where $(a_0^1,a_0^2)=(0,0),(a_1^1,a_1^2)=(-0.55,-0.45)$, and $(a_2^1,a_2^2) =(-1.25,-1.05)$.}\vspace*{-15pt}
    \end{figure}

\begin{figure}
    \begin{tikzpicture}[>=stealth,xscale=2,yscale=1.8,
        x={(.25cm,-.2cm)},y={(.9 cm,.1cm)},z={(0cm,1cm)}]
        \mystart
        \tile{0,0,0}0
        \draw($(o)+(0,.8,-.35)$)node[below]{(a) Original prism.};
        \coordinate (o) at (0,\h,-\h/\ct);
        \tile{o}1
        \draw($(o)+(0,.8,-.35)$)node[below]{(b) 1st iteration.};
        \coordinate (o) at (0,2*\h,-2*\h/\ct);
        \tile{o}2
        \draw($(o)+(0,.8,-.35)$)node[below]{(c) 2nd iteration.};
        \coordinate (o) at (0,.35*\h,-.35*\h/\ct-\v-.2);
        \tile{o}3
        \draw($(o)+(0,.8,-.35)$)node[below]{(d) 3rd iteration.};
        \coordinate (o) at (0,1.65*\h,-1.65*\h/\ct-\v-.2);
        \tile{o}4
        \draw($(o)+(0,.8,-.35)$)node[below,text width=105pt]{(e) 4th iteration.};
    \end{tikzpicture}
    \caption{\label{figtile} First few iterations on a slant prism in $\R^3$. The figures are drawn by using the data from Figure \ref{figillus}, but with a different viewpoint. The figure in (e) is intended to show a close approximation to the tile.}
    \end{figure}
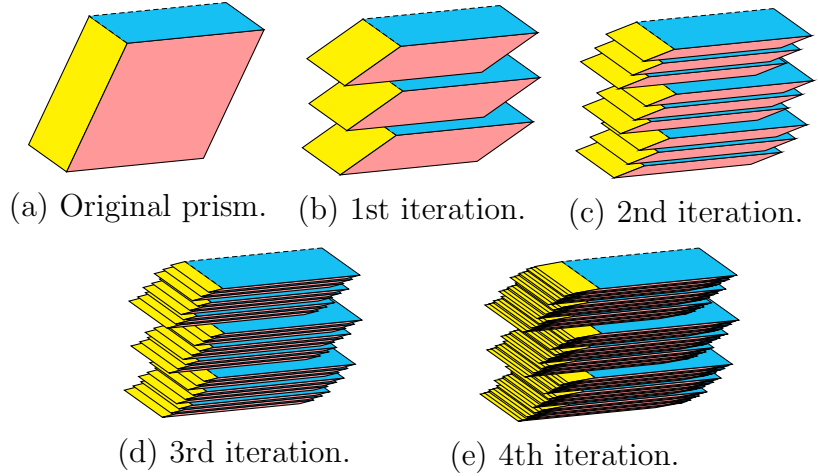

    As the proof of (3) $\Rightarrow (2)$ in Theorem \ref{theomain1} is long and complicated, we outline the main steps here, using the concrete example in Figure~\ref{figillus} as an aid. We begin with a $d$-dimensional slant prism (see Figure \ref{figillus}(a)) that is the union of convex combinations of the top surface and the bottom surface of the self-affine tile $T$ (for the definition of such a convex combination of sets, see \eqref{eqsetconvex}). Then we define five homeomorphisms on $\R^d$ that transform $\R^d$ and the original prism successively. The first is an affine transformation (see Figure \ref{figillus}(a)--(b)) that sets the  slant prism vertical. We then use three transformations: squeeze (see Figure \ref{figillus}(b)--(c)), translation (see Figure \ref{figillus}(c)--(d)), and flattening (see Figure \ref{figillus}(d)--(e)), to turn the $d$-dimensional vertical prism to a union of smaller $d$-dimensional vertical prisms. The fifth homeomorphism is another affine transformation (see Figure \ref{figillus}(e)--(f1)/f(2)) that moves the top surface so that it coincides with the top surface of the original $d$-dimensional slant prism. We call this process an \textit{iteration}. For each of these smaller $d$-dimensional slant prisms, we perform a similar iteration to get several even smaller $d$-dimensional slant prisms, keeping the top and bottom surfaces unchanged. By applying this iteration repeatedly (see Figure \ref{figtile}(a)--(e)), we can define a mapping $h$ on $\R^d$ that sends the original $d$-dimensional slant prism to the self-affine tile $T$. We remark that Figure \ref{figtile}(e) shows a close approximation to the tile. Finally, we use Brouwer's invariance of domain theorem \cite{Brouwer1912} to show that $h$ is a desired homeomorphism  on $\R^d$.
    \begin{theo}[\textrm{Brouwer \cite{Brouwer1912}}]\label{theoregion}
        Let $V\subset\R^d$ be an open set and $h$ a mapping on $V$. If $h$ is injective and continuous, then $h$ is a homeomorphism from $V$ onto $h(V)$.
    \end{theo}
    The continuity and surjectivity of $h$ can be obtained quite easily (Lemma \ref{lemhinftycontinuity}), but it is difficult to prove the injectivity (Lemma \ref{lemh0injective}). To this end, we investigate some analytic properties of each iteration (Proposition \ref{theoh}) and obtain: (1) the height of a point (i.e., the last coordinate in $\R^d$) that is not in the vertical boundary of the original $d$-dimensional slant prism will be fixed after finitely many iterations, and (2) under some mild condition, the absolute value of the height difference of two points that are in the interior of the vertical boundary of the original $d$-dimensional slant prism will not drop. These conclusions, together with a tool which we call \textit{horizontal distance}, imply the injectivity of the mapping $h$.

    The rest of this paper is organized as follows. We introduce some  notation and  concepts in Section \ref{S:notation}, and describe the construction of a quasi-hypersurface in Section 3. In Section \ref{S:section3}, we prove some analytic properties of one iteration, which lead up to Proposition~\ref{theoh}. Section \ref{S:sec4} consists of preparations for constructing the mapping $h$ that sends the self-affine tile to a tame ball. The proof that $h$ is the desired homeomorphism is given in Section \ref{S:sec5}. Finally, in Section \ref{S:questions-comments}, we state some comments on Theorem \ref{theomain2}, which is partly proved in \cite{R3}.

\section{Notation and preliminaries}\label{S:notation}
\subsection{Symbolic space}
    For an integer $m\ge2$, denote $\Sigma_m^k:=\{0,1,\ldots,m-1\}^k$, $\Sigma_m^*:=\cup_{k\ge0}\Sigma_m^k$ and $\Sigma_m^\N:=\{0,1,\ldots,m-1\}^\N$, where $\Sigma_m^0:=\{\emptyset\}$ ($\emptyset$ stands for the {\it empty word}).

    We call $\bi\in\Sigma_m^k$ a \textit{word} of length $k$, and denote its length by $|\bi|$. For $\bi=i_1\cdots i_k\in\Sigma_m^k$ and $\bj=j_1j_2\cdots\in\Sigma_m^*\cup\Sigma_m^\N$, let $\bi\bj=i_1\cdots i_kj_1j_2\cdots$. For $\bi\in\Sigma_m^*$, we denote the infinite word $\bi\bi\cdots$ by $\overline{\bi}$.

    It is well known that the map $\varphi_m:\,\Sigma_m^\N\to [0,1]$ defined as
        \begin{equation*}\label{eqvarph}
        \varphi_m(\bi):=\sum_{n\ge1}\frac{i_n}{m^n},\qquad \bi=i_1 i_2\cdots,
    \end{equation*}
    is surjective. 
    We also define $\varphi_m(\bi):=\sum_{n=1}^{|\bi|}m^{-n} i_n$ for $\bi\in\Sigma_m^*$.


\subsection{Some notation in $\R^d$}
    For $\bx\in\R^d$, we denote the $j$th component of $\bx$, $1\le j\le d$, by $x_j$ or $\omega_j(\bx)$. For $E\subset\R^d$, denote $\omega_j(E)=\{\omega_j(\bx):\,\bx\in E\}$. We let $\bone_d:=(1,\dots,1)\in\R^d$ be the vector whose components are all $1$.

    For $\bx\in\R^d$ and $y\in\R$, let
        \begin{gather*}
            y+\bx:=\big(y+\omega_1(\bx),\ldots,y+\omega_d(\bx)\big),\, y\bx:=\big(y\omega_1(\bx),\ldots,y\omega_d(\bx)\big),\\
           \|\bx\|_\infty:=\max\{|\omega_j(\bx)|:\,1\le j\le d\}.
        \end{gather*}
    We say that $\{\bx_m\}_{m\ge1}$ is $o(1)$, denoted $\bx_m=o(1)$, if
        $\lim_{m\to\infty}\|\bx_m\|_\infty=0.$
    For two vectors $\bx,\by\in\R^d$, we write $\bx\le\by$ or $\by\ge\bx$ if $\omega_j(\bx)\le\omega_j(\by)$ for all $j$ and write $\bx<\by$ or $\by>\bx$ if $\omega_j(\bx)<\omega_j(\by)$ for all $j$.

    Define a new multiplication on $\R^d$ as follows
    \begin{gather}\label{eqmultiply}
        \bx\newdot\by:=(x_1y_1,\ldots,x_dy_d).
        \end{gather}
   We also write $\bx\boxdot\by$ as $\bx\by$ for short. For convenience, we let $\bx\newdot E=\bx E=\{\bx\newdot \by:\,\by \in E\}$ for $\bx\in\R^d$ and $E\subset\R^d$. Recall $\sign(x):=0,1$, or $-1$ according to $x=0,x>0$ or $x<0$. Let $\sign(\bx)$ stand for the vector
        $$\big(\sign(\omega_1(\bx)),\ldots,\sign(\omega_d(\bx))\big),$$
   and let $\bx^{\oplus}:=\sign(\bx)\bx=(|x_1|,\ldots,|x_d|)$. For two sets $E,F\subset\R^d$, define the \textit{sum} of $E$ and $F$ as
        $$E+F:=\{e+f:\,e\in E,f\in F\}.$$
    If $F=\{f\}$, a singleton, then we also write $E+F$ simply as $E+f$. For a set $E\subset\R^d$, let $-E=\{-\bx:\,\bx\in E\}$ and let $E^\circ$ (or $\myint[] E$), $E^c$, and $\pt E$ denote the interior, complement, and boundary of $E$ respectively.

    We point out that the notation and definitions introduced on $\R^d$ ($d\ge 2$) are applicable to sets and vectors in $\R^{d-1}$. However,  if the ambient space $\R^d$ is fixed, the following two concepts are defined only for vectors in $\R^d$.
    For $\bx\in\R^d$ ($d\ge 2$), we call $\cP_h(\bx)=\big(\omega_1(\bx),\ldots, \omega_{d-1}(\bx)\big)$ the \textit{horizontal coordinates} of $\bx$ and $\cP_v(\bx)=\omega_{d}(\bx)$ the \textit{vertical coordinate} of $\bx$. We also call $\cP_v(\bx)$ the \textit{height} of $\bx$.

\subsection{Other notation\label{sectother}}
   Let $r\ge2$ be an integer and $b>0$ be a constant. For an integer $k$ satisfying $0\le k\le r$, denote
        \begin{equation}\label{eqyk}
        y_k:=kb\cdot r^{-1}.
        \end{equation}
    Clearly, $y_0=0$ and $y_r=b$. For $0<k<r$, let
        $\bu_k\in 2U^\circ:=(-1,1)^{d-1}$
    be fixed. In Section 5, we will define $\bu_k$ that have specific geometric meaning in the self-affine tile; we refer the reader to \eqref{eqbunk}. Let $\bv_0$ stand for the origin in $\R^{d-1}$ and let $\bv_k:=\bu_1+\cdots+\bu_k$ for $k\ge1$.

    We choose $\ep\in\R$ small enough such that
        \begin{equation}\label{eqep}
            0<\ep<(2r)^{-1}b.
        \end{equation}
    For $1\le k<r$, let
        $I_{\ep,k}=[y_k-\ep,y_k+\ep]$ and $I_\ep=\bigcup_{k=1}^{r-1}I_{\ep,k}$. Let  $\bx(y)$ and $\brho(y)$ be vector-valued functions defined as
        \begin{eqnarray}\label{eqbxy}
        \bx(y)&:=&
        \begin{cases}
        \bv_{k},&y\in [y_k,y_{k+1}]\backslash I_\ep,0\le k<r;\\
        (2\ep)^{-1}(\ep+y-y_k)\bu_{k}+\bv_{k-1},& y\in I_{\ep,k}, 1\le k<r,
        \end{cases}
        \end{eqnarray}
        \begin{eqnarray}\label{eqrhoy}
        \brho(y)&:=&
        \begin{cases}
        \bone_{d-1},&y\notin I_\ep;\\
        \ep^{-1}|y-y_k|\bu_k^\oplus +1-\bu_k^\oplus , &y\in I_{\ep,k}, 1\le k<r.
        \end{cases}
        \end{eqnarray}
    In particular, $\brho(y_k)=1-\bu_k^\oplus $ and $\bx(y_k)=2^{-1}\bu_k+\bv_{k-1}$. For the geometric significance of $\bx(y)$ and $\brho(y)$, see Figure \ref{figHTS} and Lemma \ref{lemHbi}. Finally, we define a constant $c$ corresponding to $\{\bu_k\}$ as
        \begin{equation}\label{eqc}
        c=c(\{\bu_k\}):=\sum_{k=1}^{r-1}\|\bu_k\|_\infty.
        \end{equation}

\section{Construction of a quasi-hyperplane\label{hyperplane}}
    The structure of the tile allows us to construct a homeomorphism from a $d$-dimensional prism onto a union of some smaller $d$-dimensional prisms (see Figure \ref{figillus}(b)--(e)). However, if we want to extend the domain of such a homeomorphism to the whole space $\R^d$, we need to slice $\R^d$ in such a way that the flattening map (see Figure \ref{figillus}(d)--(e)), which is well-defined on the prism, can be extended to $\R^d$. To achieve this, it is necessary to construct non-trivial quasi-hyperplanes to slice $\R^d$ in a desired way, and this is the main objective of this section.

\subsection{Some notation on $U$.}\label{S:3.1}
    We assume $\bu\in(-1,1)^{d-1}$ and denote $u_j:=\omega_j(\bu)$. We always assume the index $j$ is in $\{1,2,\ldots,d-1\}$ and $\bu\ne(0,\ldots,0)\in\R^{d-1}$ unless stated otherwise.

   Recall that $U=[-1/2,1/2]^{d-1}$ and $\pt U$ denotes its boundary, i.e., $\pt U=\{\bx\in U:\,\omega_j(\bx)=\pm2^{-1} \text{ for some }j\}$.
    With respect to the given vector $\bu$, the boundary of $U$ can be divided into three parts, not necessarily disjoint, as follows:
        \begin{equation}\label{eqptU}
        \pt U=\pt U^+\cup \pt U^-\cup \pt U^{\updownarrow},
        \end{equation}
    where
        \begin{gather*}
        \begin{aligned}
        &\pt U^+=\pt U^+_\bu:=\{\bx\in\pt U:\,2\omega_j(\bx)=\sign(\omega_j(\bu))\ne0\mbox{ for some } j\},\\
        &\pt U^-=\pt U^-_\bu:=-\pt U^+,\\
        &\pt U^{\updownarrow}=\pt U^{\updownarrow}_\bu:=\big\{\bx\in\pt U:\,2 \omega_j(\bx)=\pm1\text{ and }\omega_j(\bu)=0 \text{ for some } j\big\},
        \end{aligned}
        \end{gather*}
    The set $\pt U^{\updownarrow}$ is empty if none of the components of $\bu$ is zero. To explain our idea more clearly, we present the following example.
    \begin{exam} Let $I:=[-1/2,1/2]$.
    \begin{enumerate}
    \item[(1)]  If $d=4$ and $\bu=(-0.1,0.2,0)$, then \textup{(}see Figure \ref{figptU}\textup{(a,b,c))}
        \begin{gather*}
    \begin{aligned}
        &\pt U^+=\bigg(\Big\{-\dfrac12\Big\}\times I^2 \bigg)\bigcup \bigg(I\times\Big\{\dfrac12\Big\} \times I\bigg),\\
        &\pt U^-=\bigg(\Big\{\dfrac12\Big\}\times I^2 \bigg) \bigcup \bigg(I\times\Big\{-\dfrac12\Big\} \times I\bigg),\\
        &\pt U^{\updownarrow} =I^2\times\Big\{-\dfrac12,\dfrac12\Big\}.
    \end{aligned}
        \end{gather*}
    \item[(2)] If $d=4$, $u=(0.1,-0.2,-0.3)$, then $\pt U^{\updownarrow}=\emptyset$ and \textup{(}see Figure \ref{figptU}\textup{(d,e))}
            \begin{gather*}
        \begin{aligned}
            &\pt U^+=\Big(\Big\{\dfrac12\Big\}\times I^2\Big)
            \bigcup \Big(I\times \Big\{-\dfrac12\Big\}\times I\Big) \bigcup\Big(I^2\times\Big\{-\dfrac12\Big\}\Big), \\
            &\pt U^-=\Big(\Big\{-\dfrac12\Big\}\times I^2\Big)
            \bigcup \Big(I\times \Big\{\dfrac12\Big\}\times I\Big) \bigcup\Big(I^2\times\Big\{\dfrac12\Big\}\Big).
        \end{aligned}
    \end{gather*}
    \end{enumerate}
    \end{exam}
\def\ed{\end{document}}

    \begin{figure}[htbp]
    \begin{tikzpicture}[scale=1.35,>=stealth,x={(-.5cm,-.3cm)}, y={(1cm,0cm)},z={(0cm,1cm)}]
    \pgfmathsetmacro{\h}{2}
    \pgfmathsetmacro{\v}{-2.1}
    \coordinate (a) at (.5,-.5,-.5); 
    \coordinate (b) at (-.5,-.5,-.5); 
    \foreach \ad/\adtex in {%
        0/{\sign(\bu)=(-1,1,0)}, 1/{\sign(\bu)=(1,-1,-1)}}{
        \foreach \i/\col/\k/\tex/\texx/\text in
            {0/red!50/{(2*\ad-1)}/\pt U^+/(a)/(d), 1/blue!50/{(1-2*\ad)}/\pt U^-/(b)/(e)}{
            \coordinate (a) at (0,{\i*\h+\ad*\h/2},\ad*\v);
            \ifthenelse{\ad=1}{
                \draw($(a)+(0,0,-.45)$)node{\text\,$\tex$};
                    }
                {\draw($(a)+(0,0,-.45)$)node{\texx\,$\tex$};}
            \myycub{0,{\i*\h+\ad*\h/2},{\ad*\v}}[\col][{\k}]
            }
        \path[fill=gray,opacity=.5](-.5,2*\h-.5,\ad)--++(1,0,0)--++(0,1,0) --++(-1,0,0)--cycle;
        \foreach \i/\tex/\sty in {{.5,0,0}/x_1/left,{0,.5,0}/x_2/right, {0,0,.5}/x_3/left}
            \draw[->](0,-\h,.1+\ad*\v)--++(\i)node[\sty]{$\tex$};
        \draw(0,-\h,\ad*\v-.5)node[font=\tiny]{$\adtex$};
        }
    \myycub{0,2*\h,0}[white][0];
    \draw(0,2*\h,-.45)node{(c)\,$\pt U^{\updownarrow}$};
    \end{tikzpicture}
    \caption{\label{figptU} Illustration of the decomposition of $\pt U$ with respect to the vector $\bu$. The surfaces colored {red, blue, and black} are the sets $\pt U^+,\pt U^-$, and $\pt U^{\updownarrow}$ respectively.}
    \end{figure}
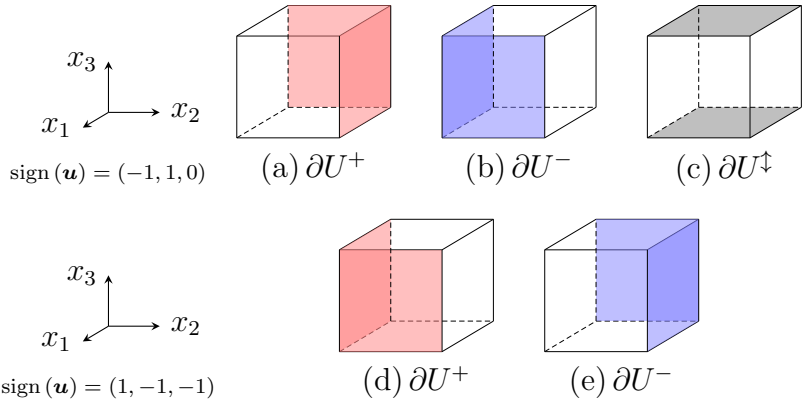

\subsection{Construction of quasi-hyperplane.}
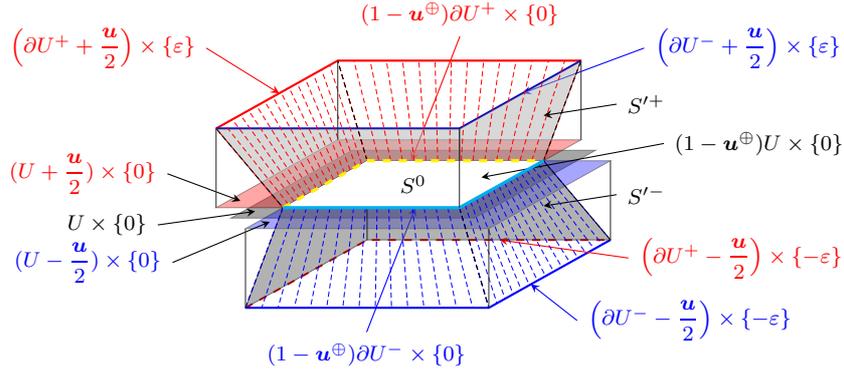
\begin{figure}
        \begin{tikzpicture}[{font=\tiny,>=stealth,xscale=3.2,yscale=3, declare function={f(\t)=(1-\rho)*\t/\h+\rho;}, x={(-.5cm,-.3cm)}, y={(1cm,0cm)},z={(0cm,1cm)}}]
        \mystart[2][.35]
        \coordinate (o) at (0,0,0);
        \coordinate (a) at (-.5*\hx,-.5*\hy,0);  
        \coordinate (c) at ($(a)+(\hx,\hy,0)$);  
        \foreach \i/\j/\tex/\col in {
            o/{.2,-.3,0}/U\times\{0\}/black,
            {-\ux,-\uy,0}/{-.4,-.5,0}/(U+\dfrac\bu2)\times\{0\}/red,
            {\ux,\uy,0}/{.6,-.1,0}/(U-\dfrac\bu2)\times\{0\}/blue}
            {\draw[fill=\col,opacity=.35]($.5*(\i)-.5*(1,1,0)$)--++(1,0,0) --++(0,1,0) --++(-1,0,0)--cycle;
            \draw[<-]($.5*(\i)+(.5-.1,.05-.5,0)$)--++(\j) node[\col,left,font=\tiny]{$\tex$};}

        \rectline{-.5*\ux,-.5*\uy,\eps}11{red,thick}{blue,thick}
        \rectline{.5*\ux,.5*\uy,-\eps}11{red,dashed,thick}{blue,thick}        
        \foreach \i/\j/\tex/\col/\sty in {
            {-.25*\hx,-.5*\uy-.5,\eps}/{0,-.3,.2}/ {\Big(\pt U^++\dfrac \bu2\Big)\times\{\ep\}}/red/left,
            {-.25*\hx,.5*\hy,\eps}/{0,.5,.2}/{\Big(\pt U^-+\dfrac \bu2\Big)\times\{\ep\}}/blue/right,
            {-\hx*.5,-.25*\hy,0}/{0,.2,\eps+.2}/{(1-\bu^\oplus )\pt U^+\times\{0\}}/red/above,
            {\hx*.5,.25*\hy,0}/{0,-.2,-\eps-.2}/{(1-\bu^\oplus )\pt U^-\times\{0\}}/blue/below,
            {-.5*\hx,.25*\hy,-\eps}/{.1,\hy*.8,-.05}/{\Big(\pt U^+-\dfrac\bu2\Big) \times\{-\ep\}}/red/right,
            {\ux,\hy*.5+\uy,-\eps}/{.2,\hy*.4,-.1}/{\Big(\pt U^--\dfrac\bu2\Big)\times\{-\ep\}}/blue/right}
        \draw[<-,\col](\i)--++(\j)node[font=\tiny,\sty]{$\tex$};

        \foreach \drt/\tex/\col/\colo in {1/S'^+/black!30/red,-1/S'^-/black!60/blue}
            {\foreach \i/\j/\k in
                {{\hx,0,0}/{1,0,0}/{-\ux,0,\eps}, {0,\hy,0}/{0,1,0}/{0,-\uy,\eps}}
            \draw[fill=\col,opacity=.5,very thin] ($(a)-{(\drt+1)}*(a)$) --++($-\drt*(\i)$)--++($\drt*(\k)$)--++($\drt*(\j)$)--cycle;
            \draw[<-]($.5*(-\hx,\hy,\drt*.4)$)--++(0,.3,.05) node[right]{$\tex$};
            \foreach \ii in {0,1,...,15}
                \foreach \jj/\len in {{1,0,0}/\hx,{0,1,0}/\hy}{
                    \coordinate (temp1) at
                        ($(a)+{(1-\drt)/2}*(\hx,\hy,0)$);
                    \coordinate (temp2) at
                        ($(temp1)-\drt*(\ux,\uy,-\eps)$);
                    \draw[dot2,\colo]($(temp1)+\drt*\ii/15*\len*(\jj)$)                       --($(temp2)+\drt*\ii/15*(\jj)$);}
            \draw[dot2]($(a)+.5*(\drt+1,1-\drt,-2*\eps)$) --++({-.5*\ux*(\drt+1)},{-.5*\uy*(1-\drt)},\eps) --++({-.5*\ux*(1-\drt)},{-.5*\uy*(\drt+1)},\eps);
            \draw[dot2]($-\drt/2*(\hx,\hy,0)$) --++($\drt*(-\ux,-\uy,\eps)$);
        }

        \rectreg{o}{\hx}{\hy}{fill=white}
        \draw[<-](0,0,0)node{$S^0$}(-.4,.08,-.09)--++(0,.75,.15) node[font=\tiny,right]{$(1-\bu^\oplus)U\times\{0\}$};
        \rectline{o}\hx\hy{very thick,yellow,dashed}{very thick,cyan}
        \foreach \i in {{0,-1,0},{-1,-1,0},{-1,0,0},{0,0,0}}
            \draw[black!80]($(c)+(\i)$)--++(0,0,\eps);
        \foreach \i in {{0,1,0},{1,1,0},{1,0,0}}
            \draw[black!80]($(a)+(\i)$)--++(0,0,-\eps);
        \end{tikzpicture}
        \caption{\label{figS} An illustration of $S^0,S^+$, $S^-$, $S'^+$ and $S'^-$. This figure is drawn with $d=3$ and $\bu=(-0.3,-0.25)$.}
    \end{figure}

    We continue to use the assumptions and notation in Section~\ref{S:3.1}. Recall that $\bx E=\{\bx\by:\,\by\in E\}$ for $\bx\in\R^n$ and $E\subset\R^n$, where $n=d-1$ or $n=d$.
    Let $\ep>0$ be small enough. Denote $S^0:=\big((1-\bu^\oplus )U\big)\times\{0\}$. Let
        $$\hat\bu_t:=(1-t)(1-\bu^\oplus )+t\times\bone_{d-1}=t\bu^\oplus +1-\bu^\oplus,\qquad0\le t\le 1.$$
    $\hat\bu_t$ is the convex combination of the vectors $1-\bu^\oplus $ and $\bone_{d-1}$ with parameter $t$; it specifies the amount of contraction in each of the $d-1$ directions, in the range determined by  $\bu$ and $\bone_{d-1}$. Define
        \begin{equation}\label{eqfourS}
        \begin{array}{rcl}
            S^+=S^+_{\bu}&:=&\dis\bigcup_{0\le t\le1}\Big(\hat\bu_t\pt U^+ +\dfrac t2\bu\Big)\times\{t\ep\},\\
            S'^+=S'^+_{\bu}&:=&\dis \bigcup_{0\le t\le1}\Big(\hat\bu_t(\pt U^-\cup \pt U^\updownarrow\big)+\dfrac t2\bu\Big)\times\{t\ep\}.
        \end{array}
        \end{equation}
    We can use the notion of convex combination to describe these sets. For $E\subset \R^{d-1}$, $\bx_1,\bx_2,\bt_1,\bt_2\in\R^{d-1}$, and $y_1,y_2\in\R$, the \textit{convex combination} of two sets $(\bt_1E+\bx_1)\times\{y_1\}$ and $(\bt_2E+\bx_2)\times\{y_2\}$ with respect to the combination coefficient $t$ is defined to be the set
            \begin{equation}\label{eqsetconvex}
            \big\{\big((1-t)(\bt_1\bx+\bx_1)+t(\bt_2\bx+\bx_2), (1-t)y_1+ty_2\big): \,\bx\in E\big\}.
            \end{equation}
    Now $S^+$ can be described as the union of the convex combinations of
    $$\big((1-\bu^\oplus )\pt U^+\big)\times\{0\}\quad\text{and}\quad\Big(\pt U^++\dfrac\bu2\Big)\times\{\ep\},
    $$
    while $S'^+$ is that of
    $$\big((1-\bu^\oplus )(\pt U^-\cup\pt U^{\updownarrow})\big)\times\{0\}\quad \text{and}\quad\Big(\pt U^-\cup\pt U^{\updownarrow}+\dfrac\bu2\Big)\times\{\ep\}.$$
    See Figure \ref{figS} for an illustration. We also point out that if we let
   $$\pi(y):=\R^{d-1}\times\{y\},\quad y\in\R,$$
    then $S'^+$ is perpendicular to the hyperplane $\pi(0)$, i.e., any normal vector of $S'^+$ is parallel to $\pi(0)$. Let $S^-:=-S^+$ and $S'^-:=-S'^+$.

\begin{figure}
    \begin{tikzpicture}[scale=3.3,>=stealth, font=\tiny]
        \pgfmathsetmacro{\x}{-.5}
        \pgfmathsetmacro{\y}{.6}
        \pgfmathsetmacro{\t}{.3}
        \pgfmathsetmacro{\ux}{(1-\t)*(1-abs(\x))+\t};
        \pgfmathsetmacro{\uy}{(1-\t)*(1-abs(\y))+\t};
        \coordinate (x1) at (.48,-.55);
        \coordinate (x2) at (.6,.5);
        \coordinate (x3) at (-.8,.2);
        \coordinate (x4) at (-.6,-.6);
        \coordinate(xp4)at($(-2,-.5*\uy)!(x4)!(1,-.5*\uy)$);
        \coordinate(xn4)at($(-.5*\ux,-2)!(x4)!(-.5*\ux,1)$);

        \draw(x3)--node[below,]{$d_t^+(\bx_3)$} ($(-.5*\ux,-2)!(x3)!(-.5*\ux,1)$);
        \draw(x1)--node[right,font=\tiny]{$d_t^{-}(\bx_1)$} ($(-2,-.5*\uy)!(x1)!(1,-.5*\uy)$);

        \draw(x2)--node[right]{$d_t^+(\bx_2)$}
        ($(-2,.5*\uy)!(x2)!(1,.5*\uy)$);
        \draw(x2)--node[above]{$d_t^{-}(\bx_2)$}
        ($(.5*\ux,-1)!(x2)!(.5*\ux,1)$);
        \draw[red,very thick](-.5*\ux,-.5*\uy)|-++(\ux,\uy);
        \draw[blue,dashed,very thick](-.5*\ux,-.5*\uy)-|++(\ux,\uy);
        \draw(0,0)node[font=\small]{$\hat\bu_t U$};
        \foreach \i/\tex/\loc in {
            x1/\bx_1/below,x2/\bx_2/below left,x3/\bx_3/left,x4/\bx_4/above right}
            \draw[fill=black](\i)circle(.01)node[\loc]{$\tex$};
        \foreach \i in {-.5,.5}
            {\draw[help lines](-1,\i*\uy)--++(2,0);
            \draw[help lines](\i*\ux,-.7)--++(0,1.5);}
        \foreach \i/\col/\tex/\sty in {1/red/+/above,-1/blue/-/below}
            \draw[\col,<-](0,\i*\uy*.5)--++(-\i*.1,\i*.3) node[\sty,font=\small]{$\hat\bu_t\pt U^{\tex}$};

        \draw(x4)--node[left]{$d_t^{-}(\bx_4)$}(xp4);
        \draw(x4)--node[below]{$d_t^+(\bx_4)$}(xn4);
        \draw(1,0)node[right, font=\small, text width=100pt]{The sign of $\bu$ is $$\sign(\bu)=(-1,1).$$};
    \end{tikzpicture}
    \caption{\label{figdt} Illustration of $d_t^+$ and $d_t^{-}$. The union of the thick lines (in red) stands for $\hat\bu_t\pt U^+$, while the union of the dashed lines (in blue) stands for $\hat\bu_t \pt U^{-}$. This figure is drawn with $\bu=(-0.5,0.6)$ and $t=0.3$.}
    \end{figure}
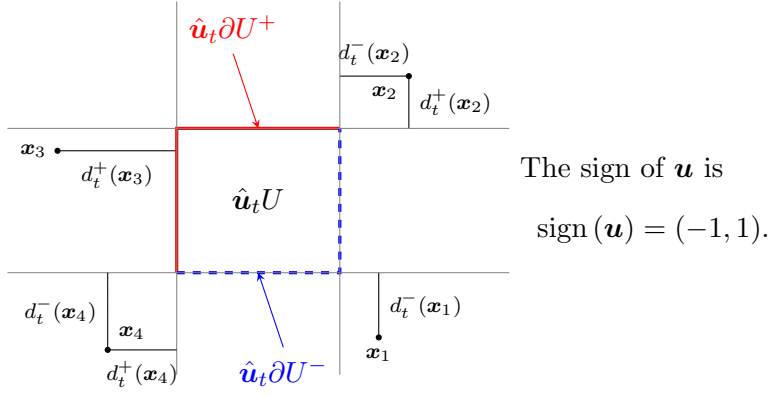

    Next we will extend the surface $S^0\cup S^+\cup S^-$ to a $(d-1)$-dimensional quasi-hyperplane $S_\bu$. As the extension is nontrivial and needs some fine structure, we introduce some additional notation. For $\bx\in\R^{d-1}$, let
        \begin{equation}\label{eqellbx}
        \ell_\bx:=\{\bx\}\times\R.
        \end{equation}
    Next, we introduce some concepts on the difference of components with respect to the given vector $\bu$, which will make the extension easier. For $0\le t\le1$, define two pseudo-distances between a point $\bx\notin\hat\bu_tU^\circ$ and the set $\hat\bu_t U$ as follows:
    \begin{gather*}
        d_t^+(\bx):=\max\Big\{|x_j|-\dfrac12\omega_j(\hat\bu_t):\, x_ju_j>0,|x_j|\ge\dfrac12\omega_j(\hat\bu_t), j=1,\dots, d-1\Big\},\\
        d_t^{-}(\bx) :=\max\Big\{|x_j|-\dfrac12\omega_j(\hat\bu_t):\, x_ju_j\le0,|x_j|\ge\dfrac12\omega_j(\hat\bu_t),j=1,\dots,d-1\Big\},
        \end{gather*}
    where the maximum of an empty set is defined to be $-\infty$. We will explain the geometric meanings of $d_t^+$ and $d_t^+$ following Lemma~\ref{lemex}, or more precisely, following Equation \eqref{eqdt+le}. We now present an example to illustrate these definitions.
    \begin{exam} Assume $d=3$, $\bu=(-0.5,0.6)$ and $t=0.3$.
    Then
        $$\hat\bu_t=(1-0.3)\big(1-(0.5,0.6)\big)+0.3(1,1)=(0.65,0.58).$$
    For
        $\bx_1=(0.48,-0.55), \bx_2=(0.6,0.5), \bx_3=(-0.8,0.2), \bx_4=(-0.6$, $-0.6),$
    we have $($see Figure \ref{figdt}$)$
        \begin{gather*}
        \begin{aligned}
        &d_t^+(\bx_1)=-\infty,\quad d_t^{-}(\bx_1)=|\omega_2(\bx_1)|-\dfrac12\omega_2(\hat\bu_t) =0.55-0.29=0.26,\\
        &d_t^+(\bx_2)=\omega_2(\bx_2)\!-\!\dfrac12\omega_2(\hat\bu_t) =0.21,\quad
        d_t^{-}(\bx_2)\!=\!|\omega_1(\bx_2)|\!-\! \dfrac12\omega_1(\hat\bu_t)\!=\!0.275,\\
        &d_t^+(\bx_3)=|\omega_1(\bx_3)|-\dfrac12\omega_1(\hat\bu_t)=0.475,\quad d_t^{-}(\bx_3)=-\infty,\\
        &d_t^+(\bx_4)=0.275,\quad d_t^{-}(\bx_4)=0.31.
        \end{aligned}
        \end{gather*}
    \end{exam}

    Recall that for two vector $\bx,\by\in\R^{d-1}$, $\bx\le\by$ (resp. $\bx<\by$) means
        $$\omega_j(\bx)\le\omega_j(\by)\quad\big(\text{resp. } \omega_j(\bx)<\omega_j(\by)\big),\quad 1\le j\le d-1.$$
    Also, $U^\circ=U\backslash \pt U=\{\bx\in U:\,-2^{-1}\bone_{d-1}<\bx<2^{-1}\bone\}$. According to the size of the quantities $d_t^+$ and $d_t^{-}$, we define the following subset of $\R^{d-1}$: 
    \begin{equation}\label{eqGt+-}
        G_t=\big\{\bx\in\R^{d-1}:\,\hat\bu_t\le2\bx^\oplus\le 2\eta\sqrt{1-t}+1+(1-t)\bu^\oplus,d_t^+(\bx)\ge d_t^-(\bx)\big\},
    \end{equation}
    where $\eta\in(\|\bu\|_\infty,1)$ is some constant. We note that $G_1=\pt G_1=\pt U^+$. See Figure \ref{fig2-2D} for an illustration of $\eta\sqrt{1-t}$ and $\hat\bu_t$ and see Figures \ref{figptG}, \ref{figptG2} and \ref{figS''} for an illustration of $G_{t}$.

\begin{figure}
    \begin{tikzpicture}[xscale=4.5,yscale=2.8,>=stealth]
        \pgfmathsetmacro{\x}{.4}
        \pgfmathsetmacro{\del}{.45}
        \pgfmathsetmacro{\v}{1}
        \pgfmathsetmacro{\t}{.40}
        \pgfmathsetmacro{\tt}{sqrt(1-\t)}
        \pgfmathsetmacro{\plen}{1.25}
        \coordinate (a) at (-.5-\x*.5,-\v);
        \coordinate (b) at (.5+\x*.5,\v);
        \coordinate (l) at (-\del-.5-\x*.5,0);
        \coordinate (r) at (\del+.5+\x*.5,0);
        \draw[very thick]($(l)+(-.15,0)$)--++(.15,0)parabola bend (a) (a)--++(\x,\v)--++(1-\x,0)--(b)parabola bend (b) (r)--++(\plen-1,0);
        \mybrace[1]{l}{\del}{$\eta$};
        \mybrace[1]{{-.5,0}}{.5*\x}{$\dfrac{\bu^\oplus }2$};
        \mybrace[1]{{-.5-\x*.5,0}}{.5*\x}{$\dfrac{\bu^\oplus }2$};
        \mybrace{.5-\x/2,0}{.5*\x}{$\dfrac{\bu^\oplus }2$};
        \mybrace{.5+\x*.5,0}{\del}{$\eta$};
        \mybrace{\tt*\del+.5+.5*\x,\t}{-\tt*\del}{$\eta\sqrt{1-t}$};
        \mybrace{.5-\x*.5+\x*\t,\t}{(-1-\x*\t+\x)}{$\hat\bu_t$};
        \mybrace[1]{.5-\x*.5+\x*\t,\t}{-\x*\t}{$t\bu^\oplus $};
        \mybrace{.5+\x*.5,\t}{(\t-1)*\x}{
        $(1\!-\!t)\bu^\oplus $};
        \mybrace{{.5,0}}{\x*.5}{$\dfrac{\bu^\oplus }2$};
        \mybrace{{-.5-\x*.5,-1}}{1}{$\hat\bu_1$};
        \mybrace[1]{{-.5+\x*.5,0}}{(1-\x)}{$\qquad\hat\bu_0$};
        \draw[fill=black](0,0)circle(.02)node[font=\small,below left]{$O$};
        \foreach \i/\tex/\loc in {\t/t\ep/above left,1/\ep/below left}
            \draw[fill=black](0,\i)circle(.02)node[\loc]{$\tex$};

       \foreach \i in {a,{\x*.5-.5,0},{.5-\x*.5,-\v},{.5+\x*.5,0}}
            \draw(\i)--++(0,\v);
        \foreach \i/\j/\sty in {-1/-/below,1/{}/below}
            \draw(\i*.5,0)node[\sty,yshift=-\i*2pt,font=\tiny]{$\j\dfrac12$}--++(0,.1);
        \mybrace[1]{{-.5+\x*.5,1}}1{$\hat\bu_1$};

        \draw[->](0,0)--++(\t*\x*.5+.08+\plen,0)node[above]{$x$};
        \draw[->](0,-\v-.05)--++(0,2*\v+.35)node[left]{$y$};
        \def\col{red}
        \mybrace[1]{{0,\t}}{\t*\x*.5}{};
        \def\mybraceht{.24}
        \draw[->,red](\x*\t*.25,\t*1.2)--++(-.5,.3) node[font=\tiny,xshift=5pt,yshift=-5pt,above left]{$\dfrac{t}2\bu^\oplus $};
        \def\col{black}
        \mybrace[1]{{\t*\x*.5,\t}}{\plen}{$\bx^{\oplus}$};
        \draw[fill=black](\plen+\t*\x*.5,\t)circle(.02)node[right]{$P$};
        \foreach \i/\tex/\loc/\j in {{.5-\x*.5,0}/A_1/above/-3pt, {.5+.5*\x,1}/P_1/above/0pt, {.5+\x*.5+\del,0}/B_1/below/0pt, {\x*.5-.5,0}/A_3/below/0pt, {-.5-.5*\x,-1}/Q_3/below/0pt, {-.5-\x*.5-\del,0}/B_3/above/-3pt}
            \draw[fill=black](\i)circle(.02) node[xshift=\j,\loc,font=\tiny]{$\tex$};
    \end{tikzpicture}
    \caption{\label{fig2-2D} An illustration of $\eta\sqrt{1-t}$ and $\hat\bu_t$ for some fixed $\bx\in\mathbb R^{d-1}$. The thick piecewise curve stands for the quasi-hyperplane $S_\bu$. The line segment $A_1P_1$, the line segment $A_3Q_3$,  the curve $A_1P_1B_1$, and the curve $A_3Q_3B_3$ represent $S^+$, $S^-$, $S''^+$, and $S''^-$, respectively. The horizontal coordinate of $P$ is $\bx+\frac t2\bu$, which is used in the proof of Lemma \ref{lemex}.}
\end{figure}
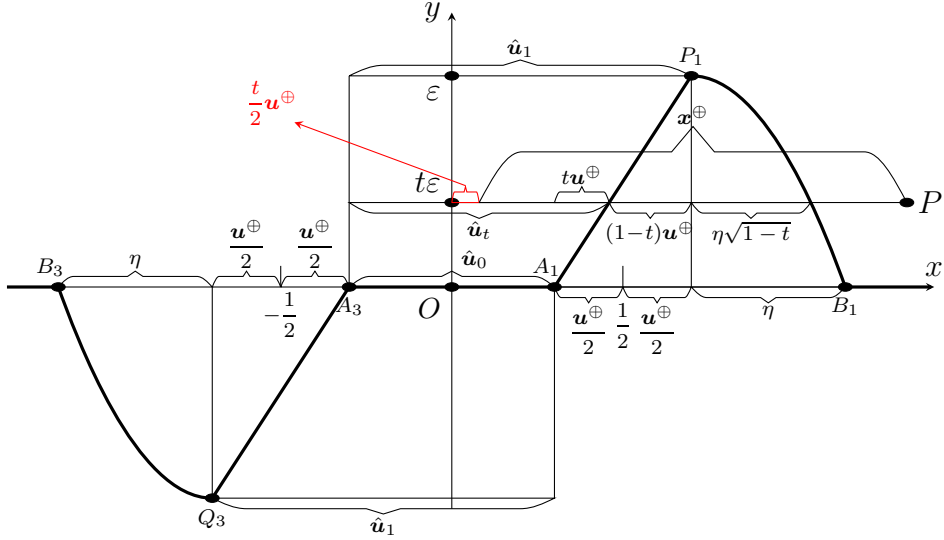
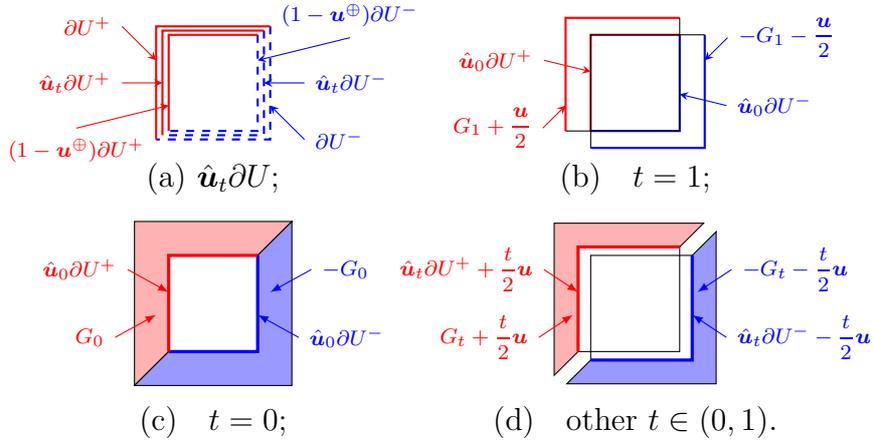
\begin{figure}
    \begin{tikzpicture}[xscale=1.5,yscale=1.5,>=latex,font=\tiny]
        \pgfmathsetmacro{\h}{-3.7} 
        \pgfmathsetmacro{\hh}{-3.7} 
        \pgfmathsetmacro{\t}{.5}   
        \pgfmathsetmacro{\del}{.3}
        \pgfmathsetmacro{\dt}{.7*\del}
        \pgfmathsetmacro{\v}{1.95} 
        \pgfmathsetmacro{\vx}{-.22} 
        \pgfmathsetmacro{\vy}{.15} 
        \pgfmathsetmacro{\bux}{(1+\vx)} 
        \pgfmathsetmacro{\buy}{(1-\vy)} 

        \pgfmathsetmacro{\x}{\vx*\t}
        \pgfmathsetmacro{\y}{\vy*\t}
        \pgfmathsetmacro{\bx}{(1+\x)}
        \pgfmathsetmacro{\by}{(1-\y)}

        \coordinate (vbr) at (.5*\bux,-.5*\buy); 
        \coordinate (val) at (-.5*\bux,.5*\buy); 

        \coordinate (br) at (.5*\bx,-.5*\by); 
        \coordinate (al) at (-.5*\bx,.5*\by); 
        \coordinate (a) at (-.5*\bux,-.5*\buy); 
        \coordinate (a1) at ($(a)+(\bux-\bx,0)$); 
        \coordinate (a2) at ($(a)+(0,\buy-\by)$); 

        \foreach \i/\j/\k in {{.5,-.5}/1/1,vbr/\bux/\buy,br/\bx/\by}
            \draw[blue,dashed,thick]($(\i)+(\h,\k+\v)$)|-++(-\j,-\k);
        \foreach \i/\j/\k in {{-.5,.5}/1/1,val/\bux/\buy,al/\bx/\by}
            \draw[red,thick]($(\i)+(\h,\v-\k)$)|-++(\j,\k);
        \foreach \j/\k in {1/1,\bux/\buy}
            \draw[red,thick]($(vbr)+(-\j,\v)$)|-++(\j,\k);
        \foreach \j/\k in {1/1,\bux/\buy}
            \draw[blue,thick]($(val)+(0,\v-\k)$)-|++(\j,\k);
        \draw($(val)+(0,\v-1)$)|-++(1,1);
        \draw($(vbr)+(-1,\v)$)-|++(1,1);

        \foreach \i/\j/\tex/\loc in
            {{-\bux*.5,-.2}/{-1.2+\bux*.5,-.2}/(1-\bu^\oplus )\pt U^+/below,
            {-.5*\bx,0}/{-.8+.5*\bx,0}/\hat\bu_t\pt U^+/left, {-.5,.2}/{-.3,.3}/\pt U^+/left}
            \draw[red,<-,>=stealth]($(\i)+(\h,\v)$) --++(\j)node[\loc]{$\tex$};
        \foreach \i/\j/\tex/\loc in
            {{\bux*.5,.2}/{1.2-\bux*.5,.2}/(1-\bu^\oplus )\pt U^-/above, {.5*\bx,0}/{.8-.5*\bx,0}/\hat\bu_t\pt U^-/right, {.5,-.2}/{.3,-.3}/\pt U^-/right}
            \draw[blue,<-,>=stealth]($(\i)+(\h,\v)$) --++(\j)node[\loc]{$\tex$};

        \foreach \i/\j/\tex/\loc in
            {
            {-\bux*.5,.1}/{-.3-.5+\bux*.5,.1}/\hat\bu_0\pt U^+/left, 
            {-.5+.5*\vx,-.2}/{-.2,-.2}/G_1+\dfrac\bu2/left}
            \draw[red,<-,>=stealth]($(\i)+(0,\v)$) --++(\j)node[\loc]{$\tex$};

        \foreach \i/\j/\tex/\loc in
            {{\bux*.5,-.1}/{.8-\bux*.5,-.1}/\hat\bu_0\pt U^-/right, {.5-.5*\vx,.2}/{.2,.2}/-G_1-\dfrac\bu2/right}
            \draw[blue,<-,>=stealth]($(\i)+(0,\v)$) --++(\j)node[\loc]{$\tex$};

        \draw[fill=red!30]($(a)+(\hh,0)$)--++(-\del,-\del)--++(0,2*\del+\buy) --++(2*\del+\bux,0)--++(-\del,-\del)--++(-\bux,0)--cycle;
        \draw[fill=blue!40]($(a)+(\hh,0)$)--++(-\del,-\del)--++(2*\del+\bux,0) --++(0,2*\del+\buy)--++(-\del,-\del)--++(0,-\buy)--cycle;
        \foreach \i/\col in {{\bux,0}/blue,{0,\buy}/red}
            \draw[\col,very thick]($(a)+(\hh,0)$)--++(\i)--($(a)+(\hh+\bux,\buy)$);
        \draw[fill=red!30](a1)--++(-\dt,-\dt)--++(0,2*\dt+\by) --++(2*\dt+\bx,0)--++(-\dt,-\dt)--++(-\bx,0)--cycle;
        \draw[fill=blue!40](a2)--++(-\dt,-\dt)--++(2*\dt+\bx,0) --++(0,2*\dt+\by)--++(-\dt,-\dt)--++(0,-\by)--cycle;
        \foreach \i/\j in {vbr/{-\bx,\by},val/{\bx,-\by}}
            \draw(\i)rectangle++(\j);
        \draw[very thick,red](a1)|-++(\bx,\by);
        \draw[very thick,blue](a2)-|++(\bx,\by);

        \foreach \i/\j\tex/\col/\loc in {{-\bx*.5-.1,-.1}/{-\del,-.2}/G_t+\dfrac t2\bu/red/left, {\bx*.5+.1,.1}/{\del,.2}/-G_t-\dfrac t2\bu/blue/right,
        {-\bx*.5+\x*\t,.1}/{-\del,.2}/\hat\bu_t\pt U^++\dfrac t2\bu /red/left,
        {\bx*.5-\x*\t,-.1}/{\del,-.2}/\hat\bu_t\pt U^--\dfrac t2\bu/blue/right}
            \draw[\col,<-](\i)--++(\j)node[\loc]{$\tex$};

        \foreach \i/\j\tex/\col/\loc in {{-.5,-.1}/{-1.2*\del,-.2}/G_0/red/left, {.5,.1}/{1.2*\del,.2}/-G_0/blue/right,
        {-\bux*.5,.1}/{-1.2*\del,.2}/\hat\bu_0\pt U^+/red/left, {\bux*.5,-.1}/{1.2*\del,-.2}/\hat\bu_0\pt U^-/blue/right}
            \draw[\col,<-]($(\i)+(\hh,0)$)--++(\j)node[\loc]{$\tex$};

        \foreach \i/\tex in {{\hh,\v-.85}/{(a) $\hat\bu_t\pt U$;}, {0,\v-.85}/{(b)\quad $t=1$;}, {\hh,-1.05}/{(c)\quad $t=0$;}, {0,-1.05}/{(d)\quad other $t\in(0,1)$.}}
            \draw(\i)node[font=\normalsize]{\tex};
    \end{tikzpicture}
    \caption{\label{figptG} An illustration of $\pt U^+$, $\pt U^-$ and $G_t,-G_t$ with none of the components of $\bu$ being zero. This figure is drawn with $\bu=(-0.22,0.15)$.}
\end{figure}

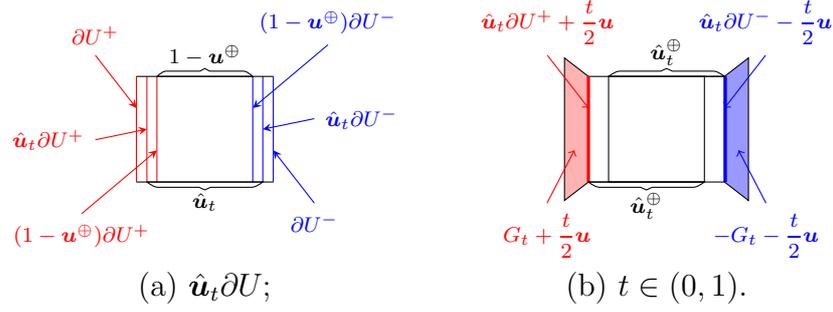
\begin{figure}
    \begin{tikzpicture}[xscale=1.8,yscale=1.4,font=\tiny]
        \pgfmathsetmacro{\h}{-3.3} 
        \pgfmathsetmacro{\t}{.51}   
        \pgfmathsetmacro{\del}{.25}
        \pgfmathsetmacro{\dt}{.7*\del}
        \pgfmathsetmacro{\vx}{-.3} 
        \pgfmathsetmacro{\bux}{(1+\vx)} 

        \pgfmathsetmacro{\x}{\vx*\t}
        \pgfmathsetmacro{\bx}{(1+\x)}

        \coordinate (vbr) at (.5*\bux,-.5); 
        \coordinate (val) at (-.5*\bux,.5); 

        \coordinate (br) at (.5*\bx,-.5); 
        \coordinate (al) at (-.5*\bx,.5); 
        \coordinate (a) at (-.5*\bux,-.5); 
        \coordinate (a1) at ($(a)+(\bux-\bx,0)$); 
        \coordinate (a2) at (a); 

        \foreach \i in {{.5,-.5},vbr,br}
            \draw[blue]($(\i)+(\h,0)$)--++(0,1);
        \foreach \i/\j/\k in {{-.5,.5},val,al}
            \draw[red]($(\i)+(\h,0)$)--++(0,-1);
        \foreach \i in {-.5,.5}\draw(-.5+\h,\i)--++(1,0);

        \foreach \i/\j/\tex/\loc in
            {{-\bux*.5,-.2}/{-.9+\bux*.5,-.6}/(1-\bu^\oplus )\pt U^+/below,
            {-.5*\bx,0}/{-.8+.5*\bx,-.1}/\hat\bu_t\pt U^+/left, {-.5,.2}/{-.3,.5}/\pt U^+/above}
            \draw[red,<-,>=stealth]($(\i)+(\h,0)$) --++(\j)node[\loc]{$\tex$};
        \foreach \i/\j/\tex/\loc in
            {{\bux*.5,.2}/{.9-\bux*.5,.6}/(1-\bu^\oplus )\pt U^-/above, {.5*\bx,0}/{.8-.5*\bx,.1}/\hat\bu_t\pt U^-/right, {.5,-.2}/{.3,-.5}/\pt U^-/below}
            \draw[blue,<-,>=stealth]($(\i)+(\h,0)$) --++(\j)node[\loc]{$\tex$};

        \draw[fill=red!30](a1)--++(-\dt,-\dt)--++(0,2*\dt+1) --++(\dt,-\dt)--cycle;
        \draw[fill=blue!40]($(a2)+(\bx,0)$)--++(\dt,-\dt)--++(0,2*\dt+1) --++(-\dt,-\dt)--cycle;
        \foreach \i/\j in {vbr/{-\bx,1},val/{\bx,-1}}
            \draw(\i)rectangle++(\j);
        \draw[very thick,red](a1)--++(0,1);
        \draw[very thick,blue]($(a1)+(\bx-\x,0)$)--++(0,1);

        \foreach \i/\j/\tex/\loc/\col in {
        {.5*\bux-\bx,.2}/{-.3,.5}/\hat\bu_t\pt U^++\dfrac t2\bu/above/red,
        {\bx-.5*\bux,.2}/{.3,.5}/\hat\bu_t\pt U^--\dfrac t2\bu/above/blue,
        {-.6,-.2}/{-.2,-.5}/G_t+\dfrac t2\bu/below/red,
        {.6,-.2}/{.2,-.5}/-G_t-\dfrac t2\bu/below/blue}
            \draw[<-,\col](\i)--++(\j)node[\loc]{$\tex$};

        \foreach \i/\tex in {{\h,-1.5}/{(a) $\hat\bu_t\pt U$;}, {0,-1.5}/{(b)  $t\in(0,1)$.}}
            \draw(\i)node[font=\normalsize]{\tex};
        \mybrace{a1}{\bx}{$\hat\bu_t^\oplus$}; \mybrace[1]{{-.5*\bux,.5}}{\bx}{$\hat\bu_t^\oplus$};
        \mybrace[1]{{-.5*\bux+\h,.5}}{\bux}{$1-\bu^\oplus $}
        \mybrace{{-.5*\bx+\h,-.5}}{\bx}{$\hat\bu_t$}
    \end{tikzpicture}
    \caption{\label{figptG2} An illustration of $\pt U^+$, $\pt U^-$,  $G_t$ and $-G_t$ with some component of $\bu$ being zero. This figure is drawn with $\bu=(-0.3,0)$.}
\end{figure}

    Let
    \begin{gather*}
        S''^+:=\bigcup_{0\le t\le 1}\Big(\pt G_t+\dfrac t2\bu\Big)\times\{t\ep\}\quad\text{and}\quad
        S''^-:=-S''^+.\
    \end{gather*}
(see Figure \ref{figS''}) and let   $W:=\big\{(\bx,0):\,\bx\notin \cP_h\big(S''^+\cup S''^-\big)\big\}$.  We now define the \textit{quasi-hyperplane determined by $\bu$}, denoted $S_\bu$, as
        \begin{equation}\label{eqSbu}
        S_\bu:=\begin{cases}
        \pi(0),&\bu=(0,\ldots,0)\in\R^{d-1},\\
        W\cup S''^+\cup S''^-,&\bu\in(-1,1)^{d-1}\setminus\{(0,\ldots,0)\}.
        \end{cases}
        \end{equation}

    \begin{figure}
    \makeatletter
    \def\ctt[#1]{(1+2*\del*sqrt(1-#1))}
    \def\ct{\@ifnextchar[{\ctt}{\ctt[\x]}}
    \makeatother
    \begin{tikzpicture}[{>=stealth,xscale=2.5,yscale=5,font=\tiny,
        x={(-.35cm,-.15cm)}, y={(1cm,0cm)},z={(0cm,1.5cm)}}]
        \mystart[2][.35]
        \Para{1} 
        \pgfmathsetmacro{\del}{.65}
        \pgfmathsetmacro{\t}{.36}
        \coordinate(a) at (\ux/2-1/2,\uy/2-1/2,0); 
        \coordinate(a1)at($-.5*(\ux+1,\uy+1,-2*\eps)$);
        \coordinate(c)at($(a)+(1-\ux,1-\uy,0)$);  

        \sectreg{o}{0}{red!30}{1}{opacity=.5} 
        \draw[red!30,opacity=.5,dot2,thick](\del+1/2+\ux/2,-\del-1/2-\uy/2,0) --++(-2*\del-1-\ux,0,0)--++(0,2*\del+1+\uy,0);
        \foreach \i/\j/\tex/\col in {
            o/{.2,-.3,0}/U\times\{0\}/black,
            {-\ux,-\uy,0}/{-.4,-.5,0}/(U+\dfrac\bu2)\times\{0\}/red,
            {\ux,\uy,0}/{.6,-.1,0}/(U-\dfrac\bu2)\times\{0\}/blue}
            \draw[fill=\col,opacity=.35]($.5*(\i)-.5*(1,1,0)$)--++(1,0,0) --++(0,1,0) --++(-1,0,0)--cycle;

        \SSplus{o}{blue!40}{-1}
        \sectreg{o}{0}{blue!40}{-1}{opacity=.5} 
        \sectreg{o}{.63}{blue!70}{-1}{}  
        \sectline{0}{0}{0}{blue}{-1} 
        \sectreggray{o}{5.3}{4.3}{-1} 
        \sectreg{o}{\t}{red}{1}{}  
        \sectline{0}{0}{0}{red}{1} 
        \SSplus{o}{red!40}{1}
        \SSplusline{o}{1}
        \SSplusline{o}{-1}

        \def\textmp{\Big(G_t+\dfrac t2\bu\Big)\times\{t\ep\}}
        \foreach \i/\j/\col/\tex/\loc in {%
            1/\t/red!80/\textmp/left,
            -1/.63/blue/-\Big(\textmp\Big)/right}
            \draw[<-,thick,black]
            ($\i/2*({1-\ux-.9-.4*(1-\i)},{\uy-(\i+1)*.4-1.2-2*\j*\uy}, {2*\j*\eps+.02})$) --++($-\i*(0,.9,-.05)$)node[text=\col,\loc]{$\tex$};

        \foreach \i/\j/\tex/\loc/\col in
            {{\ux/2-.5,\uy-1+.1,.25}/{.5-.5*\ux,0,.2}/S^+/above/red, {\ux/2,.1,-0.06}/{0,-.25,-.4}/S^-/left/blue}
            \draw[<-](\i)--++(\j)node[\col,\loc]{$\tex$};
%
        \draw[thick,red]($(a1)+(1,0,0)$)--++(-1,0,0)--++(0,1,0);
        \draw[<-,red]($(a1)+(.6,0,0)$)--++(0,-.2,.1) node[left]{$\Big(G_1+\dfrac12\bu\Big)\times\{\ep\}$};
        \draw[blue,thick]($(a)+(1,0,-\eps)$)--++(0,1,0)--++(-1,0,0);
        \draw[<-,blue]($(a)+(.4,1,-\eps)$)--++(0,.2,-.1) node[right]{$-\Big(\Big(G_1+\dfrac12\bu\Big)\times\{\ep\}\Big)$};
        \mydashline{0,0,0}

        \foreach \i in {{0,-1,0},{-1,-1,0},{-1,0,0}}
            \draw[help lines,dashed]($(c)+(\i)$)--++(0,0,\eps);
        \foreach \i in {{0,1,0},{1,1,0},{1,0,0}}
            \draw[help lines]($(a)+(\i)$)--++(0,0,-\eps);
        \rectreg{o}{\hx}{\hy}{fill=white}
        \draw(0,0,0)node[font=\tiny]{$S^0$};
        \draw[<-](0,.1,0)--++(0,1.5,.35) node[above]{$(1-\bu^\oplus)U\times\{0\}$};
        \foreach \i in {1,-1}
            \foreach \j/\jj in {{\hx,-\hy,0}/{0,-\uy,\eps}, {-\hx,-\hy,0}/{-\ux,-\uy,\eps}, {-\hx,\hy,0}/{-\ux,0,\eps}}
            \draw[dot2]($\i/2*(\j)$)--++($\i*(\jj)$);
        \foreach \i/\tex/\sty in {o/A_1/below,{0,1-\uy,0}/A_2/below, {1-\ux,1-\uy,0}/A_3/below, {1-\ux,0,0}/A_4/below}
            \draw[fill=black]($(\i)-.5*(1-\ux,1-\uy,0)$)circle(.015) node[\sty, yshift=1pt]{$\tex$};
        \foreach \i/\tex/\sty in {
            {-1-\ux-2*\del,-1-\uy-2*\del,0}/B_1/above right,
            {-1-\ux-2*\del,1+\uy+2*\del,0}/B_2/right,
            {1+\ux+2*\del,1+\uy+2*\del,0}/B_3/below,
            {1+\ux+2*\del,-1-\uy-2*\del,0}/B_4/left}
            \draw[fill=black]($.5*(\i)$)circle(.015) node[\sty]{$\tex$};
        \foreach \i/\tex/\sty in {o/P_4/above, {-1,0,0}/P_1/above,{-1,1,0}/P_2/above}
            \draw($(.5*\hx,.5*\hy-1,\eps)+(\i)$)circle(.015) node[\sty]{$\tex$};
        \foreach \i/\tex/\sty in {o/Q_3/below, {-1,0,0}/Q_2/below,{0,-1,0}/Q_4/below}
            \draw($(.5+\ux*.5,.5+\uy*.5,-\eps)+(\i)$)circle(.015) node[\sty]{$\tex$};
        \foreach \i/\j/\tex/\col in {{0,-\del*.8,0}/{0,-\del,.05}/G_0\times\{0\}/red, {\del*.7,0,0}/{\del,-.2,-.25}/-G_0\times\{0\}/blue}
            \draw[<-]($(a)+(\hx,0,0)+(\i)$)-- ++(\j)node[font=\tiny,left,\col]{$\tex$};
    \node[below,text width=.98\textwidth,font=\normalsize] at (0,0,-.5){
        \begin{itemize}
            \item $S^+$ is the union of the trapezoids $A_1A_2P_2P_1$ and $A_1A_4P_4P_1$;
            \item $S^-$ is the union of the trapezoids $A_2A_3Q_3Q_2$ and $A_3A_4Q_4Q_3$;
            \item $G_0$ is the polygon $A_4A_1A_2B_2B_1B_4$; $-G_0$ is the polygon $A_4A_3A_2B_2B_3B_4$;
            \item $G_t+\frac t2\bu$ is the section in red lying  above $S^0$;
            \item $-G_t-\frac t2\bu$ is the section in blue lying below $S^0$;
            \item $S''^+$ is the union of $S^+$ and the four surfaces $A_4B_4P_4$, $B_1B_4P_4P_1$, $B_1B_2P_2P_1$, $A_2B_2P_2$;
            \item $S''^-$ is the union of $S^-$ and the four surfaces \!$A_4B_4Q_4$,\!\! $B_3B_4Q_4Q_3$,\!\! $B_2B_3Q_3Q_2$,\!\! $A_2B_2Q_2$.
            \end{itemize}};
        \end{tikzpicture}
        \caption{\label{figS''} An illustration of $S^0,S^+$, $S^-$, $S''^+$ and $S''^-$. This figure is drawn with $d=3$ and $\bu=(-0.3,-0.3)$.
        }
    \end{figure}

    \begin{table}
    \begin{tabular}{|m{6cm}<{\centering}|m{9cm}<{\centering}|}\hline
    Set &Property\\ \hline
    $W\!=\!\big\{(\bx,0):\,\bx\!\notin\!\cP_h\big(S''^+\cup S''^-\big)\!\big\}$&
     contains $S^0$ and is the set of all points in $S$ with height zero\\ \hline
    $S''^+=\bigcup\limits_{0\le t\le 1}\Big(\pt G_t+\dfrac t2\bu\Big)\times\{t\ep\}$&containing $S^+$\\ \hline
    $S''^-=-S''^+$&a reflection of $S''^+$ about the origin\\ \hline
    \parbox[c]{6cm}{\centering$G_t=\{\bx:\,\hat\bu_t\le2\bx^\oplus\le 2\eta\sqrt{1-t}+1$\\$+(1-t)\bu^\oplus,d_t^+(\bx)\ge d_t^-(\bx)\}$}
        &$d_h(\bx,\hat\bu_t\pt U^+)\ge d_h(\bx,\hat\bu_t(\pt U^-\cup \pt U^\updownarrow))$, if we ignore the components of $\bx$ and $\bu$ with different signs on the left side of this inequality, and only consider the components of $\bx$ and $\bu$ with different signs on the right side\\ \hline
    $S^0=(1-\bu^\oplus)U\times\{0\}$& a shrunk $(d-1)$-dimensional hypercube in $\pi(0)$\\ \hline
    $S^+=\bigcup\limits_{0\le t\le1}\Big(\hat\bu_t\pt U^+ +\dfrac t2\bu\Big)\times\{t\ep\}$& a convex combination of
        $((1\!-\!\bu^\oplus)\pt U^+)\!\times\!\{0\}$ and $(\pt U^++2^{-1}\bu)\!\times\!\{\ep\}$\\ \hline
    $S^-=-S^+$&a reflection of $S^+$ about the origin\\ \hline
    \end{tabular}
\vspace{2ex}
    \caption{\label{tab1} Some properties of the sets related to $S=W\cup S''^+\cup S''^-$, an extension of $S^0\cup S^+\cup S^-$.}
    \end{table}

\subsection{Some lemmas}
    We now introduce a tool $d_h$ called \textit{horizontal distance} in $\R^d$, which helps us get some useful analytic properties of $S_\bu$ (Lemmas \ref{lemdh} and \ref{lemSbuheight}) and plays a key role in the proof of Theorem \ref{theomain1}. The horizontal distance $d_h$ between two points $P$ and $P'$ is defined as
        \begin{equation}\label{eqdhdef}
        d_h(P,P'):=\begin{cases}
            \|\cP_h(P)-\cP_h(P')\|_\infty,&\cP_v(P)=\cP_v(P');\\
            +\infty,&\text{otherwise}.
        \end{cases}
        \end{equation}
    For convenience, we let $d_h(\bx,\bx'):=\|\bx-\bx'\|_\infty$ if $\bx,\bx'\in\R^{d-1}$. For $P\in\R^d$ and $E\subset \R^d$, let
         $$d_h(P,E):=\inf\{d_h(P,Q):\,Q\in E\}.$$
    Clearly,
        $$d_h(P,E)=d_h(P,E|_{\pi(\cP_v(P))}),$$
    where $E|_{\pi(y)}$ stands for the restriction of $E$ on $\pi(y)$, i.e.,
        $$E|_{\pi(y)}:=\{P\in E:\,\cP_v(P)=y\}.$$
    We point out that $d_h$ is not a  metric on $\R^d$, since the triangle inequality may fail. Nevertheless, the following limit property holds. We omit the proof.
    \begin{lem}\label{lemdh}
        Let $\{P_n\},\{Q_n\}\subset\R^d$ be two convergent sequences in the sense of the norm $\|\cdot\|_\infty$. If for each $n\ge1$, $P_n$ and $Q_n$ have the same height, then
            $$\lim_{n\to+\infty}d_h(P_n,Q_n)=d_h\Big(\lim_{n\to+\infty}P_n, \lim_{n\to+\infty}Q_n\Big).$$
    \end{lem}

    For $\bx\in\R^d$ satisfying $2\bx^\oplus\ge\hat\bu_t$, we divide the index set $\{j:\,|2x_j|\ge \omega_j(\hat\bu_t)\}$ into three parts as follows.
        \begin{eqnarray}
        \nonumber\Lambda_1&:=&\{j:\,x_ju_j>0,\ 2|x_j|\ge\omega_j(\hat\bu_t)\},\\ \label{eqthreelambda}
        \Lambda_2&:=&\{j:\,x_ju_j<0,\ 2|x_j|\ge\omega_j(\hat\bu_t)\},\\ \nonumber\Lambda_3&:=&\{j:\,2|x_j|\ge\omega_j(\hat\bu_t)=1,\ u_j=0\}.
        \end{eqnarray}
    So, we can rewrite the pseudo-distances $d_t^+(\bx)$ and $d_t^-(\bx)$ as
            \begin{equation}\label{eqdt+le}
            \begin{array}{rcl}
            d_t^+(\bx)&:=&\max\Big\{|x_j|-\dfrac12\omega_j(\hat\bu_t):\ j\in \Lambda_1\Big\};\\[8pt]
            d_t^-(\bx)&:=&\max\Big\{|x_j|-\dfrac12\omega_j(\hat\bu_t):\, j\in \Lambda_2\cup\Lambda_3\Big\}.
            \end{array}
            \end{equation}
    These new expressions in fact give the geometric meanings of $d_t^+$ and $d_t^-$. Roughly speaking, $d_t^+(\bx)$ is the horizontal distance between $\bx$ and $\hat\bu_t\pt U^+$ if we ignore all components of $\bx$ and $\bu$ that have different signs, while $d_t^-(\bx)$ is the horizontal distance between $\bx$ and $\hat\bu_t(\pt U^+\cup\pt U^\updownarrow)$ if we only consider the components of $\bx$ and $\bu$ that have different signs. This fact will be applied to the proofs of Lemmas \ref{lemex} and \ref{lemSbuheight}, which will be used to estimate the horizontal distances of prisms and points in the next section.
    \begin{lem}\label{lemex}
        Let $t\in[0,1)$. Suppose $P=\Big(\bx+\dfrac t2\bu,t\ep\Big)$ with $2\bx^\oplus\ge\hat\bu_{t}$. \\
        \indent $(1)$ If $P$ lies  below $S_\bu$, then
            $$d_h\Big(\bx,U-\dfrac{1+t}2\bu\Big)\ge d_h(\bx,\hat\bu_t \pt U).$$
        \indent $(2)$ If $P$ lies above $S_\bu$ or $P\in S_\bu$, then
            $$d_h\Big(\bx,U+\dfrac{1-t}2\bu\Big)\ge (1-\|\bu\|_\infty)^{1/\eta} d_h(\bx,\hat\bu_t \pt U).$$
    \end{lem}
        To prove Lemma \ref{lemex}, we need the following lemma.
    \begin{lem}\label{lem1/k}
        Suppose $0<\delta<1$ and $x\ge\eta$. Then for $t\in[0,1]$, we have
            $$x\sqrt{1-t}\ge (1-\delta)^{1/\eta}\big(x\sqrt{1-t}+(1-t)\delta\big).$$
    \end{lem}
    \begin{proof}
    Note that $\eta\in(0,1)$. By the mean-value theorem, there exists some $\xi\in(0,\delta)$ such that
            $$\dfrac{1-(1-\delta)^{1/\eta}}\delta\times \dfrac1{(1-\delta)^{1/\eta}} = \dfrac1\eta{(1-\xi)^{1/\eta-1}}\dfrac1{(1-\delta)^{1/\eta}} >\dfrac1\eta.$$
    Hence,
            $$x\sqrt{1-t}\ge \eta\sqrt{1-t}\ge \dfrac{(1-\delta)^{1/\eta} \delta\sqrt{1-t}}{1-(1-\delta)^{1/\eta}}\ge \dfrac{(1-\delta)^{1/\eta}\delta(1-t)}{1-(1-\delta)^{1/\eta}}.$$
        This implies the conclusion.
    \end{proof}
    \begin{proof}[Proof of Lemma \ref{lemex}]
        When $\bu=(0,\ldots,0)$, we know that $\hat\bu_t=\bone_{d-1}$ and hence
            $$d_h\Big(\bx, U-\dfrac{1+t}2\bu\Big) =d_h\Big(\bx,U+\dfrac{1-t}2\bu\Big) =d_h(\bx,\hat\bu_t\pt U),$$
        provided $2\bx^\oplus\ge\hat\bu_t U$.
      Thus, conclusions (1) and (2) are true.

    Next, we assume $\bu\ne(0,\ldots,0)$.

        (1) Assume $P$ lies below $S_\bu$. Then $\bx\in \myint[]G_t$, which implies that $d_t^+(\bx)>d_t^-(\bx)$, namely,
            $$\max_{j\in \Lambda_1}\Big\{|x_j|-\dfrac12\omega_j(\hat\bu_t)\Big\} >\max_{j\in \Lambda_2\cup \Lambda_3}\Big\{|x_j|-\dfrac12\omega_j(\hat\bu_t)\Big\}.$$
        So, by setting $x_j'=\omega_j(\bx')$, we have
            \begin{eqnarray*}
            d_h(\bx,\hat\bu_t\pt U)&=&\min\{\|\bx-\bx'\|_\infty:\,\bx'\in \hat\bu_t\pt U\}\\
            &=&\min\!\Big\{\!\max\big\{|x_j-x_j'|:\,j\in\Lambda_1\cup \Lambda_2\cup\Lambda_3\big\}:\,\bx'\in\hat\bu_t U\!\Big\}\\
            &=&\max\!\big\{|x_j|-\dfrac12\omega_j(\hat\bu_t):\,j\in \Lambda_1\cup\Lambda_2\cup\Lambda_3\big\}\quad\\
            &=&\max\!\big\{|x_j|-\dfrac12\omega_j(\hat\bu_t):\,j\in \Lambda_1\big\}=d_t^+(\bx),
            \end{eqnarray*}
        where $\Lambda_1,\Lambda_2,\Lambda_3$ are given as in \eqref{eqthreelambda}. Recall that $\omega_j(\hat\bu_t)=t|u_j|+1-|u_j|$.
        Without loss of generality, we assume
            \begin{equation}\label{eqdt+}
             u_1>0\quad \text{and}\quad
             d_t^+(\bx)=x_1-\dfrac12\omega_1(\hat\bu_t)=x_1-\dfrac{tu_1+1-u_1}2.
            \end{equation}
        On the other hand,
            $$d_h\Big(\bx,U-\dfrac{1+t}2\bu\Big)\ge
            x_1-\Big(\dfrac12-\dfrac{1+t}2u_1\Big)=
            x_1-\dfrac{-tu_1+1-u_1}2\ge d_t^+(\bx).$$
        So the conclusion is true in this case.

        (2) Assume $P$ lies on or above $S_\bu$. When $d_t^+(\bx)>d_t^-(\bx)$, we also get $d_h(\bx,\hat\bu_t U)=d_t^+(\bx)$. Again we suppose that assumption \eqref{eqdt+} holds.
        Since $P$ lies on or above $S_\bu$, we know $\bx\notin (2\eta\sqrt{1-t}+1+(1-t)\bu^\oplus)U^\circ$ by the definition of $S_\bu$, which implies that (see Figure \ref{fig2-2D} for the location of $P$)
            $$d_t^+(\bx)\ge\eta\sqrt{1-t}+(1-t)u_1.$$
        Notice that
            \begin{eqnarray*}
            d_h\Big(\bx,U+\dfrac{1-t}2\bu\Big)&\ge& x_1-\Big(\dfrac12+\dfrac{1-t}2u_1\Big)\\
            &=&d_t^+(\bx)-(1-t)u_1\ge\eta\sqrt{1-t}.
            \end{eqnarray*}
        Applying Lemma \ref{lem1/k} for $x=(\sqrt{1-t})^{-1}d_h\big(\bx,U+2^{-1}(1-t)\bu\big)$ and $\delta=u_1$, we have
            \begin{eqnarray*}
            x\sqrt{1-t}&=&d_h\Big(\bx,U+\dfrac{1-t}2\bu\Big)\\
            &\ge& (1-u_1)^{1/\eta}\bigg(d_h\Big(\bx,U+\dfrac{1-t}2\bu\Big) +(1-t)u_1 \bigg)\\
            &\ge& (1-\|\bu\|_\infty)^{1/\eta} d_t^+(\bx)\\
            &=& (1-\|\bu\|_\infty)^{1/\eta} d_h(\bx,\hat\bu_t\pt U).
            \end{eqnarray*}
        When $d_t^+(\bx)\le d_t^-(\bx)$, we know $d_h(\bx,\hat\bu_t \pt U)=d_t^-(\bx)$. Without loss of generality, we assume $u_1\le0$ and
            $$d_t^-(\bx)=x_1-\dfrac12\omega_1(\hat\bu_t) =x_1-\dfrac{-tu_1+1+u_1}2.$$
Then
            $$d_h\Big(\bx,U+\dfrac{1-t}2\bu\Big)\ge x_1-\Big(\dfrac12+\dfrac{(1-t)u_1}2\Big)=d_t^-(\bx).$$
        So the conclusion also holds in this case since $(1-\|\bu\|_\infty)^{1/\eta}\le1$.
    \end{proof}

    \begin{lem}\label{lemSbuheight}
        The following assertions concerning $S_\bu$ hold.
    \begin{enumerate}
    \item[(1)] For each $\bx\in\R^{d-1}$, the set $S_\bu\cap \ell_\bx$ contains exactly one element and
            $$\cP_h\big(S^+\cup S^-\cup S^0\big)=(U+2^{-1}\bu)\bigcup (U-2^{-1}\bu).$$
    \item[(2)] The height of $P\in S_\bu$ is zero if $\cP_h(P)\in (1+2\eta+\bu^\oplus)U^c$ or $\cP_h(P)\in(1-\|\bu\|_\infty)U$, where $U^c=\R^{d-1}\backslash U$.
    \item[(3)] Let $P=(\bx,y)\in S_\bu$, where $\bx\in\pt U+2^{-1}\bu$ such that
            \begin{equation}\label{eqgivencondition}
            \text{either }\omega_j(2\bx-\bu)=\pm1\quad\text{ or }\quad
            \big|\omega_j(2\bx-\bu)\big|<1-2\|\bu\|_\infty.
            \end{equation}
        Then
            $$y=\begin{cases}
                \ep,&\bx\in\pt U^++2^{-1}\bu,\\[-3pt]
            0,& \bx\in \big(\pt U^-\cup\pt U^{\updownarrow}+2^{-1}\bu\big)\backslash \big(\pt U^++2^{-1}\bu\big).
            \end{cases}$$
    \item[(4)]  Let $P=(\bx,y)\in S_\bu$, where $\bx\in\pt U-2^{-1}\bu$ such that
            \begin{equation*}
            \text{either }\omega_j(2\bx+\bu)=\pm1\quad\text{ or }\quad
            \big|\omega_j(2\bx+\bu)\big|<1-2\|\bu\|_\infty.
            \end{equation*}
        Then
            $$y=\begin{cases}
                -\ep,&\bx\in\pt U^--2^{-1}\bu,\\[-3pt]
                0,&\bx\in\big(\pt U^+\cup\pt U^{\updownarrow}-2^{-1}\bu\big)
                \backslash\big(\pt U^--2^{-1}\bu\big).
            \end{cases}$$
    \end{enumerate}
    \end{lem}

    \begin{proof}
        We first consider the case $\bu=(0,\ldots,0)$. We have $S^+=S^-=\emptyset$ and $S^0=U\times\{0\}$. So, $S_\bu\cap\ell_\bx=\{(\bx,0)\}$ and $\cP_h(S^0)=U$, which implies conclusion (1). Conclusion (2) is obvious since each point in $S_\bu$ has height zero. Note that $\pt U^++2^{-1}\bu=\pt U^--2^{-1}\bu=\emptyset$ and $\pt U^{\updownarrow}+2^{-1}\bu=\pt U$. Now conclusions (3) and (4) become: If $P=(\bx,y)\in S_\bu$ with $\bx\in\pt U$, then $y=0$ when $\bx\in\pt U^{\updownarrow}$. This is obvious since $P=(\bx,y)\in S_\bu$ implies that $y=0$. Now, we assume $\bu\ne(0,\ldots,0)$.

        (1) For a given vector $\bx\in\R^{d-1}$, we denote $\bx'=\bx+2^{-1}(t-t')\bu$, where $t,t'\in[0,1]$ satisfy $t>t'$. We begin the proof with the following claim.

    \noindent\textbf{Claim.} \textit{$(i)$ If $\bx\in \hat\bu_t(U^c\cup \pt U^+)$, then $\bx'\notin\hat\bu_{t'}U$. $(ii)$ If $\bx\in (2\eta\sqrt{1-t}+1+(1-t)\bu^\oplus)U$, then $\bx'\in(2\eta\sqrt{1-t'}+1+(1-t')\bu^\oplus) U^\circ$. $(iii)$ If $d_t^+(\bx)\ge d_t^-(\bx)$, then $d_{t'}^+(\bx')>d_{t'}^-(\bx')$.}

    \noindent\textit{Proof of the claim.}  When $\bx\in\hat\bu_t \pt U^+$, there exists an index $j$ such that $x_ju_j>0$ and $|x_j|=2^{-1}(t|u_j|+1-|u_j|)$. Without loss of generality, we assume $j=1$ and $u_1>0$. Then
            $$x_1'=x_1+\dfrac{t-t'}2u_1=\dfrac12(tu_1+1-u_1)+\dfrac{t-t'}2u_1
            =\dfrac12(t'u_1+1-u_1)+(t-t')u_1,$$
        which implies $\bx'\notin\hat\bu_{t'} U$. When $\bx\notin\hat\bu_tU$, there exists an index $j$ such that $|2x_j|>t|u_j|+1-|u_j|$. So,
            \begin{eqnarray*}
              |2x_j'|&=&|2x_j+(t-t')u_j|\ge|2x_j|-(t-t')|u_j|\\
              &>&t|u_j|+1-|u_j|-(t-t')|u_j|=t'|u_j|+1-|u_j|,
            \end{eqnarray*}
        which again implies $\bx'\notin\hat\bu_{t'}U$, and thus (i) holds.

 From the fact that $\eta>\|\bu\|_\infty$, we have
            \begin{eqnarray*}
            \lefteqn{2\eta\sqrt{1-t}+(1-t)\bu^\oplus +(t-t')\|\bu\|_\infty}\\
            &=&2\eta\sqrt{1-t'}+2\eta\dfrac{t'-t}
            {\sqrt{1-t}+\sqrt{1-t'}}+(1-t)\bu^\oplus+(t-t')\|\bu\|_\infty
            \\
            &\le&2\eta\sqrt{1-t'}+(1-t')\bu^\oplus +(\|\bu\|_\infty-\eta)(t-t')\\
            &<&2\eta\sqrt{1-t'}+(1-t')\bu^\oplus,
            \end{eqnarray*}
        which implies
            \begin{equation}\label{eqptsubsetint}
            (2\eta\sqrt{1-t}+1+(1-t)\bu^\oplus)U+\dfrac{t-t'}2\bu
            \subset(2\eta\sqrt{1-t'}+1+(1-t')\bu^\oplus)U^\circ.
            \end{equation}
        Thus, (ii) holds.

        As for (iii), we let $\Lambda_1,\Lambda_2,\Lambda_3$ be given as in  \eqref{eqthreelambda}. Since $\bx\notin\hat\bu_t U$, not all sets $\Lambda_1,\Lambda_2,\Lambda_3$ are empty. The assumption  $d_t^+(\bx)\ge d_t^-(\bx)$ implies $\Lambda_1\ne\emptyset$ by \eqref{eqdt+le}. Without loss of generality, we assume
            \begin{equation*}
             u_1>0\quad\text{and}\quad
             d_t^+(\bx)=x_1-\dfrac12\omega_1(\hat\bu_t)=x_1 -\dfrac{tu_1+1-u_1}2.
            \end{equation*}
        Then
            \begin{eqnarray}\nonumber
            d_{t'}^+(\bx')&\ge& \Big(x_1+\dfrac{t-t'}2u_1\Big)-\dfrac12(t'u_1+1-u_1)\\
            \label{eqdt'>dtx}&=&x_1-\dfrac12(tu_1+1-u_1)+(t-t')u_1>d_t^+(\bx).
            \end{eqnarray}
        Note that
            $$d_t^-(\bx)=\max\Big\{|x_j|-\dfrac12\omega_j(\hat\bu_t):\, j\in\Lambda_2\cup\Lambda_3\Big\}.$$
        When $j\in\Lambda_2$ (in this case, $x_ju_j<0$),
            \begin{eqnarray*}
            \lefteqn{\Big|x_j+\dfrac{t-t'}2u_j\Big| -\dfrac12(t'|u_j|+1-|u_j|)}\\
            &=&|x_j|-\dfrac{t-t'}2|u_j|-\dfrac12(t'|u_j|+1-|u_j|)\\
            &=&|x_j|-\dfrac12(t|u_j|+1-|u_j|),
            \end{eqnarray*}
       where the first equality holds because $2|x_j|\ge t|u_j|+1-|u_j|$ and
            $$t|u_j|+1-|u_j|-(t-t')|u_j|=t'|u_j|+1-|u_j|>0.$$
        When $j\in\Lambda_3$ (in this case, $u_j=0$),
            $$\Big|x_j+\dfrac{t-t'}2u_j\Big|-\dfrac12(t'|u_j|+1-|u_j|)
            =|x_j|-\dfrac12(t|u_j|+1-|u_j|)$$
        So, for $j\in\Lambda_2\cup\Lambda_3$, we have
            $|x_j'|-\dfrac12\omega_j(\hat\bu_{t'}) =|x_j|-\dfrac12\omega_j(\hat\bu_t),$
        which implies
            $$d_t^-(\bx)= d_{t'}^-(\bx').$$
        Hence, (iii) holds by \eqref{eqdt'>dtx}, which completes the proof of the claim.

        From the expression of $G_t$ (see \eqref{eqGt+-}), we get
            \begin{gather*}
            \pt G_t=E_1\bigcup E_2\bigcup E_3,\\
            \myint[] G_{t'}=\{\bx\in\mathbb R^{d-1}:\,d_{t'}^+(\bx)> d_{t'}^-(\bx),\hat\bu_{t'}<2\bx^\oplus <2\eta\sqrt{1-t'}+1+(1-t')\bu^\oplus\},
            \end{gather*}
        where
        $$
        \begin{aligned}
            &E_1:=\{\bx\in \hat\bu_t \pt U:\,d_{t}^+(\bx)\ge d_t^-(\bx)\}\subset(2\eta\sqrt{1-t}+1 +(1-t)\bu^\oplus)U^\circ,\\
            &E_2:=\{\bx\in(2\eta\sqrt{1-t}+1+(1-t)\bu^\oplus)\pt U:\,d_t^+(\bx)\ge d_t^-(\bx) \},\\
            &E_3:=\{\bx:\,d_t^+(\bx)=d_t^-(\bx), \hat\bu_{t'}<2\bx^\oplus<2\eta\sqrt{1-t'}+1 +(1-t)\bu^\oplus \}.
        \end{aligned}$$
        Now, the claim implies that $\bx'\in \myint[] G_{t'}$ when $\bx\in E_1\cup E_2\cup E_3$. This implies the following relationship among the $G_t$: for $t,t'\in[0,1]$ with $t>t'$,
        \begin{equation}\label{eqGtsubset}
        \pt G_t+\dfrac{t}2\bu \subset
        \myint[] G_{t'}+\dfrac {t'}2\bu.
        \end{equation}
        Now, it follows from the construction of $S_\bu$ that $\ell_\bx\cap S_\bu$ is a singleton.
        The latter part of the conclusion in (1) is obvious since
            $\cP_h(S^+)$ is the closure of $(U+2^{-1}\bu)\backslash(1-\bu^\oplus )U$, $\cP_h(S^-)$ is the closure of $(U-2^{-1}\bu)\backslash(1-\bu^\oplus )U$, and $\cP_h(S^0)=(1-\bu^\oplus )U$.

        (2) The height of $P\in S_\bu$ is not equal to zero only if $P\in S''^+\cup S''^-$, the horizontal coordinate of which is in
            $$(1+2\eta+\|\bu\|_\infty)U\backslash (1-\|\bu\|_\infty)U.$$
        So, the conclusion follows.

        We now show (3). By (1), the given condition $\bx\in\pt U+2^{-1}\bu$ implies that $P\in S^+\cup S^0$. Since $S^0$ lies in $\pi(0)$ and $S^+$ is  the union of the convex combinations $\big((1-\bu^\oplus )\pt U^+\big)\times\{0\}$ and $(\pt U^++2^{-1}\bu)\times\{\ep\}$,  the height of $P$ is equal to $\ep$ if and only if $\bx\in \pt U^++2^{-1}\bu$. Next, we assume $\bx\in\pt U^-\cup\pt U^{\updownarrow}+2^{-1}\bu$ but $\bx\notin \pt U^++2^{-1}\bu$. If the height of $P$ is positive, then $P\in S^+$ and $0<y<\ep$. By the expression of $S^+$ (see \eqref{eqfourS})
        there exist $\by\in\pt U^+$ and $t\in(0,1)$ such that
            $$\bx=\hat\bu_t \by+\dfrac t2\bu.$$
       Since $\by\in\pt U^+$, we choose some index $j$ such that $2\omega_j(\by)=\sign(u_j)\ne0$. So,
            \begin{eqnarray*}
            \big|\omega_j\big(2\bx-\bu\big)\big|&=&
            \big|(t|u_j|+1-|u_j|)\cdot2\omega_j(\by)+tu_j-u_j\big| \\
            &=&\big|2tu_j+\sign(u_j)-2u_j\big|=\big|1-2(1-t)|u_j|\big|\\
            &\in&
            \big(1-2|u_j|(1-t),1\big)\subset \big(1-2\|\bu\|_\infty,1\big).
            \end{eqnarray*}
        This contradicts condition \eqref{eqgivencondition}  and thus (3) holds.

        The proof of (4) is analogous to that of (3). This completes the proof.
    \end{proof}

\section{A key homeomorphism on $\R^d$\label{S:section3}}
    \subsection{Definitions of transformations}\label{subsection3_1}
    In this subsection we present four transformations on $\R^d$, which will be used to construct the desired homeomorphism on $\R^d$. Before doing this, we introduce some notation, in addition to those in Subsection 2.3. Let $b>0$ be a fixed positive constant as in Subsection 2.3. We choose $\ep$ stated as in Subsection 2.3 small enough such that the union $\cup_{k=0}^r S_k$ is disjoint, where
        $S_0=\pi(0),S_r=\pi(b)$
    and
        \begin{equation}\label{eqSk}
            S_k=S_{\bu_k}+\Big(\dfrac12\bu_k+\bv_{k-1},y_k\Big),\quad 1\le k<r.
        \end{equation}
   Hence $S_{\bu_k}$ is defined as in \eqref{eqSbu} with $\bu$ replaced by $\bu_k$ and $\eta$ is a constant satisfying
        \begin{equation}\label{eqetarequire}
        \max_{0<k<r}\{\|\bu_k\|_\infty\}<\eta<1.
        \end{equation}
    For $1\le k<r$, $S^0_k,S^+_k,S^-_k,S'^+_k$, and $S'^-_k$ are defined similarly. For convenience, we set
        $$S_0^0=U\times\{0\},\quad S_r^0=(U+\bv_{r-1})\times\{b\},\quad\text{and}\quad S_0^+=S_0^-=S_r^+=S_r^-=\emptyset.$$
 We use $P$ to denote $(\bx,y)\in\R^{d-1}\times[0,b]$, unless stated otherwise.
    Let $\bt\in\R^{d-1}$ and $\bx(y)$, $\brho(y)$ be defined as in \eqref{eqbxy}, \eqref{eqrhoy} respectively with $\bu_1,\ldots,\bu_{r-1}$ given as in Subsection \ref{sectother}. Define an \textit{affine transformation} $f_\bt$, a \textit{squeezing transformation} $f_S$, and a \textit{translation} $f_T$ as follows
       \begin{eqnarray}
        f_\bt(P)&:=&P+\Big(\dfrac {y\bt}b,0\Big),\label{eqfB}\\
       f_S(P)&:=&\big(\brho(y),1\big)P, \label{eqfS}\\
       f_T(P)&:=&P+\big(\bx(y),0\big). \label{eqfT}
       \end{eqnarray}
    Figure~\ref{figillus}(a)--(d) show how these maps transform the prism.

    We introduce some additional notation. Assume quasi-hyperplanes $S_k$ are given as in Section~\ref{hyperplane}. For a vertical line $\ell_\bx$ (given as in \eqref{eqellbx}), denote by $Q_k=Q_k(\bx)$ the intersection of $\ell_\bx$ and $S_k$ (see Lemma \ref{lemSbuheight}(1) for the existence and uniqueness of $Q_k$) and denote $y_k'=y_k'(\bx)=\omega_d(Q_k)$ for $0\le k\le r$. Clearly $y_0'=0$ and $y_r'=b$. We now define the fourth transformation $f_F$, called a \textit{flattening transformation}, as follows. Recall $P=(\bx,y)\in\R^{d-1}\times[0,b]$. Let $k=k(P)$ be equal to $r$ or the unique integer that satisfies $y_k'\le y<y_{k+1}'$ according to $y=b$ or $y<b$, where $y_k$ is as in \eqref{eqyk}. Define $f_F$ on $\R^{d-1}\times[0,b]$ as (see Figure \ref{figprojection})
        \begin{equation}\label{eqfR}
        f_F(P):=\big(\bx,\varphi(P)\big),
        \end{equation}
    where %
        \begin{equation}\label{eqvarphidef}
        \varphi(P):=
        \begin{cases}
            \dfrac{b}r\cdot\dfrac{y-y_k'(\bx)}{y_{k+1}'(\bx)-y_k'(\bx)}+y_k, &y<b;\\
            b,&y=b.
        \end{cases}
        \end{equation}
    Figure~\ref{figillus}(d)--(e) show the restriction of $f_F$ to the figure in (d).

    We comment on the functions $\varphi$ and $f_F$. Let $X_k\subset\R^d$ be the horizontal slab-shaped closed region  bounded on the top and bottom by the irregular surfaces $S_k$ and $S_{k+1}$. Then $f_F$ maps each $X_k$ onto the horizontal slab-shaped region $\R^{d-1}\times[y_k,y_{k+1}]$, and is piecewise linear on each vertical line segment. More precisely, on each vertical line segment
    $$\{\bx\}\times[0,b]=\bigcup_{k=0}^{r-1}\{\bx\} \times[y_k',y_{k+1}'],$$
    $f_F$ maps the $k$th line segment $\{\bx\}\times[y_k',y_{k+1}']$ to the line segment $\{\bx\}\times[y_k,y_{k+1}]$. Furthermore, we get the following estimation on the height difference between $P$ and its image $f_F(P)$,  which will be used to prove the continuity of the desired homeomorphism  (Lemma \ref{lemcountinuity}).
    \begin{lem}\label{lemhtdiff1}
        Let $f_F$ be defined as above and let $\ep>0$ be as in \eqref{eqep}. Then, for each $P\in\R^{d-1}\times[0,b]$, we have
        \begin{equation}\label{eqvarphi}
            \big|\cP_v\big(f_F(P)\big)-\cP_v(P)\big|\le \ep.
        \end{equation}
    \end{lem}
    \begin{proof}
        Let $P=(\bx,y)\in\R^{d-1}\times[0,b]$. We may assume the height of $P$ is in $[y_k',y_{k+1}']$ for some $0\le k<r$,  since $\cP_v\big(f_F(P)\big)=\cP_v(P)=b$ when $y=b$. Note that
            $$\cP_v\big(f_F(P)\big)-\cP_v(P)=\varphi(P)-\cP_v(P)
            =\dfrac{b}r\cdot\dfrac{y-y_k'(\bx)}{y_{k+1}'(\bx)-y_k'(\bx)} +y_k-y$$
        is a linear function of the variable $y$. The maximum in  absolute value of such a linear function can only be obtained at the end-points of the interval $[y_k',y_{k+1}']$, which are
        $$\varphi\big((\bx,y_k')\big)-y_k'=y_{k}-y_k'\text{ and }
        \varphi\big((\bx,y_{k+1}')\big)-y_{k+1}'=y_{k+1}-y_{k+1}'.$$
        From the construction of $S_k$ and $S_{k+1}$, for all $\bx\in\R^{d-1}$, we know
            $$y_n'\in(y_n-\ep,y_n+\ep),\quad n=k,k+1.$$
        Hence, \eqref{eqvarphi} follows and the proof is complete.
    \end{proof}
    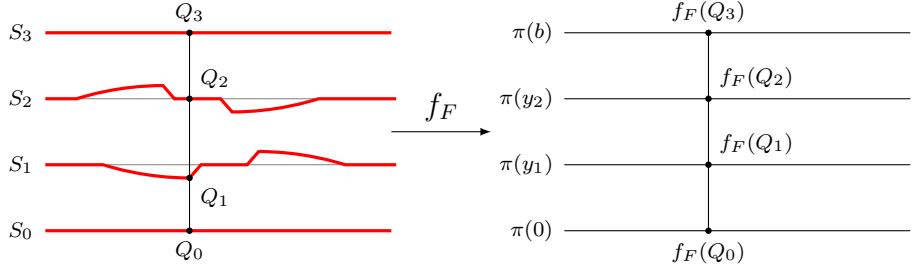
\begin{figure}
    \centering{}
    \begin{tikzpicture}[scale=1.75,declare function={f(\t)=-2*(\t-1)+1;}]
        \pgfmathsetmacro{\y}{.1}
        \pgfmathsetmacro{\v}{1.5}
        \pgfmathsetmacro{\h}{1.3}
        \pgfmathsetmacro{\fc}{3}
        \foreach \i/\tex/\sty in {
            0/\pi(0)/below,1/\pi(y_1)/above right,2/\pi(y_2)/above right,3/\pi(b)/above}{
            \draw(\fc*\h,\i*\v/3)node[font=\tiny,left]{$\tex$}--++(2*\h,0);
            \draw[fill=black](5*\h/6+\fc*\h,\i*\v/3)circle(.02) node[font=\tiny,\sty]{$f_F(Q_\i)$};}
        \draw(\fc*\h+5*\h/6,0)--++(0,\v);
        \foreach \i in {1,2}
            \draw[help lines](0,{\i*\v/3})--++(2*\h,0);
        \draw[->,>=latex](2*\h,\v/2)--node[above]{$f_F$}++(\h/2-\y+.2,0);

        \foreach \i in {0,3}
            \draw[red,very thick,text=black] (0,\i*\v/3)node[font=\tiny,left]{$S_\i$} --++(2*\h,0);
        \foreach \i/\j in {1/0,2/-.2}
            \draw[red,very thick,text=black](0,{\i*\v/3}) node[left,font=\tiny]{$S_\i$}--++(\h/3+\j,0) parabola bend ++ (.5*\h,{-f(\i)*\y}) ++(.5*\h,{-f(\i)*\y}) --++(\h/15,{f(\i)*\y})--++(4*\h/15,0) --++(\h/15,{f(\i)*\y}) parabola bend ++(0,0)++(.5*\h,{-f(\i)*\y}) --++(\h/3-\j-.05,0);
        \draw(5*\h/6,0)--++(0,\v);
        \foreach \i/\tex/\sty in {\v/3/above,{2*\v/3}/2/above right, {\v/3-\y}/1/below right, 0/0/below}
        \draw[fill=black](\h/2+\h/3,\i)circle[radius=.02] node[font=\tiny,\sty]{$Q_\tex$};
    \end{tikzpicture}
    \caption{\label{figprojection} An illustration of $f_F$. The vertical line segment $Q_0Q_3=\{\bx\}\times[0,b]$ intersects $S_i$ at $Q_i$, $i=0,\ldots,3$. The last coordinate of $f_F(Q_i)$ is $\varphi(Q_i)$. This figure is drawn with $r=3$.}
    \end{figure}

    \subsection{Some lemmas}
    Define a map $f$ as
        \begin{equation}\label{eqhomf}
        f(P)=f_F\circ f_{T}\circ f_{S}(P),\quad P\in\R^{d-1}\times[0,b].
        \end{equation}
    Denote $f_T\circ f_S$ by $f_{TS}$ for short.
    One can easily check that the mapping $f$ is a homeomorphism  on $\R^{d-1}\times[0,b]$. In this subsection we prove some properties of the several mappings on $\R^{d-1}\times[0,b]$ that we have defined.

    In view of the complexity of $f_F$ or $f$ on the set $\pt U\times[0,b]$, we will decompose $\pt U\times[0,b]$. For $\bi=i_1\cdots i_{d-1}\in\{-1,0,1\}^{d-1}$, let
    \begin{equation*}
        \begin{aligned}
        &U_\bi:=\{\bx\in U:\,\omega_j(\bx)=2^{-1}i_j\mbox{ if } i_j\ne0 \text{ and }\omega_j(\bx)\in [-2^{-1},2^{-1}]\text{ if }i_j=0\};\\
         &\myint[] U_\bi=U_\bi^\circ=\{\bx\in U_\bi:\,
         \omega_j(\bx)\in(-1/2,1/2)\text{ if } i_j=0\};\\
         &\pt U_\bi=U_\bi\backslash \myint[] U_\bi
    \end{aligned}
    \end{equation*}
    and
        \begin{gather*}
         H_\bi:=U_\bi\times[0,b],\quad \myint H_\bi:=\myint[] U_\bi\times[0,b],\quad
         \pt_v H_\bi:=\pt U_\bi\times[0,b].
    \end{gather*}
    The following facts are clear. First, $U=U_\bi$ and $\pt U=\pt U_\bi$ when $\bi=0^{d-1}$. Second, $\pt U_\bi=U_\bi$, $\myint[] U_\bi=\emptyset$ (and hence $\myint H_\bi=\emptyset$) if no component of $\bi$ is zero. Finally,
        $$\pt U_\bi=\{\bx\in U_\bi:\,\exists\text{ }
            j \text{ s.t. }
            \omega_j(\bx)=\pm2^{-1}\text{ and }i_j=0\}$$
    if $\bi$ has at least one zero component. We call $\pt_v H_\bi$ the \textit{vertical boundary} of $H_\bi$, which is different from $\pt H_\bi$ in the set $U_\bi^\circ\times\{0,b\}$.
    For simplicity, denote $H_{0^{d-1}}$ by $H$. For a set $E\subset \R^{d-1}\times[0,b]$, we use $E^{TS}$ to denote the image of $E$ under $f_{TS}$, e.g.,
        $$H^{TS}=f_{TS}(H).$$

    The following lemma is obvious; we omit its proof.
    \begin{lem}\label{lemfBT}
        Let $f_\bt$ and $f_T$ be defined as above. Then for $P,P'\in\R^{d-1}\times[0,b]$, we have
            $$d_h(P,P')=d_h\big(f_\bt(P),f_\bt(P')\big) =d_h\big(f_T(P),f_T(P')\big).$$
        Consequently, For any compact set $E\subset\R^d$ with $\cP_v(E)=[0,b]$, we have
            \begin{gather*}
            d_h\big(P,E\big)=d_h\big(f_\bt(P),f_\bt(E)\big)
            =d_h\big(f_T(P),f_T(E)\big),\quad P\in\R^{d-1}\times[0,b].
            \end{gather*}
        In particular, for $\bi\in\{-1,0,1\}^{d-1}$ and $P\in\R^{d-1}\times[0,b]$,
            $$d_h(P,\pt_v H_\bi)=d_h\big(f_\bt(P),f_\bt(\pt_v H_\bi)\big)=d_h\big(f_T(P),f_T(\pt_v H_\bi)\big).$$
    \end{lem}
    Recall the constant $c=\sum_{k=1}^{r-1}\|\bu_k\|_\infty$.
    \begin{lem}\label{lemfS}
        Use the above notation. If $P,P'$ are both in $R^{d-1}\times[0,b]$, then
            \begin{equation}\label{eqpp'}
            d_h\big(f_S(P),f_S(P')\big)\ge (1-c)d_h(P,P').
            \end{equation}
        Consequently, for $\bi\in\{-1,0,1\}^{d-1}$ and for $P\in H_\bi$, we have
            \begin{equation}\label{eqfSPHbi}
            d_h\big(f_S(P),f_S(\pt_vH_\bi)\big) \ge (1-c)d_h(P,\pt_v H_\bi).
            \end{equation}
        Moreover, for $P\in\R^{d-1}\times[0,b]$,
            \begin{equation}\label{eqdhfS}
            d_h\big(f_S(P),f_S(\pt_vH)\big)\ge (1-c)d_h(P,\pt_vH).
            \end{equation}

    \end{lem}
    \begin{proof}
       If $\cP_v(P)\ne\cP_v(P')$, \eqref{eqpp'} holds since both sides are equal to $+\infty$. So, we assume $P=(\bx,y),P'=(\bx',y)\in \R^{d-1}\times[0,b]$. Hence,
            \begin{eqnarray*}
            d_h\big(f_S(P),f_S(P')\big)&=&\|\brho(y)(\bx-\bx')\|_\infty
            \ge(1-c)\|\bx-\bx'\|_\infty\\
            &=&(1-c)d_h(P,P').
            \end{eqnarray*}
    Let $\bi\in\{-1,0,1\}^{d-1}$ be fixed. For $P=(\bx,y)\in H_\bi$, it follows from the definition of $f_S$ that
        \begin{eqnarray*}
          d_h\big(f_S(P),f_S(\pt_vH_\bi)\big)&=&
          \inf\big\{d_h\big(f_S(P),f_S(P')\big):\,P'\in \pt_vH_\bi\big\} \\
          &\ge&\inf\big\{(1-c)d_h(P,P'):\,P'\in\pt_vH_\bi\big\}\\
          &=&(1-c)d_h(P,\pt_vH_\bi).
        \end{eqnarray*}
    The inequality in \eqref{eqdhfS} can be proved similarly. This completes the proof.
    \end{proof}

\begin{figure}
    \begin{tikzpicture}[{xscale=4,yscale=3, x={(-.45cm,-.05cm)}, y={(.3cm,-.1cm)},z={(0cm,1.65cm)}}]
        \mystart[2][.08]
        \pgfmathsetmacro{\del}{.45}
        \Sextension{gray!80}[1]{4.3}{4.3}

        \draw(0,0,-.1)node[below]{(a1) $H^{TS}$.};
    \end{tikzpicture}

    \begin{tikzpicture}[{>=stealth,xscale=1.5,yscale=3.5, x={(-.35cm,-.05cm)}, y={(1cm,-.0cm)},z={(0cm,1cm)}}]
        \mystart[2][.1]
        \pgfmathsetmacro{\del}{.45}
        \coordinate (oo1) at (-\ux/2,-\uy/2,1/3);
        \coordinate (oo2) at (-\ux-\uux/2,-\uy-\uuy/2,2/3);
        \def\colorgray{gray}
        \coordinate (center) at (0,2.2*\h,0);
    \foreach \i/\sto/\ax/\ay/\az in
        {1/oo1/{-\ux/2}/{-\uy/2}/{1/3},
            2/oo2/{-\ux-\uux/2}/{-\uy-\uuy/2}/{2/3}}{
        \Para{\i}
        \coordinate (ctmp) at ($(center)+(\sto)$);
        \pgfmathsetmacro{\pix}{6.3}
        \pgfmathsetmacro{\piy}{3.3}
        \rectreg{ctmp}{\pix}{\piy}{fill=gray!50}
        \sectreg{ctmp}{0}{red!50}{1}{opacity=1} 
        \lateral*{0,0,1}{$(ctmp)-(0,0,\eps)$}[\i/3-\eps][opacity=.5]
        \lateral{0,0,1}{$(ctmp)+(0,0,\eps)$}[\i/3+\eps][opacity=.75]

        \sectreg{ctmp}{0}{blue!50}{-1}{opacity=1}
        \path[fill=white]($(ctmp)+(-\pix/2,\piy/9,0)$)--++(0,\piy/2-\piy/9,0) --++(\pix/5*2,0,0)--cycle;
        \path[fill=white]($(ctmp)+(\pix/2,-\piy/9,0)$)--++(0,-\piy/2+\piy/9,0) --++(-\pix/5*2,0,0)--cycle;
        \foreach \j/\colone/\coltwo in {1/red!20/blue!20,-1/blue!20/red!20}
            {
            \path[fill=\colone]($(ctmp)+\j*(\parhx/2,-\parhy/2,0)$) --++($\j*(0,-\paray,0)$)--++ ($\j*(\del+\parax+.5-\parhx/2,-\del-.5+\parhy/2,0)$) --cycle;
            \path[fill=\coltwo]($(ctmp)+\j*(\parhx/2,-\parhy/2,0)$) --++($\j*(\paray,0,0)$)--++
            ($\j*(\del+.5-\parhx/2,-\del-\paray-.5+\parhy/2,0)$) --cycle;
            }
    }
    \foreach \i/\ix/\iy in {0/0/0,1/-\ux/-\uy,2/-\ux-\uux/-\uy-\uuy}{
        \coordinate(centmp) at($(center)+(\ix,\iy,\i/3)$);
        \foreach \j/\sty in {{-.5,-.5,0}/dot2,{.5,-.5,0}/very thin, {.5,.5,0}/very thin, {-.5,.5,0}/very thin}
            \draw[\sty]($(centmp)+(\j)$)--++(0,0,1/3);
        \ifnum\i=0
            \rect{centmp}11{fill=gray}
            \rect{$(centmp)+(0,0,1/3)$}11{opacity=.5,fill=blue!30}[very thin]
            \rect{$(center)+(o3)$}\hx\hy{fill=white}[dot2]
            \else\ifnum\i=1
                \rect{centmp}11{opacity=.5,fill=blue!20}
                \rect{$(centmp)+(0,0,1/3)$}11{opacity=.5,fill=blue!30}[very thin]            \rect{$(center)+(o6)$}\rx\ry{fill=white}[dot2]
                \else
                    \rect{centmp}11{opacity=.5,fill=blue!20}
                    \rect{$(centmp)+(0,0,1/3)$}11{fill=cyan!50}[very thin]
                \fi
           \fi
        }
        \draw[->](0,1.5,.5)--node[above]{flatten}++(0,.45*\h,0);

        \Sextension{gray!80}[1]{6.3}{3.3}\shrinkprismcenter{1}{0,0,0}
        \foreach \i/\j/\tex/\sty in {o1/0/O_1/below, o2/\eps/O_2/left,
        o3/0/O_3/right, o4/-\eps/O_4/right, o5/\eps/O_5/left, o6/0/O_6/right,
        o7/-\eps/O_7/left, o8/0/O_8/above}{
            \draw(\i)node[font=\tiny,\sty]{$\tex$};
            \draw[fill=black]($(\i)+(center)+(0,0,\j)$)circle(.02);}
        \draw[dashed](center)node[below,font=\tiny]{$f_F(O_1)$} --++(0,0,1/3)--++(-\ux,-\uy,0) --++(0,0,1/3)--++(-\uux,-\uuy,0) --++(0,0,1/3)node[above,font=\tiny]{$f_F(O_8)$};
        \foreach \i/\j/\tex/\sty in {{0,0,1/3}/{0,1,-.25}/2/right, {-\ux,-\uy,1/3}/{0,1,.2}/4/right, {-\ux,-\uy,2/3}/{0,1,.25}/5/right, {-\ux-\uux,-\uy-\uuy,2/3}/{0,-1,.2}/7/left}
            \draw[->]($(center)+(\i)$)--++(\j)node[\sty,font=\tiny]{$O_\tex'$};
        \foreach \i/\tex in {o3/O_3,o6/O_6}
            \draw[->]($(center)+(\i)$) --++(0,-1,-.2)node[font=\tiny,below]{$f_F(\tex)$};;
        \draw(0,0,-.2)node[below]{(a2) $H^{TS}$.};
        \foreach \i/\j in {1/.35,2/.7}
            \draw(0,1.5+.25*\h,\j)node[font=\tiny]{$\pi(y_\i)$};

        \draw($(center)+(0,0,-.2)$)node[below]{(b) $f(H^{TS})$.};
    \end{tikzpicture}

        \caption{\label{figHTS} The sets $H^{TS}$ and $f_F(H^{TS})$, along with the quasi-hyperplanes. (a1) and (a2) are the same set $H^{TS}$, with different viewpoints. The center axis of $H^{TS}$ is the dashed line $O_1O_2\cdots O_8$, the coordinates of which are $\big(\bx(y),y\big)$ and $0\le y\le b$. The restriction of $H^{TS}$ on the hyperplane $\pi(y)$ is a $(d-1)$-dimensional hypercube with edge length vector $\brho(y)$. The line segment $O_2'f_F(O_3)O_4'$ is in the plane $\pi(y_1)$, while the line segment $O_5'f_F(O_6)O_7'$ is in the plane $\pi(y_2)$.
        This figure is drawn with $d=3$,  $\bu_1=(-0.3,-0.2)$, $\bu_2=(-0.5,-0.35)$, $b=1$, and $\ep=0.1$.}
    \end{figure}
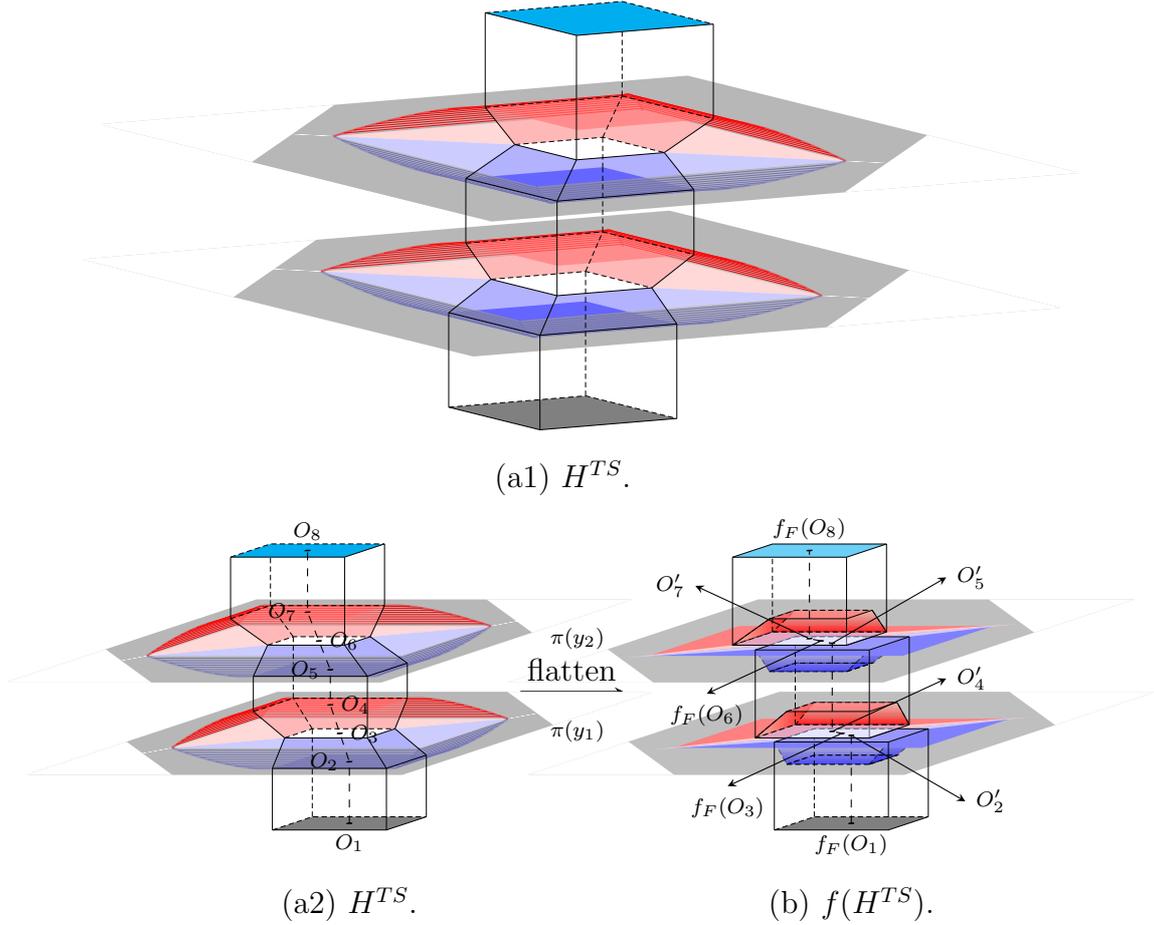

   We now discuss the relationship between the horizontal distance from $P$ to $\pt_v H^{TS}$ as defined in \eqref{eqdhdef} and the horizontal distance between their images under the transformation $f_F$. The transformation $f_F$ is  determined by the quasi-hyperplanes defined in Subsection \ref{subsection3_1}, or more precisely in Section \ref{hyperplane}. In order to state our idea clearly, we first comment on $H^{TS}$ and $f_F(H^{TS})=f(H)$ (see Figure \ref{figHTS} for an illustration). The set $H^{TS}$ (see Figure \ref{figHTS} (a1,a2)) is a $d$-dimensional polytope (whose height is in $[0,b]$), the center axis of which is $\bx(y)$, $0\le y\le d$, and the restriction of which on the hyperplane $\pi(y)$, $0\le y\le b$, is a $(d-1)$-dimensional hypercube with edge length vector being $\brho(y)$. The $d$-dimensional polytope $f_F(H^{TS})$ (see Figure \ref{figHTS}(b)) is a union of $r$ small vertical $d$-dimensional prisms (the height of the $k$th small $d$-dimensional prism is in $[y_k,y_{k+1}]$), the restriction of which on the hyperplane $\pi(y)$, $0\le y\le b$, is a unit hypercube or a union of two unit hypercubes according to $y\not\in\{y_1,\ldots,y_{r-1}\}$ or not.

    We state a lemma concerning the boundary of $H_\bi$. We omit the proof since the conclusion follows directly from the definitions of $f_T,f_S,f_F$ and $\brho(\cdot)$, $\bx(\cdot)$.

    \begin{lem}\label{lemHbi}
        For each $\bi\in\{-1,0,1\}^{d-1}$, the sets $H_\bi^{TS}$ and $f_F(H_\bi^{TS})$ can be rewritten as
        \begin{gather}
            \label{eqHbiTS}
            H_\bi^{TS}=\bigcup_{k=0}^{r-1} \bigcup_{y\in[y_k,y_{k+1}]}\big({\bm \brho}(y)U_\bi+\bx(y)\big) \times\{y\},\\
            \label{eqHbiRTS}
            f_F(H_\bi^{TS})=\bigcup_{k=0}^{r-1} (U_\bi+\bv_k)\times[y_k,y_{k+1}].
        \end{gather}
    \end{lem}

    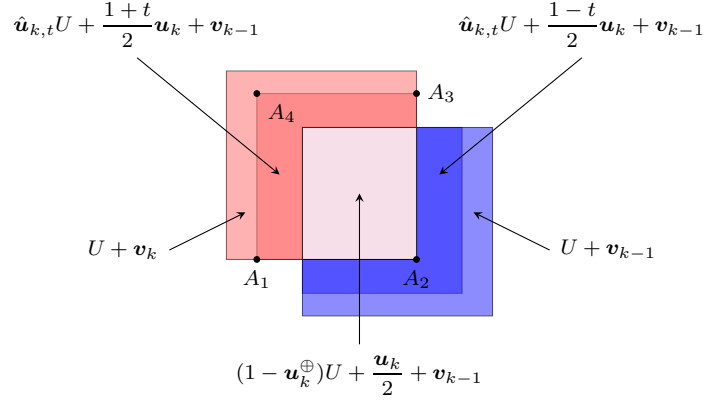
\begin{figure}
        \begin{tikzpicture}[scale=2.5,>=stealth,font=\tiny]
            \pgfmathsetmacro{\x}{.2} 
            \pgfmathsetmacro{\y}{.15} 
            \pgfmathsetmacro{\t}{.4}  
            \pgfmathsetmacro{\h}{0}  
            \pgfmathsetmacro{\xt}{\x*\t}
            \pgfmathsetmacro{\yt}{\y*\t}
            \pgfmathsetmacro{\bx}{(1-2*\x)}
            \pgfmathsetmacro{\by}{(1-2*\y)}
            \pgfmathsetmacro{\bxt}{(1-2*\xt)}
            \pgfmathsetmacro{\byt}{(1-2*\yt)}
            \coordinate (a1) at (.5-\x,-.5+\y);
            \coordinate (b1) at ($-1*(a1)$);
            \coordinate (a2) at ($(b1)+(\bxt,-\byt)$);
            \coordinate (b2) at ($-1*(a2)$);
            \coordinate (a3) at ($(b1)+(1,-1)$);
            \coordinate (b3) at ($-1*(a3)$);
            \coordinate (c) at (-.5+\x,.5-\y);

            \foreach \i/\col in {1/blue!90,-1/red!60}
                \foreach \j/\k in {a2/{-\bxt,\byt},a3/{-1,1}}
                    \draw[fill=\col,opacity=.5]($\i*(\j)$) rectangle++($\i*(\k)$);

            \foreach \j/\k/\tex/\sty in {{-.5*\bxt,.1}/{-.75,.62}/ \hat\bu_{k,t}U+\dfrac{1+t}2\bu_k+\bv_{k-1}/above, {-.6,-.1}/{-.4,-.2}/U+\bv_{k}/left}
                \draw[<-](\j)--++(\k)node[\sty]{$\tex$};
            \foreach \j/\k/\tex/\sty in {{.5*\bxt,.1}/{.75,.62}/\hat\bu_{k,t}U+\dfrac {1-t}2\bu_k+\bv_{k-1}/above,{.6,-.1}/{.4,-.2}/U+\bv_{k-1}/right}
                \draw[<-](\j)--++(\k)node[\sty]{$\tex$};
           \draw[fill=white,opacity=.8](a1)rectangle++(-\bx,\by);
           \draw[<-](0,0)--++(0,-.8)node[below]{$(1-\bu_k^\oplus )U+\dfrac{\bu_k}2+\bv_{k-1}$};

           \foreach \i/\tex/\sty in {{-\bxt,-\byt}/A_1/below,{\bxt,-\byt}/A_2/below, {\bxt,\byt}/A_3/right, {-\bxt,\byt}/A_4/below right}
                \draw[fill=black]($.5*(\i)+1.5*(-\xt,\yt)$) circle(.015)node[\sty]{$\tex$};
        \end{tikzpicture}
        \caption{\label{figproj} The projection of $H^{TS}\cap (\R^{d-1}\times[y_{k}-\ep,y_k+\ep])$ on the plane $\pi(y_k)$. If $P=(\bx,y)\notin H^{TS}$ with $y\in(y_k,y_k+\ep)$, then $\bx=\cP_h(P)$ is not in the rectangle $A_1A_2A_3A_4$, which is just $\hat\bu_{k,t}U+2^{-1}(t-1)\bu_k+\bv_k$. This figure is drawn with $\bu_k=(-0.4,0.3)$ and $t=0.4$, $d=3$.}
    \end{figure}

    In the following three lemmas we establish some properties of the interior, exterior, and vertical boundary of $H$.

    Recall that $\myint H^{TS}=f_{TS}(\myint H).$
    \begin{lem}[Interior]\label{leminterior}
        Let $f$ be given as in \eqref{eqhomf}. For $P\in \myint H$, we have
            $$d_h\big(f(P),f(\pt_v H)\big) \ge (1-c)d_h(P,\pt_v H).$$
    \end{lem}
    \begin{proof}
        We first claim that
            \begin{equation}\label{eqdhfR1}
            d_h\big(f_F(P),f_F(\pt_v H^{TS})\big) \ge d_h(P,\pt_v H^{TS}), \text{ for }P\in \myint H^{TS}.
            \end{equation}
        Suppose $P=(\bx,y)\in \myint H^{TS}$ and let $y':=\omega_d\big(f_F(P)\big)$. We prove the claim for different values of $y$. If $y=0$ or $b$, \eqref{eqdhfR1} holds since $f_F$ is an identity on $\pi(0)\cup \pi(b)$. Now, let $y=y_k$ for some $k\in\{1,2,\dots,r-1\}$. Since $P\in S_k^0$, we have (see Figure \ref{figproj})
            \begin{equation}\label{eqbxin}
            \bx\in \brho(y_k)U+\bx(y_k) =(1-\bu_k^\oplus )U+\dfrac{\bu_k}2+\bv_{k-1}=(U+\bv_{k-1})\cap (U+\bv_k).
            \end{equation}
       The fact
            $$\cP_h\Big(f_F(H^{TS})\cap \pi(y_k)\Big)=(U+\bv_{k-1})\cup (U+\bv_k), $$
        together with \eqref{eqbxin}, implies that
            \begin{eqnarray*}
            d_h\big(f_F(P),f_F(\pt_v H^{TS})\big)&=&d_h\Big(\bx,\pt \big((U+\bv_{k-1})\cup (U+\bv_k)\big)\Big)\\
            &=&d_h\Big(\bx,\pt \big((U+\bv_{k-1})\cap (U+\bv_k)\big)\Big)\\
            &=&d_h\big(\bx,\brho(y_k)\pt U+\bx(y_k)\big),
            \end{eqnarray*}
        where we have used the expressions for $\brho(y)$ and $\bx(y)$ (see \eqref{eqbxy} and \eqref{eqrhoy}) and the fact that $f_F$ is an identity on $S_k^0$ (see Figures \ref{figprojection} and \ref{figproj}). Now, the conclusion follows from
            $$d_h(P,\pt_v H^{TS})=d_h\big(\bx,\brho(y_k)\pt U+\bx(y_k)\big).$$
        If $y\notin I_\ep$, by \eqref{eqbxy} and \eqref{eqrhoy} again we know
            $\brho(y)U=U$ and $\bx(y)=\bv_k,$
        where $k$ is the integer satisfying $y_k<y<y_{k+1}$. This implies
            $$d_h(P,\pt_v H^{TS})=d_h\big(f_F(P),f_F(\pt_vH^{TS})\big)$$
        by \eqref{eqHbiTS}, \eqref{eqHbiRTS} and the fact $y'\in(y_k,y_{k+1})$.

       Next, we assume $y\in I_{\ep,k}$ for some $k\in\{1,\ldots,r-1\}$ and $y\ne y_k$. By  symmetry, we assume $y_k<y<y_k+\ep$. The definition of $f_F$ yields $y'\in(y_k,y_k+\ep)$, which implies
            $$d_h\big(f_F(P),f_F(\pt_v H^{TS})\big)=d_h(\bx,\pt U+\bv_k).$$
        Setting $t=\ep^{-1}(y-y_k)\in [0,1]$, we see that
            $$(t\bu_k^\oplus +1-\bu_k^\oplus )U+\dfrac{t- 1}2\bu_k
            \subset\big(t\bu_k^\oplus +1-\bu_k^\oplus +(1-t)\bu_k^\oplus \big)U=U.$$
        Adding $\bv_k$ to both sides, we get $\brho(y) U+\bx(y)\subset U+\bv_k$, which leads to \eqref{eqdhfR1} since
            $$d_h(P,\pt_v H^{TS})=d_h\big(\bx,\brho(y)\pt U+\bx(y)\big).$$
        Now, \eqref{eqdhfR1}, together with Lemmas \ref{lemfBT} and \ref{lemfS}, yields
        \begin{eqnarray*}
          d_h\big(f(P),f(\pt_vH)\big)&\ge&d_h\big(f_{TS}(P),f_{TS}(\pt_vH)\big)\\
          &=&d_h\big(f_S(P),f_S(\pt_vH)\big)\\
          &\ge&(1-c)d_h(P,\pt_vH),
        \end{eqnarray*}
       which completes the proof.
    \end{proof}

    \begin{lem}[Exterior]\label{lemexterior}
        Let $f$ be given as in \eqref{eqhomf}. Let $P\in\R^{d-1}\times[0,b]$. If $P\notin H$, then
            $$d_h\big(f(P),f(\pt_v H)\big) \ge (1-c)^{1/\eta+1}d_h(P,\pt_v H).$$
    \end{lem}
    \begin{proof}
        We first prove the following inequality
        \begin{equation}\label{eqdhfF2}
        d_h\big(f_F(P),f_F(\pt_v H^{TS})\big) \ge (1-c)^{1/\eta}d_h(P,\pt_v H^{TS}),\quad\text{for } P\notin H^{TS}.
        \end{equation}
      Since $f_F$ is an identity on $\pi(0)$ and $\pi(b)$, we assume $P=(\bx,y)\notin H^{TS}$ satisfying $y\in[y_k,y_{k+1})$ for $0\le k<r$. If $y\notin I_\ep$, we see that $\brho(y)U+\bx(y)=U+\bv_k$. By \eqref{eqHbiTS} and \eqref{eqHbiRTS}, we know
            \begin{eqnarray*}
            d_h\big(f_F(P),f_F(\pt_v H^{TS})\big)&=&d_h(\bx,U+\bv_k)=d_h\big(\bx,\brho(y)U+\bx(y)\big)\\ &=&d_h(P,\pt_v H^{TS}).
            \end{eqnarray*}
        Next we assume $y\in I_{\ep}$, which implies $y\in I_{\ep,k}$ if $y\le y_k+\ep$ with $k>0$ and $y\in I_{\ep,k+1}$ if $y>y_{k+1}-\ep$ with $y<r$.  By symmetry, we consider the case $y_k\le y\le y_k+\ep$ and $k>0$. From \eqref{eqSk}, we know that
            $$S_{u_k}=S_k-\Big(\dfrac12\bu_k+\bv_{k-1},y_k\Big),$$
        and
            $$S^+_{u_k}=S_k^+-\Big(\dfrac12\bu_k+\bv_{k-1},y_k\Big),\quad S'^+_{u_k}=S_k'^+-\Big(\dfrac12\bu_k+\bv_{k-1},y_k\Big).$$
        Denote $P'=P-(2^{-1}\bu_k+\bv_{k-1},y_k)$. Since $P\notin H^{TS}$ and $y_k\le y\le y_k+\ep$, we know $\cP_h(P')\notin \hat\bu_{k,t}U+2^{-1}t\bu_k$, where $t=y'/\ep=(y-y_k)/\ep$ and $\hat\bu_{k,t}=t\bu_k^\oplus +1-\bu_k^\oplus $ (see Figure \ref{figproj}). Thus, we denote $\cP_h(P')$ by $\bx'+2^{-1}t\bu_k$. Then
            $$\bx'+\dfrac t2\bu_k=\bx-\dfrac12\bu_k-\bv_{k-1}.$$
        So the horizontal distance on the right side of \eqref{eqdhfF2} is
            \begin{eqnarray*}
            d_h(P,\pt_v H^{TS})&=&d_h(P,S_k^+\cup S_k'^+)
            =d_h(P',S^+_{\bu_k}\cup S'^+_{\bu_k})
            \\&=&d_h\Big(\bx'+\dfrac t2\bu_k,\hat\bu_{k,t}\pt U+\dfrac t2\bu_k\Big).\\
            &=&d_h(\bx',\hat\bu_{k,t}\pt U).
            \end{eqnarray*}%
       If $P$ lies below $S_k$, which means that $P'$ lies below $S_{\bu_k}$, we see that $f_F(P)$ lies below $\pi(y_k)$. Hence \setcolsep%
            \begin{eqnarray*}
            d_h\big(f_F(P),f_F(\pt H_v^{TS})\big)&=&
            d_h(\bx,U+\bv_{k-1})\\
            &=&d_h\Big(\bx-\bv_{k-1}-\dfrac {1+t}2\bu_k,U-\dfrac {1+t}2\bu_k\Big).\\
            &=&d_h\Big(\bx',U-\dfrac {1+t}2\bu_k\Big).
            \end{eqnarray*}%
      If $P$ lies above $S_k$ or $P$ is in $S_k$, which means $P'$ lies above $S_{\bu_k}$ or $P'$ is in $S_{\bu_k}$, we know $f_F(P)$ lies above $\pi(y_k)$ or $f_F(P)$ is in $\pi(y_k)$. Thus,
        so\setcolsep
            \begin{eqnarray*}
            d_h\big(f_F(P),f_F(\pt H^{TS})\big)&=&
            d_h(\bx,U+\bv_{k})\\
            &=&d_h\Big(\bx-\bv_{k-1}-\dfrac {1+t}2\bu_k,U+\dfrac {1-t}2\bu_k\Big)\\
            &=&d_h\Big(\bx',U+\dfrac {1-t}2\bu_k\Big).
            \end{eqnarray*}
        Applying Lemma \ref{lemex}, we get the inequality \eqref{eqdhfF2}. Now, the result follows from Lemmas \ref{lemfBT} and \ref{lemfS} and \eqref{eqdhfF2}.
    \end{proof}
    \begin{lem}[Boundary]\label{lembound}
        Let $0^{d-1}\ne\bi\in\{-1,0,1\}^{d-1}$. If $P\in \myint H_\bi$ and the height of $f(P)$ is not in $\{y_0,\ldots,y_r\}$, then
            $$d_h\big(f(P),f(\pt_v H_\bi)\big)\ge (1-c)d_h(P,\pt_v H_\bi).$$
    \end{lem}
    \begin{proof}
        Let $P=(\bx,y)\in\myint H_\bi$, $P'=f_{TS}(P):=(\bx',y)$ and denote by $y'$ the height of $f(P)$.
        Then $\bx'=\brho(y)\bx+\bx(y)$. Since $f_F$ is an identity on $\cup_{k=0}^rS_k^0$, Lemma \ref{lemfBT} and \eqref{eqfSPHbi} in Lemma \ref{lemfS} imply $y'=y_k$ if $y=y_k$ for any $k$. So, by the assumption that the height of $f(P)$ is not in $\{y_0,\ldots,y_r\}$, we assume $y\in(y_k,y_{k+1})$ for some $0<k<r$, which implies $y'\in(y_k,y_{k+1})$. If $y\notin I_\ep$, we know that $\bx(y)=\bv_k$ and $\brho(y)=\bone_{d-1}$ by \eqref{eqbxy} and \eqref{eqrhoy}. Hence, the restriction of $f(H_\bi)$ on $\pi(y')$ is $(U_\bi+\bv_k)\times\{y'\}$ by \eqref{eqHbiRTS}. Thus
            \begin{eqnarray*}
            d_h\big(f(P),f(\pt H_\bi)\big)&=&d_h(\bx',\pt U_\bi+\bv_k)
            \\ &=&d_h\big(\brho(y)\bx+\bx(y),\pt U_\bi+\bx(y)\big)\\
            &=&d_h(\bx,\pt U_\bi)=d_h(P,\pt H_\bi).
            \end{eqnarray*}
        If $y\in I_\ep$, we may assume, in  view of symmetry, that $y\in(y_k,y_k+\ep)$. For such $y$, we claim that $\brho(y)U_\bi+\bx(y)\subset U_\bi+\bv_k$. If the claim holds, then
            \begin{eqnarray*}
            d_h\big(f(P),f(\pt_v H_\bi)\big)&=&d_h\big(\brho(y)\bx+\bx(y),\pt U_\bi+\bv_k\big)\\
            &\ge&d_h\big(\brho(y)\bx+\bx(y),\brho(y)\pt U_\bi+\bx(y)\big)\\
            &=&\brho(y)d_h(\bx,\pt U_\bi)\ge(1-c)d_h(P,\pt_v H_\bi).
            \end{eqnarray*}
        The first inequality above holds since $\brho(y)\bx+\bx(y)\in\rho(y)U_\bi^\circ+\bx(y)$ (see Figure \ref{figptUbi}). Thus, the conclusion holds.
        \begin{figure}
             \pgfmathsetmacro{\h}{1}
             \pgfmathsetmacro{\t}{.22}
             \pgfmathsetmacro{\eps}{.5}
             \pgfmathsetmacro{\ux}{.35}
             \pgfmathsetmacro{\uy}{.3}
             \pgfmathsetmacro{\uz}{.4}
            \begin{tikzpicture}[scale=3,>=stealth]

             \draw[fill=gray!40](-.5,0)rectangle++(1,\eps);
             \path[fill=white](.5,\eps)--++(0,-\eps)--++(-\ux,0)--cycle;
             \draw[red](.5,\eps)--++(-\ux,-\eps);
             \def\mybraceht{.07}
             \mybrace{-.5,\t*\eps}{(1-\ux+\t*\ux)}{$\brho(y)U_\bi+\bx(y)$}
             \def\mybraceht{.05}
             \mybrace[1]{-.5,\t*\eps}1{$U_\bi+\bv_k$};
             \draw(-.5+\h+\ux,-.5+\uy)rectangle++(1,1);
             \draw[fill=gray!40](-.5+\h+\ux,-.5+\uy)rectangle++ (1-\ux+\t*\ux,1-\uy+\t*\uy);
             \foreach \i/\tex in {0/(a)\,$d=3$,1/(b)\,$d=4$}
                \draw(\i*\h+\i*.45-.2,-.35)node[below]{\tex};
             \foreach \i/\tex/\loc in {-.5/P_1=P_1'/left, .5-\ux+\t*\ux/P_2/below right, .5/P_2'/right,-.3/{}/right}
                \draw[fill=black](\i,\t*\eps)circle(.02)node[\loc,font=\tiny]{$\tex$};
             \draw[<-,red](-.3,1.5*\t*\eps)--++(0,\eps)node[above]{$P$};
            \draw[red](-.5+\ux+\h,.5+\t*\uy)-|(.5+\t*\ux+\h,-.5+\uy);
            \foreach \i/\tex/\loc in {{-.5+\ux,-.5+\uy}/P_1=P_1'/below, {.5+\t*\ux,-.5+\uy}/P_2/below,{-.5+\ux,.5+\t*\uy}/P_4/left, {.5+\t*\ux,.5+\t*\uy}/P_3/right,{.5+\ux,-.5+\uy}/P_2'/below, {.5+\ux,.5+\uy}/P_3'/right,{-.5+\ux,.5+\uy}/P_4'/left}
                \draw[fill=black]($(\i)+(\h,0)$)circle(.02) node[\loc,font=\tiny]{$\tex$};
            \foreach \i/\tex in {{.2,.2}/\brho(y)U_\bi+\bx(y), {.4,.65}/U_i+\bv_k}
                \draw($(\i)+(\h,0)$)node[font=\tiny]{$\tex$};
            \draw[fill=black](\h+.2,0)circle(.02)node[font=\tiny,above,text=red]{$P$};
        \end{tikzpicture}\,
        \begin{tikzpicture}[scale=2.8,x={(-.35cm,-.25cm)}, y={(1cm,0cm)},z={(0cm,1cm)}]
        \pgfmathsetmacro{\rx}{(1-\ux)};
        \pgfmathsetmacro{\ry}{(1-\uy)};
        \pgfmathsetmacro{\rz}{(1-\uz)};
        \pgfmathsetmacro{\rtx}{(1-\ux+\t*\ux)};
        \pgfmathsetmacro{\rty}{(1-\uy+\t*\uy)};
        \pgfmathsetmacro{\rtz}{(1-\uz+\t*\uz)};

        \coordinate (a) at (.5-.5*\ux,-.5+.5*\uy,-.5+.5*\uz); 
        \coordinate (b) at (.5-.5*\ux,-.5+.5*\uy-\t*\uy,-.5+.5*\uz); 
        \coordinate (c) at (.5-.5*\ux,-.5-.5*\uy,-.5+.5*\uz);

        \foreach \i/\j/\jj/\jjj/\col/\color in {
            c/1/1/1/gray!30/gray!45,b/\rtx/\rty/\rtz/red!35/red!50}
            {\path[fill=\col]($(\i)+(-\j,0,\jjj)$)--++(\j,0,0) --++(0,0,-\jjj)--++(-\j,0,0)--cycle;
            \path[fill=\col]($(\i)+(-\j,0,\jjj)$)--++(0,\jj,0) --++(0,0,-\jjj) --++(0,-\jj,0)--cycle;            \path[fill=\color]($(\i)+(-\j,0,\jjj)$)--++(\j,0,0) --++(0,\jj,0) --++(-\j,0,0)--cycle;}
        \foreach \i/\tex/\loc/\text/\loct in {%
            {0,0,0}/P_4/below/Q_4/left, {0,\rty,0}/P_1=P_1'/below/Q_1/right, {-\rtx,\rty,0}/P_2/below right/Q_2/above right, {-\rtx,0,0}/P_3/above/Q_3/above} {
                \draw[fill=black]($(b)+(\i)$)circle(.02) node[font=\tiny,\loc]{$\tex$};
                \draw[fill=black]($(b)+(\i)+(0,0,\rtz)$)circle(.02) node[font=\tiny,\loct]{$\text$};
                }
        \foreach \i/\tex/\loc/\text/\loct in {{0,0,0}/P_4'/below left/Q_4'/left, {-1,1,0}/P'_2/right/Q'_2/right, {-1,0,0}/P'_3/above left/Q'_3/left,{0,1,0}/{}/right/Q_1'/above}
            {\draw[fill=black]($(c)+(\i)$)circle(.02) node[font=\tiny,\loc]{$\tex$};
            \draw[fill=black]($(c)+(\i)+(0,0,1)$)circle(.02) node[font=\tiny,\loct]{$\text$};}

        \foreach \i/\j/\jj/\jjj in {c/1/1/1, b/\rtx/\rty/\rtz}
        \foreach \k in {{0,0,0},{0,\jj,0},{-\j,\jj,0}}
            {\draw[thin]($(\i)+(\k)$)--++(0,0,\jjj);
            \draw[dashed,very thin]($(\i)+(-\j,0,0)$)--++(0,0,\jjj);}

        \draw[very thin,dashed](c)--++(-1,0,0)--++(0,1,0);
        \draw[thin](c)--++(0,1,0)--++(-1,0,0);
        \draw[thin]($(c)+(0,0,1)$)--++(-1,0,0)--++(0,1,0)--++(1,0,0)--cycle;
        \foreach \i/\sty in {0/dashed,\rtz/dashed}
            {\draw($(b)+(0,0,\i)$)--++(0,\rty,0)--++(-\rtx,0,0);
            \draw[very thin,\sty]($(b)+(0,0,\i)$)--++(-\rtx,0,0)--++(0,\rty,0);}
        \draw($(c)+(0,.6,-.17)$)node[below]{(c) $d=5$};
        \foreach \i/\tex in {1.15/U_\bi+\bv_k,.1/\brho(y)U_\bi+\bx(y)}
            \draw($(c)+(0,.6,\i)$)node[font=\tiny]{$\tex$};
        \coordinate (p) at (0,0,0);\
        \draw[fill=black](p)circle(.02)node[above,text=red,font=\tiny]{$P$};

        \end{tikzpicture}
        \caption{\label{figptUbi} An illustration of Lemma \ref{lembound}. Figure (a) is the intersection of $H_\bi^{TS}$ and $(U_\bi+\bv_k)\times[y_k,y_k+\ep]$ for the case $d=3$. Figures (b) and (c) are the intersections of $H_\bi^{TS}$ and $(U_\bi+\bv_k)\times\{y\}$, where $(b)$ is for $d=4$ and (c) is for $d=5$. The line $P_1P_2$ in (a), the rectangle $P_1P_2P_3P_4$ in (b), the cuboid $P_1P_2P_3P_4$-$Q_1Q_2Q_3Q_4$ in (c) are $\brho(y)U_\bi+\bx(y)$, while the line $P'_1P'_2$ in (a), the rectangle $P'_1P'_2P'_3P'_4$ in (b) and the cuboid $P_1'P_2'P'_3P'_4$-$Q'_1Q'_2Q'_3Q'_4$ in (c) are $U_\bi+\bv_k$.
        The figures are drawn with $y=y_k+0.22\ep$ and $\bu=(u_1,0.35)$, $\bi=10$ in (a), $\bu=(u_1,0.35,0.3)$, $\bi=100$ in (b) and $\bu=(u_1,0.35,0.3,0.4)$, $\bi=1000$ in (c), where $u_1$ is a constant in $(0,1)$.
        }
        \end{figure}
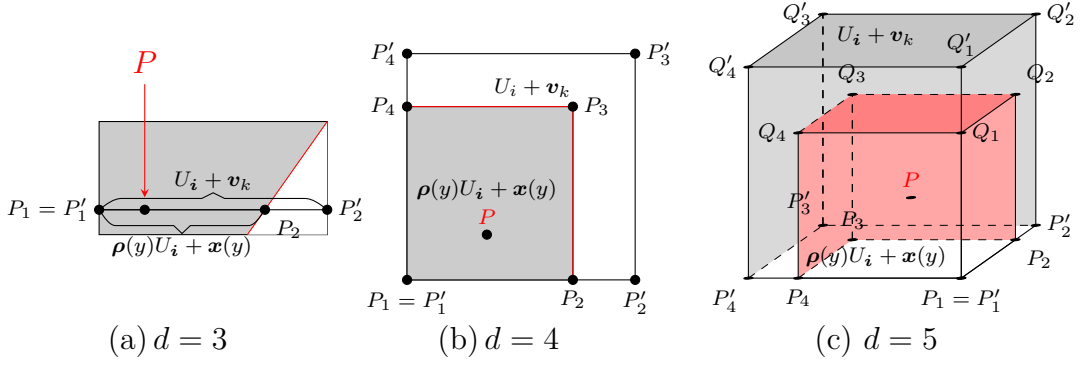

        Next, we prove the claim by considering the components of horizontal coordinates of $\brho(y)U_\bi+\bx(y)$ and $\bu_i+\bv_k$. Let $u_{kj}:=\omega_j(\bu_k)$ and let
            $$\Lambda_1=\{j:\,i_ju_{kj}<0\},\quad \Lambda_2=\{j:\,i_j\ne0,u_{kj}=0\},\quad \Lambda_3=\{j:i_j=0\}.$$
        Since $P'\notin S^+_k$ (using the given condition $\cP_v\big(f(P)\big)\notin\{y_0,\ldots,y_r\}$), we know that $\Lambda_1,\Lambda_2,\Lambda_3$ form a partition of $\{1,\ldots,d-1\}$. So, for each $\by\in U_\bi$, we have
            $$2\omega_j (\by)=-\sign(u_{kj}),\quad 2\omega_j(\by)=i_j\quad \text{or}\quad 2\omega_j(\by)\in[-1,1],$$
        according to $j\in\Lambda_1$, $j\in\Lambda_2$ or $j\in\Lambda_3$.
        It follows that
            $$\omega_j\big(\brho(y)U_\bi+\bx(y)\big) =\{2^{-1}i_j+\omega_j(\bv_k)\}=\omega_j(U_\bi+\bv_k)$$
        or
            $$\omega_j\big(\brho(y)U_\bi+\bx(y)\big)\subset [-1/2,1/2]+\omega_j(\bv_k)=
            \omega_j(U_\bi+\bv_k),$$
        according to $j\in\Lambda_1\cup\Lambda_2$ or $j\in\Lambda_3$.
    This proves the claim and completes the proof of the lemma.
    \end{proof}

    For $P,P'\in\R^d$ and $E\subset \R^d$, let
        $$d_v(P,P'):=\big|\cP_v(P)-\cP_v(P')\big|\quad\text{and}\quad
        d_v(P,E):=\inf\{d_v(P,P'):\,P'\in E\}.$$
    A word $\bi=i_1\ldots i_{d-1}$ in $\{-1,0,1\}^{d-1}$ will determine a vector $(i_1,\ldots,i_{d-1})\in\R^{d-1}$. We use $\bi$ to denote such a vector if there is no confusion. The following lemma concerns the horizontal distance and the height difference of two points and their images under $f$. 

    \begin{lem}[Distance between two points]\label{lemdistance}
        Let $P,P'\in\R^{d-1}\times[0,b]$. Suppose $c<4^{-1}$.
\begin{enumerate}
\item[(1)]\textup{(Not in $\pt_vH$)} If the horizontal distance between $P$ and $\pt_vH$ satisfies
            $$d_h(P,\pt_vH)\ge\begin{cases}
                \dfrac12(1+2\eta+4c)^2-\dfrac12,&P\notin H;\\
                2c,&P\in H,
            \end{cases}$$
        then $f(P)$ has the same height as $P$. If $d_h(P',\pt_vH)$ also satisfies the inequality stated above, then
            \begin{equation}\label{eqdfPP'1}
                d_h\big(f(P),f(P')\big)\ge (1-c)^{1/\eta}d_h(P,P').
            \end{equation}
 \item[(2)]\textup{(In $\pt_vH$)} Let $0^{d-1}\ne\bi\in\{-1,0,1\}^{d-1}$ and $P,P'\in H_\bi$. Assume there exists $k$, $0\le k<r$, such that the heights of $f(P),f(P')$ are both in the interval $(y_k,y_{k+1})$.
        If none of the components of $\bi$ is zero, then
            \begin{equation}\label{eqdiffht}
                d_v\big(f(P),f(P')\big)\ge d_v(P,P').
            \end{equation}
        If some component of $\bi$ is zero and $P,P'\in\myint H_\bi$ satisfy
            $d_h(P, \pt_vH_\bi)>2c$ and $d_h(P', \pt_vH_\bi)>2c,$
        then \eqref{eqdfPP'1} and \eqref{eqdiffht} also hold. In particular, if $d_v\big(f(P),\pi(0)\cup\pi(b)\big)<b/r$, then
            $$d_v\big(f(P),\pi(0)\cup\pi(b)\big)\ge d_v(P,\pi(0)\cup\pi(b)).$$
\end{enumerate}
%
%
%
    \end{lem}

    \begin{proof}
        Let $P=(\bx,y)$ and $P'=(\bx',y')\in\R^{d-1}\times[0,b]$.

        (1) We first consider the case $P\notin H$. This implies $\bx\in(1+2\eta+4c)^2U^c$. As $0<c<1/4$ and $\eta>0$,  we have $(1-c)(1+2\eta+4c)>1+3c-4c^2\ge1+2c>1$.
        For $1\le k<r$ and $y\in[y_k,y_{k+1}]$, we have
            \begin{gather*}
            \cP_h\Big(f_{TS}\big((1+2\eta+4c)^2U^c \times\{y_k\}\big)\Big)
            =\brho(y)(1+2\eta+4c)^2U^c+\bx(y)\\
            \begin{array}{rcl}
            &\subset&\dfrac{\brho(y)}{1-c}(1+2\eta+4c)U^c+\bx(y)\\
            &\subset&(1+2\eta+4c)U^c+\bx(y)\\
            &\subset&(1+2\eta+4c-2\|\bx(y)\|_\infty)U^c
            \subset(1+2\eta+2c)U^c.
            \end{array}
            \end{gather*}
        On the other hand $\varphi(\bz,y_k)=y_k$ for $0\le k\le r$ if $\bz\in \cP_h(\hat S_k^0)$, where
            $$\hat S_n^0=\Big((1+2\eta+\|\bu_n\|_\infty)U^c+\dfrac12\bu_n+\bv_{n-1}\Big) \times\{y_n\},\quad 0<n<r,$$
        and $\hat S_0^0=\hat S_1^0$, $\hat S_r^0=\hat S_{r-1}^0$. This implies,  by the definitions of $\varphi$ and $f_F$, that $f_F$ is an identity on $\{\bz\}\times[y_k,y_{k+1}]$ when $\bz\in\cap_{n=0}^r\cP_h\big(\hat S_k^0\big)$. For $1\le k<r$,
            \begin{equation}\label{eqcPSk0tmp}
            (1+2\eta+\|\bu_k\|_\infty)U+\dfrac{\bu_k}2+\bv_{k-1}
            \subset(1+2\eta+2c)U,
            \end{equation}
        which implies
            $$\bigcap_{k=0}^r\cP_h\big(\hat S_k^0)\big) \supset (1+2\eta+2c)U^c.$$
    Consequently, $f$ does not change the height of each point in $(1+2\eta+4c)^2U^c\times[0,b]$.

        Next we consider the case $P\in H$. The given hypothesis implies $\bx\in (1-4c)U$. For each $1\le k<r$,
            \begin{equation}
\brho(y_k)(1-4c)U+\bx(y_k)
            \subset(1-4c)U+\dfrac12\bu_k+\bv_{k-1}
            \subset(1-2c)U.\label{eqsub}
            \end{equation}
        On the other hand, from Lemma \ref{lemSbuheight}(2) and \eqref{eqSk}, we see that
            \begin{eqnarray}
            \nonumber
            \bigcap_{n=0}^r\cP_h\big(S_n^0\big) &=&\bigcap_{n=1}^{r-1}\big((1-\bu_n^+)U+\dfrac12\bu_n+\bv_{n-1}\Big)
            \\ &\supset&\bigcap_{n=1}^{r-1}\Big(1-2\sum_{m=1}^n\|\bu_m\|_\infty\Big)U
            \supset(1-2c)U.\label{eqsup}
            \end{eqnarray}
        \eqref{eqsub} and \eqref{eqsup} imply that $f$ does not change the height of $P\in(1-4c)U\times[0,b]$.

        Now we show \eqref{eqdfPP'1}. The above discussion says that the heights of $f(P)$ and $f(P')$ are $y$ and $y'$ respectively. If $y\ne y'$, the left side of \eqref{eqdfPP'1} is $+\infty$, and thus \eqref{eqdfPP'1} holds trivially. If $y=y'$, a straightforward computation yields
            \begin{eqnarray*}
            d_h\big(f(P),f(P')\big)
              &=&d_h\big(\brho(y)\bx+\bx(y),\brho(y)\bx'+\bx(y)\big)\\
              &\ge& \brho(y)d_h\big(\bx,\bx')\ge(1-c)^{1/\eta}d_h(P,P'),
            \end{eqnarray*}
        since $0<\eta<1$ and $0\le c<1$.

        \begin{figure}
        \newcommand{\myrectangle}[3][]
            {\draw[#1](#2)--++(#3,0,0)--++(0,#3,0)--++(-#3,0,0)--cycle}
        \newcommand{\mypathrectangle}[3][]
            {\path[#1](#2)--++(#3,0,0)--++(0,#3,0)--++(-#3,0,0)--cycle}
        \newcommand{\myhzline}[4][]
            {\draw[#1](#2)--++(0,#3,0)--++(-#4,0,0);}
        \newcommand{\myzhline}[4][]
            {\draw[#1](#2)--++(-#4,0,0)--++(0,#3,0);}
        \begin{tikzpicture}[{>=stealth,scale=2.65, declare function={f(\t)=(1-\rho)*\t/\h+\rho;}, x={(-.35cm,-.25cm)}, y={(1cm,0cm)},z={(0cm,1cm)}}]
        \pgfmathsetmacro{\h}{1.75}
        \pgfmathsetmacro{\eps}{.2}
        \pgfmathsetmacro{\v}{1}
        \pgfmathsetmacro{\vv}{(\v-\eps)}
        \pgfmathsetmacro{\t}{.6}
        \pgfmathsetmacro{\rho}{.7}      
        \pgfmathsetmacro{\rhot}{\rho+\t-\rho*\t}      
        \pgfmathsetmacro{\u}{(1-\rho)}  
        \coordinate (a) at (-.5*\rho,-.5*\rho,0);  
        \coordinate (y) at (-\u,-\u,\eps); 
        \coordinate (t) at ($(a)+(-\u,-\u,\v)$); 
        \foreach \i in {0,1}
            {\myrectangle[fill=red!30,opacity=.5]{$(a)+(-\u,2*\i*\h-\u,0)$}1;
            \myrectangle[fill=blue!30,opacity=.5]{$(t)+(0,\i*\h,0)$}1;}
        \foreach \i/\j/\jj/\tex in {0/\rho/\rho/{(a)\, $0<k<r-1$}, 1/1/\rho/{(b)\,$k=0$}, 2/\rho/1/{(c) $k=r-1$}}
            {\myrectangle[fill=gray!50]{$(a)+(\rho-\j,\rho+\i*\h-\j,0)$}\j;
            \myrectangle[fill=gray!50]{$(t)+(0,\i*\h,0)$}\jj;
            \draw(-.5*\u,\i*\h-.5*\u,-.25)node[below]{\tex};
            }
        \foreach \i/\j/\colo/\colt in {0/0/blue/red,1/2/blue/red}
            {\myhzline[\colo]{$(t)+(1,\j*\h,\eps-\v)$}11;
            \myzhline[\colt,dashed]{$(t)+(1,\j*\h,\eps-\v)$}11;
            \myhzline[\colo]{$(t)+(1,\i*\h,\vv-\v)$}11;
            \myzhline[\colt,dashed]{$(t)+(1,\i*\h,\vv-\v)$}11;
            }
        \foreach \i/\len in {0/\rho,1/1,2/\rho}
            {
            \draw[very thin, dashed]($(t)+(0,\i*\h,0)$)--++(0,0,-\v);
            \draw[fill=black,thin](-.5*\u,-.5*\u+\i*\h,0)circle(.02) node[right,font=\tiny]{$O_b$}
            --++(0,0,\v)circle(.02)node[left,font=\tiny]{$O_t$};
            \foreach \j in {{1,0,0},{1,1,0},{0,1,0}}
                \draw[very thin]($(\j)+(t)+(0,\i*\h,0)$)--++(0,0,-\v);
            \foreach \ii/\jj in {{0,0,0}/3,{-\len,0,0}/2,{0,-\len,0}/4, {-\len,-\len,0}/1}
            \draw[fill=black]($(a)+(\rho,\rho+\i*\h,0)+(\ii)$)circle(.02) node[font=\tiny,below]{$A_\jj$};}
        \foreach \i/\len in {0/\rho,1/\rho,2/1}
            \foreach \ii/\jj in {{0,0,0}/1,{0,\len,0}/2,{\len,\len,0}/3, {\len,0,0}/4}
            \draw[fill=black]($(t)+(0,\i*\h,0)+(\ii)$)circle(.02) node[font=\tiny,above]{$B_\jj$};
        \foreach \i in {0,1}
            {
            \draw[fill=blue!30,opacity=.5]($(t)+(\rho,\rho+\i*\h,0)$) --++(-\rho,0,0)--++(0,\u,-\eps)--++(1,0,0)--cycle;
            \draw[fill=blue!30,opacity=.5]($(t)+(\rho,\rho+\i*\h,0)$) --++(0,-\rho,0)--++(\u,0,-\eps)--++(0,1,0)--cycle;
            \draw[fill=red!30,opacity=.5]($(a)+(0,2*\i*\h,0)$) --++(\rho,0,0)-- ++(0,-\u,\eps)--++(-1,0,0)--cycle;            \draw[fill=red!30,opacity=.5]($(a)+(0,2*\i*\h,0)$) --++(0,\rho,0)--++(-\u,0,\eps)--++(0,-1,0)--cycle;
            \draw[<-]($(t)+(.5,\i*\h+1,-.1)$)--++(0,.1,-.15) node[font=\tiny,below,yshift=3pt]{$S^-_{k\!+\!1}$};            \draw[<-]($(a)+(.2,2*\i*\h+.45,\eps)$)--++(-.5,-.3,.05) node[font=\tiny,above,yshift=-3pt]{$S^+_{k}$};
            \draw[fill=black](-.5*\u,-.5*\u+\i*\h,\v-\eps)circle(.02) node[right,font=\tiny]{$O_2$};
            \draw[fill=black](-.5*\u,-.5*\u+2*\i*\h,\eps)circle(.02) node[right,font=\tiny]{$O_1$};
            \foreach \j/\sty/\tex in {1/-8pt/P,2/8pt/Q}
                \foreach \ii/\num/\loc in {{0,0,0}/1/-5pt,{0,1,0}/2/5pt, {1,1,0}/3/5pt, {1,0,0}/4/-5pt}
                    \draw[fill=black]($(t)+(\ii) +(0,\j*\i*\h,{(\j-2)*\eps+(\j-1)*(\eps-\v)})$)circle(.02) node[yshift=\sty,xshift=\loc,text=black, font=\tiny]{$\tex_\num$};
            }
        \end{tikzpicture}
    \caption{\label{fig3I}
        The restrictions of $H^{TS}$ and $f_F(H^{TS})$ on $\R^{d-1}\times[y_k,y_{k+1}]$, where (a) is for the case $k\ne0,r-1$; (b) and (c) are  for the cases $k=0$ and $k=r-1$ respectively.
        The line segment $O_bO_t$ is in the line $\ell_{\bv_k}$, which is the center axis of $(U+\bv_k)\times[y_k,y_{k+1}]$. The ranges of heights of $O_bO_1,O_1O_2$ and $O_2O_t$ in (a)--(c) are $I^{(b)},I^{(m)}$ and $I^{(t)}$ respectively.}
    \end{figure}
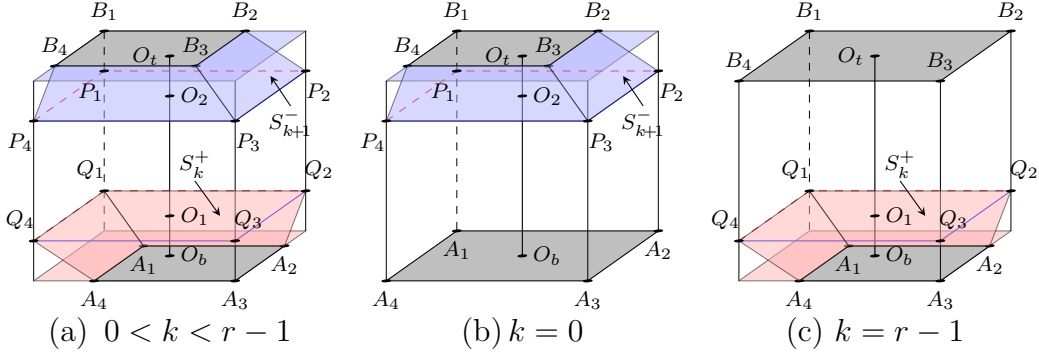

         (2) The condition $\bi\ne0^{d-1}$ implies that $P\in\pt_vH$. We may assume $y\in[y_k,y_{k+1}]$ for some $k$. Let $P_1=f_{TS}(P)=(\bx_1,y)$ and let
            $$I^{(b)}=\begin{cases}\emptyset,&k=0;\\ [y_k,y_k+\ep),&k>0,\end{cases}\quad
            I^{(t)}=\begin{cases}\emptyset,& k=r-1; \\ (y_{k+1}-\ep,y_{k+1}],&k<r-1,\end{cases}$$
        and $I^{(m)}=[y_k+\ep-\delta_k(0)\ep, y_{k+1}-\ep+\delta_{k+1}(r)\ep]$, where $\delta_x(\cdot)$ is the Dirac function on $\R$; see Figure \ref{fig3I}.
        From the given condition on the height of $f(P)$, we see that
            \begin{equation}\label{eqby-bvk}
            \bx_1-\bv_k\in
            \begin{cases}
            (\pt U_{\bu_k}^-\cup U_{\bu_k}^0)\backslash (\pt U_{\bu_k}^+),
            & y\in I^{(b)};\\
            (\pt U_{\bu_{k+1}}^+\cup \pt U_{\bu_{k+1}}^0)\backslash (\pt U_{\bk_{k+1}}^-), &y\in I^{(t)};\\
            \pt U,&y\in I^{(m)}.
            \end{cases}
            \end{equation}

Let $t_1=\ep^{-1}(y-y_k)$ and $t_2=\ep^{-1}(y_{k+1}-y)$. For $j\in\Lambda:=\{j:\,i_j\ne0\}$, it follows from $\bx_1=\brho(y)\bx+\bx(y)$, $\bx\in U_\bi$ and \eqref{eqby-bvk} that
            \begin{eqnarray*}
            \lefteqn{\omega_j\big(2\bx_1-2\bv_k\big)} \\
            &=&
            \begin{cases}
            \omega_j\big((t_1\bu_k^\oplus +1-\bu_k^\oplus )i_j+(t_1-1)\bu_k\big), &y\in I^{(b)};\\
            \omega_j\big((t_2\bu_{k+1}^\oplus+1-\bu_{k+1}^\oplus)i_j-(t_2-1)\bu_{k+1}\big),
            &y\in I^{(t)};\\
            \omega_j\big(2(\bx+\bv_k)-2\bv_k\big),&y\in I^{(m)}
            \end{cases}\\
            &=&i_j=\pm1.
            \end{eqnarray*}
      If $\bi$ has some a component, the hypothesis $d_h(P_1,\pt_vH_\bi)>2c$ yields $|\omega_j(2\bx)|\le1-4c$ for  $j\notin\Lambda$, which implies that
            \begin{eqnarray*}
                \big|\omega_j(2\bx_1-2\bv_k)\big| &=&\big|\omega_j\big(2\brho(y)\bx +2\bx(y)-2\bv_k\big)\big|\\
                &\le&\|\brho(y)\|_\infty\omega_j(2\bx^{\oplus}) +2\big(\omega_j(\bu_k^\oplus )+\omega_j(\bu_{k+1}^\oplus)\big)
                \\&\le&1-4c+2c\\
                &\le& \min\big\{1-2\|\bu_k\|_\infty,1-2\|\bu_{k+1}\|_\infty\big\},
            \end{eqnarray*}
        where we set $\bu_0=\bu_{r}=\bv_0$.
        Recall that
            $$S_n=S_{\bu_n}+\dfrac12\bu_n+\bv_{n-1}
            =S_{\bu_n}-\dfrac12\bu_n+\bv_n,\quad 1\le n<r.$$
        Denote by $y_k',y_{k+1}'$ the heights of the intersections of $\ell_{\bx_1}$ and $S_k, S_{k+1}$ respectively. Applying Lemma \ref{lemSbuheight}(3,4) to such $\ell_\bx$ and $S_k, S_{k+1}$, we have
           \begin{equation}\label{eqcPvQ}
              y_k'=
              \begin{cases}
                 y_k,&\text{if }\Lambda_1=\emptyset;\\
                 y_k+\ep,&\text{if }\Lambda_1\ne\emptyset,
              \end{cases}\quad
              y_{k+1}'=
              \begin{cases}
                 y_{k+1},&\text{if }\Lambda_2=\emptyset;\\
                 y_{k+1}-\ep,&\text{if }\Lambda_2\ne\emptyset,
              \end{cases}
           \end{equation}
        where $\Lambda_1=\{j:\,\omega_j(\bu_k)i_j>0\}$ and $\Lambda_2=\{j:\,\omega_j(\bu_{k+1})i_j<0\}$.
        So if $\bi$ has a zero component, the height of $f(P)$, given by\setcolsep
            \begin{equation}\label{eqcPvfP}
            \cP_v\big(f(P)\big)=\dfrac{b}r\cdot\dfrac{y-y_k'}{y_{k+1}'-y_k'}+y_k
            \in [y_k,y_{k+1}]
            \end{equation}
        is determined only by $y$, $\bi$, $\bu_k$ and $\bu_{k+1}$, and is independent of $\bx$ (i.e., independent of the horizontal coordinates of $Q,Q'$). Assume the height of $f(P')$ is given as in \eqref{eqcPvfP} with $y$ replaced by $y'$. Then
            \begin{eqnarray*}
              d_v\big(f(P),f(P')\big) &=&
              \Big|\dfrac{b}r\cdot\dfrac{y-y_k'}{y_{k+1}'-y_k'} -\dfrac{b}r\cdot\dfrac{y'-y_k'}{y_{k+1}'-y_k'}\Big|\\
              &=&\Big|\dfrac{b/r}{y_{k+1}'-y_k'}(y-y')\Big|\\
              &\ge& |y-y'|=d_v(P,P')
            \end{eqnarray*}
        since $|y_{k+1}'-y_k'|\le b/r$ by \eqref{eqcPvQ} and \eqref{eqep}. We prove \eqref{eqdiffht} for case $\bi$ has a zero component and the case $\bi$ does not have any zero component. Now, \eqref{eqdfPP'1} follows from Lemmas \ref{lemfBT} and \ref{lemfS}. To prove the rest of the conclusion, we let $P'\in S_0^0$ or $P'\in S^0_r$ according to whether the height of $f(P)$ is less than $y_1$ or larger than $y_{r-1}$ and then apply \eqref{eqdiffht} for such $P$ and $P'$. This completes the proof.
    \end{proof}

\begin{figure}
    \begin{tikzpicture}[xscale=.95,yscale=1.8,>=stealth, declare function={f(\x)=\x*\tone;g(\x)=-\x*\x*\eps/\et/\et+\eps;
    invf(\x)=\x*\ttwo;}]
      \pgfmathsetmacro{\v}{-2}    
      \pgfmathsetmacro{\h}{4.5}    
      \pgfmathsetmacro{\eps}{.1}   
      \pgfmathsetmacro{\et}{.3}    
      \pgfmathsetmacro{\ett}{.2}   
      \pgfmathsetmacro{\uo}{-.2}    
      \pgfmathsetmacro{\ut}{-.15}   
      \pgfmathsetmacro{\vt}{(\uo+\ut)} 
      \pgfmathsetmacro{\tone}{.2}   
      \pgfmathsetmacro{\ttwo}{\tone-\vt}
      \pgfmathsetmacro{\r}{0}       

      \foreach \i/\tex in {1/\pi(0),2/\pi(b)}
        {
        \draw({-\et-.5-\ett+f(\i-1)},\i-1)node[font=\tiny,left]{$\tex$} --++(1+2*\ett+2*\et,0);
        \draw(\i-1.5,0)--++({f(1)},1);
        \draw[fill=red!30]({-.5+f(\i/3-\eps)},\i/3-\eps) --++(1,0)--++({f(2*\eps)},2*\eps)--++(-1,0) --cycle;
        \draw[domain=-\et:0]({-\et-.5-\ett+f(\i/3)},\i/3) node[font=\tiny,left]{$\widetilde{f_{TS}^{-1}(S_\i)}$}  --++(.2,0) plot({\x-.5+g(\x)*\tone+f(\i/3)},{g(\x)+\i/3}) --++({-f(\eps)},-\eps)--++(1,0)--++({-f(\eps)},-\eps) plot({-(-\et-\x)+.5-g(-\et-\x)*\tone+f(\i/3)},{-g(-\et-\x)+\i/3}) --++(\ett,0);}

      \foreach \i/\tex in {1/\pi(0),2/\pi(b)}
        {
        \draw({-\et-.5-\ett+\h},\i-1)node[font=\tiny,left]{$\tex$} --++(1+2*\ett+2*\et,0);
        \draw(\i-1.5+\h,0)--++(0,1);
        \draw[fill=red!30](-.5+\h,\i/3-\eps)rectangle++(1,2*\eps);
        \draw[domain=-\et:0](-\et-.5-\ett+\h,\i/3) node[font=\tiny,left]{$f_{TS}^{-1}(S_\i)$} --++(.2,0) plot({\x-.5+\h},{g(\x)+\i/3}) --++(0,-\eps)--++(1,0)--++(0,-\eps) plot({-(-\et-\x)+.5+\h},{-g(-\et-\x)+\i/3})--++(\ett,0);}

      \foreach \i/\j/\tex in {1/\uo/\pi(0),2/\ut/\pi(b)}
        {
        \draw({-\et-.5-\ett+2*\h},\i-1)node[font=\tiny,left]{$\tex$} --++(1+2*\ett+2*\et,0);
        \draw(\i-1.5+2*\h,0)--++(0,1);
        \draw[white](2*\h+.5,1/3-\eps)--++(0,2*\eps);
        \foreach \ii in {-1,1}
            \draw[white,thick,fill=white](\ii*.5+2*\h,\i/3-\eps)  node[font=\tiny,left]{$f_S(S_\i)$} --++(0,2*\eps);
        \draw[fill=red!30](-.5+2*\h,\i/3-\eps)--++(-\j*.5,\eps)--++(\j*.5,\eps) --++(1,0)--++(\j*.5,-\eps)--++(-\j*.5,-\eps)--cycle;
        \draw[domain=-\et:0](-\et-.5-\ett+2*\h,\i/3)  node[font=\tiny,left]{$f_T^{-1}(S_\i)$}--++(.2,0) plot({\x-.5+2*\h},{g(\x)+\i/3}) --++(-\j*.5,-\eps)--++(1+\j,0)--++(-.5*\j,-\eps) plot({-(-\et-\x)+.5+2*\h},{-g(-\et-\x)+\i/3})--++(\ett,0);}

      \foreach \i/\j/\k/\tex in {1/\uo/\uo/\pi(0),2/\ut/\vt/\pi(b)}
        {
        \draw({-\et-.5-\ett+1.3*(\i-1)*\vt+(1-\r)*2*\h},\i-1+\v) node[font=\tiny,left]{$\tex$}--++(1.2+2*\ett+2*\et,0);
        \draw({\i-1.5+(1-\r)*2*\h},\v)--++(0,{1/3+(1-\i)*\eps}) --++(\uo,\eps)--++(0,1/3-\eps)--++(\ut,\eps) --++(0,{1/3-(2-\i)*\eps});
        \draw[fill=red!30]({-.5+(\i-1)*\uo+(1-\r)*2*\h},\i/3+\v)--++(\j,\eps)--++(1,0) --++(0,-\eps)--++(-\j,-\eps)--++(-1,0)--cycle;
        \draw[domain=-\et:0]({-\et-.5-\ett+\k+(1-\r)*2*\h},\i/3+\v)  node[font=\tiny,left]{$S_\i$}--++(.2,0) plot({\x-.5+\k+(1-\r)*2*\h},{g(\x)+\i/3+\v}) --++(-\j,-\eps)--++(1+\j,0)--++(-\j,-\eps) plot({-(-\et-\x)+.5+(\i-1)*\uo+(1-\r)*2*\h},{-g(-\et-\x)+\i/3+\v}) --++(\ett,0);
      }

      \foreach \i/\j/\k/\tex in {1/\uo/\uo/\pi(0),2/\ut/\vt/\pi(b)}
        {
        \draw({-\et-.5-\ett+1.3*(\i-1)*\vt+\h},\i-1+\v)  node[font=\tiny,left]{$\tex$} --++(1.2+2*\ett+2*\et,0);
        \draw(\i-1.5+\h,\v)--++(0,1/3) --++(\uo,0)--++(0,1/3)--++(\ut,0) --++(0,1/3);
        \draw[fill=red!30]({-.5+(\i-1)*\uo+\h},\i/3+\v)--++(\j,0) --++(-\j,\eps)--++(1+\j,0) --++(0,-\eps)--++(-\j,0)--++(\j,-\eps)--++(-1-\j,0)--cycle;
        \draw[domain=-\et:0]({-\et-.5-\ett+\h+\k},\i/3+\v) node[font=\tiny,left]{$f_F(S_\i)$} --++(1-\k+2*\ett+2*\et,0);
      }
      \foreach \i/\j/\k/\tex in {1/\uo/\uo/\pi(0),2/\ut/\vt/\pi(b)}
        {
        \draw({-\et-.5-\ett+f(\i-1)+\r*2*\h},\i-1+\v)  node[font=\tiny,left]{$\tex$} --++(1+2*\ett+2*\et,0);
        \draw({\i-1.5+\r*2*\h},\v)--++({invf(1/3)},1/3)--++(\uo,0)--++({invf(1/3)},1/3) --++(\ut,0)--++({invf(1/3)},1/3);
        \draw[fill=red!30]({-.5+(\i-1)*\uo+\r*2*\h+invf(\i/3)},\i/3+\v) --++(\j,0)--++({-\j+invf(\eps)},\eps)--++(1+\j,0) --++({-invf(\eps)},-\eps)--++(-\j,0)--++({-invf(\eps)+\j},-\eps) --++({-1-\j},0)--cycle;
        \draw[domain=-\et:0]({-\et-.5-\ett+\r*2*\h+f(\i/3)},\i/3+\v)  node[font=\tiny,left]{$f_{\bt'}\!\circ\! f_F(S_\i)$} --++(1+2*\ett+2*\et,0);
      }
      \foreach \i/\tex in {1/f_\bt,2/f_S}
        \draw[->](\i*\h-\h+1,.5)--node[above]{$\tex$}++(1.2,0);
      \foreach \i/\tex in {1/f_F,2/f_{\bt'}}
        \draw[->]({(2-\i)*.2+2.2+2*\h-\i*\h-1.2*\r-3*\h*\r+2*\h*\i*\r}, .5+\v)--node[above]{$\tex$}++({-1.1+2.6*\r},0);
      \foreach \i/\tex in {1/h,2/f_T}
        \draw[->]({(\i-1)*2*\h},-.4)--node[right]{$\tex$} ++({(1.5-\i)*4*\r*\h},\v+1.6);
      \foreach \i/\tex in {0/\widetilde H,1/H,2/f_S(H)}
        \draw(\i*\h,-.25)node{$\tex$};
      \foreach \i/\tex in {0/f_{TS}(H),1/f(H),2/h(\widetilde H)}
        \draw({(2-\i)*\h},-.33+\v)node{$\tex$};
        \foreach \j in {0,\v}
        \foreach \i in {1,2}
            {\draw[fill=black](\i-1.5,\j)circle(.02) node[font=\tiny,below]{$P_\i$};
            \draw[fill=black](\i-1.5+\tone,\j+1)circle(.02) node[font=\tiny,above]{$Q_\i$};}
    \end{tikzpicture}\vspace*{-10pt}
    \caption{\label{fighom2} An illustration of $h$, where $\bt'=\bt-\bv_{r-1}$.}
    \end{figure}
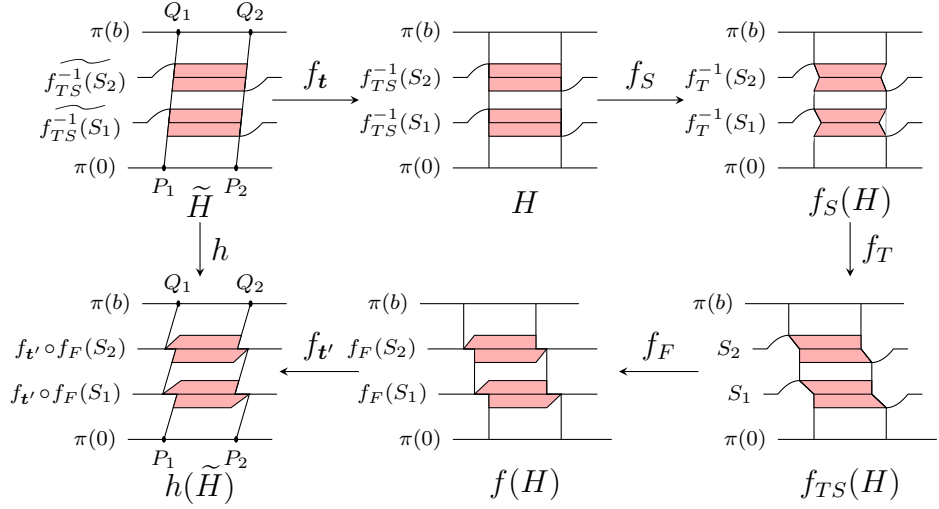

    Define a homeomorphism  $h$ on $\R^d$ as follows (see Figure \ref{fighom2}):
        \begin{equation}\label{eqhomh}
        h(P):=\begin{cases}
          f_{\bt-\sum\bu_{k}}\circ f\circ f_{-\bt}(P),&\cP_v(P)\in[0,b];\\
          P,&\cP_v(P)\notin(0,b).
        \end{cases}
        \end{equation}
    We will denote $h$ by $h_{\bt,\{\bu_k\}}$ when we want to emphasize that $h$ is determined by the vectors $\bt,\bu_1,\ldots,\bu_{r-1}$.

    For a set $E\subset\R^d$,  we denote  $f_{\bt}(E)$ by $\widetilde E$, e.g.,
        \begin{equation}\label{eqwtHbi}
        \begin{array}{rcl}
        \widetilde H_\bi=f_{\bt}(H_\bi)
        &=&\Big\{\Big(\dfrac{y\bt}b,y\Big):\,0\le y\le b\Big\}+U_\bi\times\{0\}\\[5pt]
        &=&\big\{t(\bt,b):\,0\le t\le 1\big\}+U_\bi\times\{0\},
        \end{array}
        \end{equation}
        $\myint \widetilde H=f_{\bt}(\myint H)$, $\pt_v\widetilde H_\bi=f_{\bt}(\pt_v H_\bi)$, and so on. 
        Noting that $0<\eta<1$, $f_\bt^{-1}=f_{-\bt}$, 
        we get the following useful proposition.
    \begin{prop}\label{theoh}
        Assume $\bt$, $\bu_1,\ldots,\bu_{r-1}$ are $(d-1)$-dimensional vectors satisfying $\{\bu_k\}\subset2U^\circ$ and $\eta\in(\max\{\|\bu_k\|_\infty\},1)$. Then the function $h$, defined as in \eqref{eqhomh} with respect to the  vectors $\bt$, $\bu_1,\ldots,\bu_{r-1}$, is a homeomorphism  on $\R^d$. Moreover, if $c=\sum\|\bu_k\|_\infty<4^{-1}$, then the following statements hold.

\begin{enumerate}
\item[(1)] \indent\textup{(\textbf{General property})} The mapping $h$ is an identity on $\R^{d-1}\times\big((-\infty,0]\cup [b,+\infty)\big)$ satisfying
            \begin{equation}\label{eqh(H)}
            h(\widetilde H)=\bigcup_{k=0}^{r-1} \bigcup_{y_k\le y\le y_{k+1}} \bigg(U+\sum_{m=1}^{k}\bu_m +
            \dfrac yb\Big(\bt-\sum_{m=1}^{r-1}\bu_{m}\Big)\bigg)\times\{y\}
            \end{equation}
        (where we set an empty sum to be zero).
        \item[(2)]        \indent\textup{(\textbf{Not on boundary})} For $P\notin\pt_v\widetilde H$, we have
            $$d_h\big(h(P),h(\pt_v \widetilde H)\big) \ge (1-c)^{2/\eta}d_h(P,\pt_v \widetilde H).$$
        \item[(3)]            \indent\textup{(\textbf{Boundary})} Assume $0^{d-1}\ne\bi\in\{-1,0,1\}^{d-1}$.
        If $P\in \myint\widetilde H_\bi$ is such that the height of $f(P)$ is not in $\{y_0,\ldots,y_r\}$, then
            $$d_h\big(h(P),h(\pt_v\widetilde H_\bi)\big)\ge (1-c)^{2/\eta}d_h(P,\pt_v \widetilde H_\bi).$$
\item[(4)] \indent\textup{(\textbf{Distance between points})} Assume $P,P'\in\R^{d-1}\times[0,b]$.\\
\item[(i)]  \textup{(\textbf{Not in $\pt_v\widetilde H$})} If the horizontal distance between $P$ and $\pt_v\widetilde H$ satisfies
            $$d_h(P,\pt_v\widetilde H)\ge\begin{cases}
                \dfrac12(1+2\eta+4c)^2-\dfrac12,&P\notin\widetilde H;\\
                2c,&P\in \widetilde H,
            \end{cases}$$
        then $h(P)$ has the same height as $P$. If $d_h(P',\pt_v\widetilde H)$ also satisfies the inequality stated above, then
            \begin{equation}\label{eqdhPP'2}
                d_h\big(h(P),h(P')\big)\ge (1-c)^{2/\eta}d_h(P,P').
            \end{equation}
    \item[(ii)]  \textup{(\textbf{In $\pt_v\widetilde H$})}  Let $0^{d-1}\ne\bi\in\{-1,0,1\}^{d-1}$ and $P,P'\in\widetilde H_\bi$. Assume there exists $k$, $0\le k<r$, such that the heights of both $h(P)$ and $h(P')$ are in the interval $(y_k,y_{k+1})$. If none of the components of $\bi$ is zero, then
            \begin{equation}\label{eqdifhht}
                d_v\big(h(P),h(P')\big)\ge d_v(P,P').
            \end{equation}
        If some component of $\bi$ is zero and $P,P'\in\myint\widetilde H_\bi$ satisfy
            $$d_h(P, \pt_v\widetilde H_\bi),d_h(P', \pt_v\widetilde H_\bi)>2c,$$
        then \eqref{eqdhPP'2} and \eqref{eqdifhht} also hold. In particular, if $d_v\big(h(P),\pi(0)\cup\pi(b)\big)<b/r$, then
            $$d_v\big(h(P),\pi(0)\cup\pi(b)\big)\ge d_v(P,\pi(0)\cup\pi(b)).$$
\end{enumerate}
    \end{prop}
    \begin{proof}
        That the mapping $h$ is a homeomorphism  on $\R^d$ follows directly from our construction. Conclusion (1) is a  consequence of Lemma \ref{lemHbi} and the definition of $h$. Conclusion (2) is a corollary of Lemmas \ref{lemfBT}, \ref{leminterior}, and \ref{lemexterior}. Conclusion (3)
        follows from Lemmas \ref{lemfBT} and \ref{lembound}. Lemmas \ref{lemfBT} and \ref{lemdistance} yield conclusions (4)(i) and (4)(ii). The proof is complete.
    \end{proof}

\section{Preparation for the construction
    \label{S:sec4}}
    In this section, we introduce some notation concerning the self-affine tile $T(A,\D)$, some of which come from \cite{Rn}. Let $(A,\D)$ be the self-affine pair as in \eqref{eqexprad}. For convenience, we also use $r$ to denote $p_d$. Let $\cL_{d-1}:=\prod_{j=1}^{d-1}[0,|p_j|)\cap \Z^d$. For each dight $(\biota,k)\in\D$ in Theorem \ref{theomain1} we add a perturbation $\ba_k=(a^{(1)}_{k},\ldots,a^{(d-1)}_{k})\in\R^{d-1}$ to the horizontal coordinate $\biota$. Using a translation, we assume $\ba_0=(0,\ldots,0)$, which implies that the origin of $\R^d$ is in the new digit set. The new digit set can be expressed in a more general form:
    \begin{equation}\label{digit}
        \widetilde{\D}=\big\{(\biota+\ba_{k},k):\, 0\le k<|r|,\biota\in\cL_{d-1}\big\}.
    \end{equation}
Theorem \ref{theomain1} can be adjusted to a general form as follows.
    \begin{theo}\label{theo1}
        For $d\ge 2$, let $A$ and $\widetilde{\D}$ be as in \eqref{eqexprad} and \eqref{digit} respectively,  and let $T=T(A,\widetilde{\D})$ be the self-affine tile.
    Then the following statements are equivalent.
    \begin{enumerate}
        \item $T^\circ$ is connected;
        \item $T$ is a $d$-dimensional tame ball;
        \item the following inequality holds:
        \begin{equation}\label{eqnonstandard}
            \max_{0\le k<|r|-1\atop 1\le j\le d-1}
            \left\{
                \left|
                    \dfrac{a_{k+1}^{(j)}-a_{k}^{(j)}}{p_j}-\dfrac{\sign(r)\cdot a^{(j)}_{|r|-1}+s_j} {p_j(p_j-\sign(r))}
                \right|
            \right\}
            <1.
       \end{equation}
    \end{enumerate}
    \end{theo}

   The following lemma, proved in \cite{Rn}, explains why we use this perturbation strategy.%
    \begin{lem}\label{csa}
        Let $A$ and $\widetilde{\D}$ be as in \eqref{eqexprad} and \eqref{digit} respectively. Then there exists another self-affine pair $(A^2,\D')$ such that $T(A^2,\D')$ is a translation of $T(A,\widetilde{\D})$ and the condition
        \eqref{eqnonstandard} for $(A,\widetilde{\D})$ is equivalent to its analogue for $(A^2,\D')$.
    \end{lem}
    \begin{rema}
    From the above lemma, we assume throughout the rest of this paper that the matrix $A$ in the self-affine pair $(A,\widetilde{\D})$ satisfies \eqref{eqexprad} with positive eigenvalues,  the digit set $\widetilde{\D}$ has the form \eqref{digit}, and $(0,\ldots,0)\in\widetilde{\D}$.
    \end{rema}

    Recall that the multiplication of vectors in $\R^d$ is defined as in \eqref{eqmultiply}, i.e.,
        $$\bx\boxdot\by=\bx\by:=\big(\omega_1(\bx)\omega_1(\by),\ldots, \omega_{n}(\bx)\omega_{n}(\by)\big),\quad \bx,\by\in\R^d.$$
 From the form of the matrix $A$, we see that
    \begin{equation*}\label{a}
       A^{-n}=\begin{pmatrix}
           A_1^{-n}&\bs\boxdot\bs_n\\
           0&r^{-n}
              \end{pmatrix},\quad n\ge 1,
    \end{equation*}
    where $\bs_n=(s^{(1)}_{n},\ldots,s^{(d-1)}_{n})$ and
    \begin{gather}\label{st}
        s^{(j)}_{n}:=\lim_{q\to p_j}\dfrac{q^{-n}-r^{-n}}{r-q}=
        \begin{cases}
            \dfrac{p_j^{-n}-r^{-n}}{r-p_j}, & p_j\ne r, \smallskip\\
            \dfrac{n}{r^{n+1}},                 & p_j=r.
        \end{cases}
    \end{gather}
    For $\bk=(k_1,\ldots,k_{|\bk|})\in \Sigma_r^*\cup\Sigma_r^{\mathbb N}$, let
    \begin{equation}
        \ba(\bk):=\big(a_1(\bk),\ldots,a_{d-1}(\bk)\big),\quad
        \bs(\bk):=\big(s_1(\bk),\ldots,s_{d-1}(\bk)\big),
    \end{equation}
    where
    \begin{gather}\label{eqajsj}
        a_j(\bk):=\sum_n\frac{a^{(j)}_{k_n}}{p_j^n},\quad s_j(\bk):=\sum_n s^{(j)}_{n}k_n, \quad 1\le j<d,
        \end{gather}
    and let
    $$
    \bb(\bk):=\ba(\bk)+\bs\boxdot \bs(\bk).
    $$
    Here the empty sum is defined to be zero. Hence $\bb(\emptyset)=\ba(\emptyset)=\bs(\emptyset)=(0,\ldots,0)\in\R^{d-1}$. Let
        \begin{equation}\label{eqIFS}
            S_{\biota,k}(P):=A^{-1}\big(P+(\biota+\ba_k,k)\big), \quad  \biota\in\cL_{d-1},\quad  0\le k<r,
        \end{equation}
    be the iterated function system associated with $(A,\widetilde{\D})$.  Then it follows from the definition of $T$ that  
        \begin{eqnarray}
            \label{eqTIFS}
        T&=&\bigcup_{\biota\in \cL_{d-1},0\le k\le r-1}S_{\biota,k}(T)\\
         &=&\Big\{\Big(\bb(\bk)+ \big(\varphi_{p_1}(\bi_1),\ldots,\varphi_{p_{d-1}}(\bi_{d-1})\big), \varphi_r(\bk)\Big):\, \bi_j\in{\Sigma_{p_j}^{\mathbb N}},\ \bk\in{\Sigma_r^{\mathbb N}}\Big\} \nonumber\\
         &=&\Big\{\big(\bb(\bk),\varphi_r(\bk)\big):\,\bk\in\Sigma_r^{\mathbb N}\big\}
         +[0,1]^{d-1}\times\{0\}.\label{eqtrd}
        \end{eqnarray}
Define
            $$S_{\biota_1\cdots \biota_n,k_1\cdots k_n}:=S_{\biota_1,k_1}\circ \cdots \circ S_{\biota_n,k_n},\quad n\ge1.$$

    We now introduce the symbol $T_\bk$ which will help us simplify the proof of Theorem \ref{theo1}. Recall that the convex combination of two sets is defined in Subsection \ref{hyperplane}, or more precisely in \eqref{eqsetconvex}. For $n\ge1$ and $\bk\in\Sigma_r^{n-1}$, let $T_\bk$ be the union of the convex combinations of the top and bottom surfaces of the following iterates of $T$:
        $$ S_{\biota_1\cdots \biota_n,\bk}(T),\quad \biota_i\in\cL_{d-1},\ 1\le i\le n.$$
    In other words, $T_\bk$ is the union of the convex combination of two hypercubes, namely,
        $$\big([0,1]^{d-1}+\bb(\bk)\big)\times \big\{\varphi_r(\bk)\big\}$$
        \text{ and }
        $$\big([0,1]^{d-1}+\bb(\bk\overline{r-1})\big)\times\Big\{ \varphi_r(\bk)+\dfrac1{r^{n-1}}\Big\}.$$
    Hence (see Figure \ref{figTbk}(a) for $d=3$ and (c) for $d=2$), \setcolsep
        \begin{eqnarray}\nonumber
        \lefteqn{T_\bk=\Big\{
            \Big((1-t)\bb(\bk\bar0)\!+\!t\bb(\bk\overline{r-1}),\dfrac t{r^{n-1}}\Big):\,0\!\le\! t\!\le\! 1\Big\} \!+\![0,1]^{d-1}\!\times\!\{\varphi_r(\bk)\}}\\
            &=&\Big\{t\Big(\bb(0^{n-1}\overline{r-1}), \dfrac1{r^{n-1}}\Big):\, 0\!\le\! t\!\le\! 1\Big\} \!+\![0,1]^{d-1}\!\times\!\{0\}\!+\!\big(\bb(\bk),\varphi_r(\bk)\big) \label{eqTbk}
            \\&=&\Big\{\big(\bb(\bk),\varphi_r(\bk)\big)\!+\! t\Big(\bb(0^{n-1}\overline{r-1}), \dfrac1{r^{n-1}}\Big):\, 0\!\le\! t\!\le\! 1\Big\} \!+\![0,1]^{d-1}\!\times\!\{0\} \label{eqTbkT}
        \end{eqnarray}
    Comparing \eqref{eqtrd} and \eqref{eqTbkT}, the expressions for $T$ and $T_\bk$, we obtain the following  lemma.
    \begin{lem}\label{lemTk}
        For $n\ge1$ and $\bk\in\Sigma_r^{n-1}$, let $T_\bk$ is defined as above. Then
            $$T=\lim_{n\to\infty}\bigcup_{\bk\in\Sigma^n_r} T_{\bk}$$
        in the Hausdorff distance.
    \end{lem}
    \begin{proof}
        Note that $T$ is the attractor of the iterated function system $\{S_{\biota,k}\}$ defined as in \eqref{eqIFS}. It follows that
            $$T=\lim_{n\to\infty}\bigcup_{\mathop{\biota_i\in\cL_{d-1}, 0\le k_i<r,} \atop 1\le i\le n}S_{\biota_1\cdots \biota_n,k_1\cdots k_n}(E)$$
        in the Hausdorff distance\cite{Falc}, where $E$ is any compact subset of $\R^d$.
        In particular, the above limit also exists if we set $E=T_\emptyset$. On the other hand, by iterating $n$ times, one obtains, for $\bk\in\Sigma_r^n$,
            $$\bigcup_{\mathop{\biota_i\in\cL_{d-1},}\atop 1\le i\le n} S_{\biota_1\cdots \biota_n,\bk}(T_\emptyset)
            =T_{\bk}.$$
The conclusion follows.
    \end{proof}

    We will employ the homeomorphism  defined in \eqref{eqhomh} to construct
    a desired homeomorphism  on $\R^d$. For this purpose, we need some additional notation, some of which has been introduced in \cite{R3,Rn}. For a set $E\subset\R^d$, let $\|E\|_\infty:=\sup\{\|P\|_\infty:\,P\in E\}$. Denote $a:=\max\{\|\ba_k\|_\infty:\,1\le k<r\}$, $p:=\min\{p_j:\,1\le j\le d\}$ and $s:=\|\bs\|_\infty$, where $\bs$ is given as in \eqref{eqexprad}. For $n\ge1$, let
        \begin{gather}
        \ep_n:=\dfrac1{3r^{2n}},\quad
            c_n:=\sum_{k=1}^{r-1}\|\bu_{n,k}\|_\infty,\label{eqcn}\\
        \bt_n:=\bb(0^{n-1}\overline{r-1}),\label{eqbtn}\\
        Y_n:=\Big\{\dfrac k{r^{n-1}}:\,0\le k\le r^{n-1}\Big\} =\{\varphi_r(\bk):\, \bk\in\Sigma_r^{n-1}\}\cup\{1\},\label{eqYn}
        \end{gather}
    where
        \begin{equation}\label{eqbunk}
        \bu_{n,k}:=\bb\big(0^{n-1}k\big) -\bb\big(0^{n-1}(k-1)\overline{r-1}\big),\quad\text{ for }0<k<r,
        \end{equation}
    and let
        \begin{equation}
        \eta_n:=\max\biggl\{\big\|\{\bu_{n,k}\}_{k=1}^r\big\|_\infty \Big(\dfrac12+\dfrac1{2\big\|\{\bu_{n,k}\}_{k=1}^r\big\|_\infty}\biggr), \dfrac1{n+1}\bigg\}. \label{eqetan}
        \end{equation}
    We point out here that the sequence $\{\eta_n\}$ can be replaced by any sequence satisfying $\|\{\bu_{1,k}\}\|_\infty<\eta_1<1$, and decreasing to $0$ slowly so that the homeomorphism  $h_{\bt_n,\{\bu_{n,k}\}}$ is well-defined.

    For $n\ge1$ and $\bk\in\Sigma_r^{n-1}$, define
        \begin{equation}\label{eqbwbk}
        \bw(\bk):=\big(2^{-1}\bone_{d-1}+\bb(\bk),\varphi_r(\bk)\big)
        \end{equation}
    which is the center of the hypercube
        $$T_\bk|_{\pi(\varphi_r(\bk))}=U\times\{0\}+\bw(\bk).$$
    We comment on the vectors $\bt_n$ and $\bu_{n,k}$ (see Figure \ref{figTbk}). For $n\ge1$, $0<k<r$, $\bu_{n,k}$ is the difference of the centers of two $(d-1)$-dimensional unit hypercubes in $T|_{\pi(\varphi(\bk k))}$ (where $\bk$ is any word in $\Sigma_r^{n-1}$), namely, $U+\bb(\bk k)$ and $U+\bb(\bk(k-1)\overline{r-1})$. In fact,
        $$\bb(\bk k)-\bb\big(\bk(k-1)\overline{r-1}\big) =\bb(0^{n-1}k)-\bb\big(0^{n-1}(k-1)\overline{r-1}\big)=\bu_{n,k}.$$
    For each $\bk\in\Sigma_r^{n-1}$, $\bt_n$ is the vector corresponding to an affine transformation defined as in \eqref{eqfB} with $b=r^{1-n}$ in a suitable rectangular coordinate system (moving the origin to $\bw(\bk)$), the inverse of which transforms $T_{\bk}$ into a vertical $d$-dimensional prism, namely,
        $$f_{-\bt_n}\big(T_{\bk}-\bw(\bk)\big) =U\times[0,r^{1-n}].$$
    Both $\bu_{n,k}$ and $\bt_n$ are determined by the self-affine pair $(A,\widetilde{\D})$ but are independent of $\bk\in\Sigma_r^{n-1}$.

    \def\myrect#1#2{\draw[fill=#2,opacity=.5](#1)--++(1,0,0)--++(0,1,0) --++(-1,0,0)--cycle}
    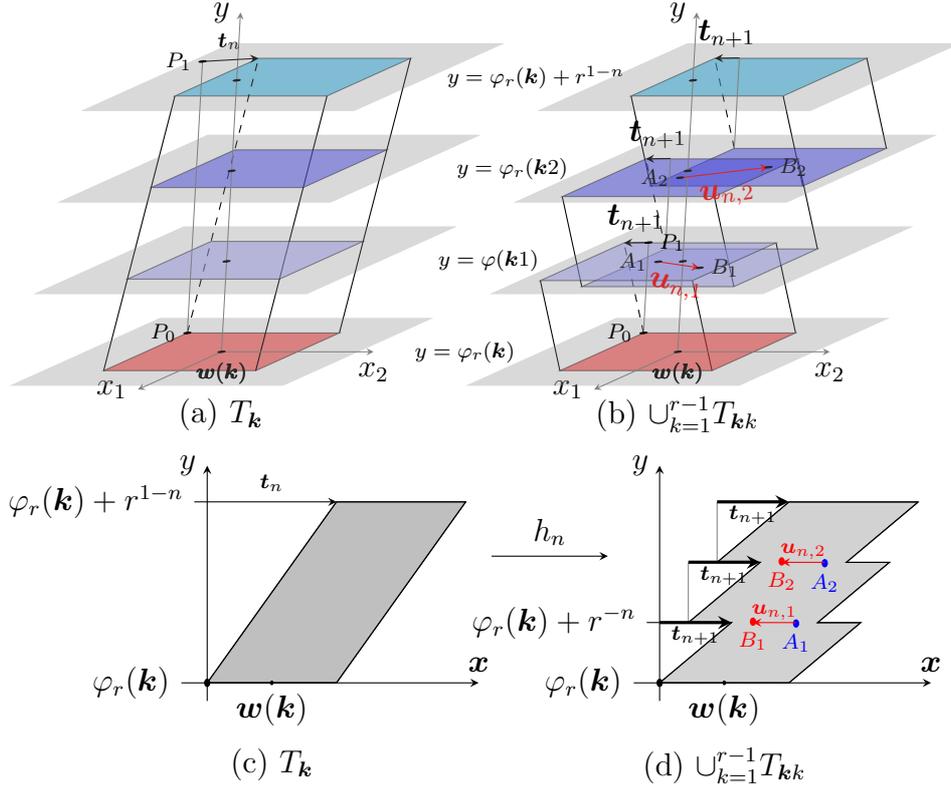
\begin{figure}
    \begin{tikzpicture}[{>=stealth,xscale=2,yscale=2,
        declare function={%
            xone=\xt*.5-\xo;
            xtwo=2*\xo-\xt*.5;
            x(\k)=xone*\k*\k+xtwo*\k;
            yone=\yt*.5-\yo;
            ytwo=2*\xo-\yt*.5;
            y(\k)=yone*\k*\k+ytwo*\k;
            bbx(\k)=x(\k)/\r+\so*\k/\r/\r;
            btx(\k)=bbx(\k)+x(2)/\r/(\r-1)+\so*(2*\r-1)/\r/\r/(\r-1);
            bby(\k)=y(\k)/\r+\st*\k/\r/\r;
            bty(\k)=bby(\k)+x(2)/\r/(\r-1)+\st*(2*\r-1)/\r/\r/(\r-1);},
            x={(-.55cm,-.25cm)}, y={(1cm,0cm)},z={(.1cm,1.8cm)}}]
        \pgfmathsetmacro{\r}{3}
        \pgfmathsetmacro{\h}{3}
        \pgfmathsetmacro{\so}{.35}   
        \pgfmathsetmacro{\st}{-.25}  
        \pgfmathsetmacro{\xz}{0}   
        \pgfmathsetmacro{\xo}{.4} 
        \pgfmathsetmacro{\xt}{-.52}
        \pgfmathsetmacro{\yz}{0}   
        \pgfmathsetmacro{\yo}{.1}
        \pgfmathsetmacro{\yt}{.4} 
        \pgfmathsetmacro{\eps}{.4}  

        \foreach \k in {0,1,2}
            {
            \foreach \i/\sty in {{0,0,0}/dashed,{1,0,0}/black, {0,1,0}/black,{1,1,0}/black}
            \draw[\sty]($({bbx(\k)},{bby(\k)},\k/3)+(\i)$)--++ ({btx(\k)-bbx(\k)},{bty(\k)-bby(\k)},1/3);
            }
        \foreach \i in {0,-\h}
            {\myrect{{bbx(0)},{bby(0)+\i},0}{red};
            \myrect{{btx(2)},{bty(2)+\i},1}{cyan};}

        \foreach \i/\sty in {{0,0,0}/dashed,{1,0,0}/black, {0,1,0}/black,{1,1,0}/black}
            \draw[\sty]($({bbx(0)},{bby(0)-\h},0)+(\i)$)--++ ({btx(2)-bbx(0)},{bty(2)-bby(0)},1);

        \foreach \k/\col in {1/blue!40,2/blue!80}
            {\foreach \i in {
                {{bbx(\k)},{bby(\k)},\k/3},
                {{btx(\k-1)},{bty(\k-1)},\k/3}}
                \myrect{\i}{\col};
            \myrect{{btx(2)*\k/3},{bty(2)*\k/3-\h},\k/3}{\col};}
        \foreach \j in {0,-\h}
            {\foreach \i/\loc/\tex in {{1,0,0}/left/x_1,{0,1,0}/below/x_2, {0,0,1.25}/left/y}
            \draw[help lines,->](.5,.5+\j,0) node[font=\tiny,below,text=black]{$\bw(\bk)$} --++(\i)node[\loc,black]{$\tex$};
            \foreach \k in {0,1,2,3}
                \draw[fill=black](.5,.5+\j,\k/3)circle(.02);
            }
        \draw[help lines](0,-\h,0)--++(0,0,1);
        \draw[->](0,-\h,1)--node[sloped,above,font=\tiny] {$\bt_n$}({btx(2)},{bty(2)-\h},1);
        \foreach \k in {0,1,2}
            {\draw[help lines]({bbx(\k)},{bby(\k)},\k/3)--++(0,0,1/3);
            \draw[->]({bbx(\k)},{bby(\k)},\k/3+1/3) --node[sloped,above]{$\bt_{n+1}$} ({btx(\k)},{bty(\k)},\k/3+1/3);}
        \foreach \k in {1,2}
            {\draw[<-,red]({bbx(\k)+.5},{bby(\k)+.5},\k/3)-- node[sloped,below]{$\bu_{n,\k}$} ({btx(\k-1)+.5},{bty(\k-1)+.5},\k/3);
            \foreach \kk/\sty/\ki in {{{bbx(\k)+.5},{bby(\k)+.5},\k/3}/right/B, {{btx(\k-1)+.5},{bty(\k-1)+.5},\k/3}/left/A}
                \draw[fill=black](\kk)circle(.02)node[\sty,font=\tiny] {$\ki_\k$};
            }

        \foreach \k/\tex in {0/\varphi_r(\bk),1/\varphi(\bk1), 2/\varphi_r(\bk2),3/\varphi_r(\bk)+r^{1-n}}
            {\foreach \i in {0,-\h}
                \path[fill=gray,opacity=.3]({btx(2)*\k/3-\eps},{bty(2)*\k/3-\eps+\i},\k/3) --++(1+2*\eps,0,0) --++(0,1+2*\eps,0)--++(-1-2*\eps,0,0)--cycle;
            \draw({btx(2)*\k/3+.5},{bty(2)*\k/3-1.3+\eps},\k/3) node[text=black,font=\tiny] {$y=\tex$};}
        \foreach \i in {0,-\h}
            \draw[fill=black](0,\i,0)circle(.02)node[font=\tiny,left]{$P_0$};
        \foreach \i/\j/\sty in {0/.333/right,-\h/1/left}
            \draw[fill=black](0,\i,\j)circle(.02)node[\sty,font=\tiny]{$P_1$};
        \foreach \i/\tex in {-\h/(a) $T_\bk$, 0/(b) $\cup_{k=1}^{r-1}T_{\bk k}$}
            \draw(0,\i+.25,-.3)node{\tex};
    \end{tikzpicture}
    \begin{tikzpicture}[>=stealth,xscale=1.7,yscale=2.4]
    \pgfmathsetmacro{\v}{-3.5}
    \pgfmathsetmacro{\h}{5/9}
    \pgfmathsetmacro{\hh}{1/3}
    \draw[fill=gray!50](\v,0)--++(1,0)--++(1,1)--++(-1,0)--cycle;
    \draw[fill=gray!35](0,0)--++(\h,.333)--++(-\hh,0)--++(\h,.333) --++(-\hh,0) --++(\h,.333)--++(1,0)--++(-\h,-.3333)--++(\hh,0) --++(-\h,-.3333)--++(\hh,0)--++(-\h,-.3333)--cycle;
        \draw[->](-1.3,.7)--node[above]{$h_n$}++(.9,0);
    \foreach \i in {0,1}
        {
        \draw(\i,\i)--++(1,0);
        \draw(\i+\v,\i)--++(1,0);
        \draw(\i+\v,0)--++(1,1);
        \draw(\i,0)--++(\h,.333)--++(-\hh,0)--++(\h,.333)--++(-\hh,0) --++(\h,.333);}
    \foreach \i/\j/\tex in {\h/1/\bu_{n,1},2*\h-\hh/2/\bu_{n,2}}
        {\draw[red,fill=red](\i-\hh+.5,\j*.333+.005) circle(.02)node[below,font=\tiny]{$B_\j$};
        \draw[blue,fill=blue](\i+.5,\j*.333-.005) circle(.02)node[below,font=\tiny]{$A_\j$};
        \draw[red,->](\i+.5,.333*\j)-- node[above,font=\tiny,yshift=-3pt]{$\tex$}++(-\hh,0);}
    \foreach \i in {1,2}
        \draw[help lines](\i*\h-\i*\hh,\i*.3333)--++(0,.333);
    \foreach \i/\tex in {0/(d) $\cup_{k=1}^{r-1}T_{\bk k}$,\v/(c) $T_{\bk}$}
        {%
        \draw[fill=black](\i+.5,0)circle(.01)node[below]{$\bw(\bk)$};
        \draw[->](\i-.2,0)--++(2.3,0)node[above]{$\bx$};
        \draw[->](\i,-.1)--++(0,1.3)node[left]{$y$};
        \draw[fill=black](\i,0)circle(.02) node[left,xshift=-10pt]{$\varphi_r(\bk)$};
        \draw(\i+.5,-.3)node[below]{\tex};
        }
    \draw(\v,1)--++(-.1,0)node[left]{$\varphi_r(\bk)+r^{1-n}$};
    \draw(0,\hh)--++(-.1,0)node[left]{$\varphi_r(\bk)+r^{-n}$};
    \foreach \i in {\v}
        \draw[<-](1+\i,1)--node[above,font=\tiny]{$\bt_n$}++(-1,0);
    \foreach \i/\j in {0/1,\h-\hh/2,2*\h-2*\hh/3}
        \draw[->,very thick](\i,\j*\hh)-- node[below,font=\tiny,yshift=3pt]{$\bt_{n+1}$}++(\h,0);
    \end{tikzpicture}
    \newcommand{\myoverrightarrow}[2][18pt]
        {#2\hspace*{-#1}\raise9pt\hbox{$\longrightarrow$}\,\,\,}
        \vspace*{-10pt}
     \caption{\label{figTbk} %
        An illustration of $\bt_n$, $\{\bu_{n,k}\}$ and $\bw(\bk)$ and the homeomorphism  $h_n$ on $T_\bk$, where $n\ge1$ and $\bk\in\Sigma_r^{n-1}$. The vectors $\bu_{n,k}$ in (b) and (d) are $\myoverrightarrow[25pt]{A_k B_k}$, where $A_k$ is the center of $U+\bb(\bk (k-1)\overline{r-1})$, while $B_k$ is the center of $U+\bb(\bk k)$, $k=1,2$. The line segment $P_0P_1$ (in (a) and (b)) is in the line $\ell_{\bb(\bk)}$. The vector $\bt_n=\bb(0^{n-1}\overline{r-1}) =\bb(\bk\overline{r-1})-\bb(\bk).$ This figure is drawn with $r=3$; moreover, $d=3$ in (a) and (b), while $d=2$ in (c) and (d).}
    \end{figure}

    For $n\ge1$ and $\bk\in\Sigma_r^{n-1}$,  we can use the symbols introduced above to rewrite the expression for $T_\bk$ defined as in \eqref{eqTbk} as
        \begin{equation}\label{eqTbk2}
        \begin{array}{rcl}
            T_\bk&=&\Big\{t\Big(\bt_n,\dfrac1{r^{n-1}}\Big):0\le t\le1\Big\}+U\times\{0\}+\bw(\bk)\\[5pt]
            &=&\Big\{\Big(\dfrac{y}{r^{1-n}}\bt_n,y\Big):\,0\le y\le \dfrac1{r^{n-1}}\Big\}+U\times\{0\}+\bw(\bk).
            \end{array}
        \end{equation}
    Combining \eqref{eqwtHbi} and \eqref{eqTbk2}, we know that
        \begin{equation}\label{eqTbk3}
            T_\bk=\widetilde H+\bw(\bk)
        \end{equation}
    with $r^{1-n},\bt_n$ replacing $b,\bt$ respectively. In view of the relation between $T_\bk$ and $\widetilde H$, we let
        \begin{equation}\label{eqTbkbi}
            T_{\bk,\bi}:=\widetilde H_\bi+\bw(\bk),
            \quad \myint T_{\bk,\bi}:=\myint\widetilde  H_\bi+\bw(\bk),
            \quad \pt_v T_{\bk,\bi}:=\pt_v\widetilde  H_\bi+\bw(\bk).
        \end{equation}
    These relations allow us to define a homeomorphism  $h_n$ under the condition
        $$\|\{\bu_{1,k}\}\|_\infty<1,$$
    by using the function $h_{\bt_n,\{\bu_{n,k}\}}$ defined in \eqref{eqhomh}, with $\ep=\ep_n$, $\eta=\eta_n$ and $\bu_{k}=\bu_{n,k}$. We first point out that if all $\bu_{1,k}$, $1\le k<r$, are in $2U^\circ$, then homeomorphism  $h_{\bt_n,\{\bu_{n,k}\}}$ is well-defined. In fact, what we need to check are that $\ep_n$ satisfies \eqref{eqep} (which is obvious from the first formula in \eqref{eqcn}) and that $\|\{\bu_{n,k}\}\|_\infty<\eta_n<1$. The definition of $\eta_n$ in \eqref{eqetan}, together with the monotonicity of $\{\|\bu_{n,k}\|_\infty\}_{n\ge1}$, $1\le k<r$, implies that such an $\eta_n$ satisfies the requirement. Hence, $h_{\bt_n,\{\bu_{n,k}\}}$ is well-defined.

    Now we defined $h_n$ as follows. We let $h_n$ be an identity on  $\R^{d-1}\times(\R\backslash[0,1])$. For each slab $\R^{d-1}\times([0,r^{1-n}]+\varphi_r(\bk))$, where  $\bk\in\Sigma_r^{n-1}$, we move the origin of $\R^d$ to $\bw(\bk)$, which sends $T_\bk$ to $\widetilde H$, and then move the origin back after applying $h_{\bt_n,\{\bu_{n,k}\}}$ on $\R^d$, resulting in the map $h_n$ on the slab $\R^{d-1}\times([0,r^{1-n}]+\varphi_r(\bk))$ (see Figure \ref{figTbk}), namely,
        \begin{equation}\label{eqhn}
        h_n(P)=\begin{cases}
          P,&\cP_v(P)\notin [0,1];\\
          h_{\bt_n,\{\bu_{n,k}\}}\big(P-\bw(\bk)\big)+\bw(\bk),&
          \cP_v(P)\in[0,r^{1-n}]+\varphi_r(\bk),\bk\in\Sigma_r^{n-1}.
          \end{cases}
        \end{equation}
     We point out that if $n\ge1$ and $\cP_v(P)\in [0,r^{1-n}]+\varphi_r(\bk)$ with $\bk\in\Sigma_r^{n-1}$, then
        $$h_n(P)-\bw(\bk)=h_{\bt_n,\{\bu_{n,k}\}}\big(P-\bw(\bk)\big).$$
    From the definition of $h_n$, we obtain the following lemma.
    \begin{lem}\label{lemhn}
       Let $n\ge1$, $\|\{\bu_{n,k}\}\|_\infty<1$, and $y\in Y_n$. Then $h_n$ is an identity on the hyperplane $\pi(y)$ and maps $\R^{d-1}\times[y,y+r^{1-n}]$ onto itself. In addition, for $\bk\in\Sigma_r^{n-1}$,
            \begin{equation}\label{eqhnTbk}
                h_n(T_\bk)=\bigcup_{k=0}^{r-1} T_{\bk k}.
            \end{equation}
    \end{lem}
    \begin{proof}
        The first part of the conclusion can be obtained easily by using the definition of $h_n$; we will prove \eqref{eqhnTbk}. Fix a word $\bk\in\Sigma_r^{n-1}$. Let $h_{\bt_n,\{\bu_{n,k}\}}$ be defined as in \eqref{eqhomh} with respect to $\ep=\ep_n$ and $b=r^{1-n}$. We first list some formulas which will simply our proof. By \eqref{eqbwbk}, we get
            \begin{equation}\label{eqbwbkk}
            \bw(\bk k)-\bw(\bk)=\big(\bb(0^{n-1}k),\varphi_r(0^{n-1}k)\big)
            =\big(\bb(0^{n-1}k),y_k\big),
            \end{equation}
        where $y_k=kb/r=kr^{-n}$. By \eqref{eqbunk} and the fact that $\bb(\bk\bk')=\bb(\bk)+\bb(0^{|\bk|}\bk')$ for $\bk,\bk'\in\Sigma_r^*$, we get
        \begin{eqnarray}
            \nonumber
            \sum_{m=1}^{k}\bu_{n,m}
            &=&\sum_{m=1}^{k}\Big(\bb(0^{n-1}m)- \bb\big(0^{n-1}(m-1)\overline{r-1}\big)\Big)\\
            \label{eqbunmk}&=&\bb\big(0^{n-1}k\big) -k\bb\big(0^{n}\overline{r-1}\big).
            \end{eqnarray}
        This, together with the expressions for $\bt_n$ (see \eqref{eqbtn}), implies
        \begin{eqnarray}
            \nonumber
            \bt_{n}-\sum_{k=1}^{r-1}\bu_{n,k}
            &=&\bb(0^{n-1}\overline{r-1})-\bb\big(0^{n-1}(r-1)\big)+(r-1) \bb\big(0^{n}\overline{r-1}\big)\\
            \label{eqtntn+1}
            &=&r\bb(0^n\overline{r-1})=r\bt_{n+1}.
            \end{eqnarray}
        For $t\in[0,1]$, let $y=y_k+tr^{-n}=y_k+tb/r\in[y_k,y_k+b/r]$. It follows from \eqref{eqbunmk} and \eqref{eqtntn+1} that
            \begin{eqnarray}
            \nonumber
            t\cdot\bt_{n+1}+\bb(0^{n-1}k)
            &=&\dfrac{y}{b}\Big(\bt_n-\sum_{m=1}^{r-1}\bu_{n,m}\Big) -\dfrac{y_k}b\times r\cdot \bt_{n+1}+\bb(0^{n-1}k)\\
            \nonumber&=&\dfrac{y}{b}\Big(\bt_n-\sum_{m=1}^{r-1}\bu_{n,m}\Big) +\bb(0^{n-1}k)-k\bb(0^n\overline{r-1})\\
            \label{eqttn+1}&=&\dfrac{y}{b}\Big(\bt_n-\sum_{m=1}^{r-1}\bu_{n,m}\Big) +\sum_{m=1}^{k}\bu_{n,m}.
            \end{eqnarray}

        We now turn to the proof of \eqref{eqhnTbk}. By the relationship between $T_\bk$ and $\widetilde H$ and the expression  \eqref{eqTbk2} for $T_\bk$, we see that \eqref{eqhnTbk} is equivalent to
        \begin{equation}
          h_{\bt_n,\{\bu_{n,k}\}}(\widetilde H)\!=\! \bigcup_{k=0}^{r-1}\!\Big(\!\Big\{ t\Big(\bt_{n+1},\dfrac1{r^{n}}\Big):0\!\le\! t\!\le\! 1\!\Big\}\!+\!U\!\times\!\{0\}\!+\!\bw(\bk k)\!-\!\bw(\bk)\!\Big).\label{eqh(H)2}
        \end{equation}
        So, from \eqref{eqbwbkk} and \eqref{eqttn+1}, we see that the $k$th term in the union on the right side of \eqref{eqh(H)2} is (recall that $b=r^{1-n}$)
       \begin{eqnarray*}
       \lefteqn{\Big\{ t\Big(\bt_{n+1},\dfrac1{r^{n}}\Big):\,0\le t\le 1\Big\} +U\times\{0\}+\big(\bb(0^nk),y_k\big)}\qquad\\
          &=&\bigcup_{0\le t\le 1}
          \big(U+t\cdot\bt_{n+1}+\bb(0^nk)\big)\times \Big\{y_k+\dfrac t{r^n}\Big\}\\
          &=&\bigcup_{y_k\le y\le y_k+r^{-n}} \Big(U+\dfrac{y}{b}\Big(\bt_n-\sum_{m=1}^{r-1}\bu_{n,m}\Big) +\sum_{m=1}^{k}\bu_{n,m}\Big)\times\{y\}.
        \end{eqnarray*}
       The last line is just the $k$th term in the union in \eqref{eqh(H)}. Hence, \eqref{eqh(H)2} holds and the proof is complete.
    \end{proof}

    In the following lemma we prove some analytic properties of the quantities $\bu_{n,k}$, $\bt_n$ and $c_n$, $\bw(\bk)$. Recall that  $a:=\max\{\|\ba_k\|_\infty\}$, $s:=\|\bs\|_\infty$ and $p:=\min\{p_j:\,1\le j\le d\}>1$.
    \begin{lem}\label{lemlambda}
       Let $\bu_{n,k}$, $\bt_n$ and $c_n$, $\bw(\bk)$ be given as in \eqref{eqbunk}, \eqref{eqbtn},  \eqref{eqcn},  and \eqref{eqbwbk}. Then the following statements hold.
    \begin{enumerate}
    \item[(1)] For all $k\in\N$, $\|\bu_{n,k}\|_\infty<1$ if and only if
                \begin{equation}\label{eqmaxjk}
                \max_{\mathop{1\le k<r,}\atop 1\le j<d} \bigg\{\Big|\dfrac{a_{k+1}^{(j)}-a_k^{(j)}}{p_j^n}- \dfrac{a^{(j)}_{r-1}+s_j} {p_j^n(p_j-1)}\Big|\bigg\}<1.
                \end{equation}
    \item[(2)] $\|\bt_n\|_\infty\le p^{1-n}(a+s)$ and $\bw(\bk)$ are bounded. Moreover, for all $n$ and all $\bk\in\Sigma_r^*$, we have
            \begin{equation}\label{eqtnwk}
                \|\bt_n\|_\infty+\|\bw(\bk)\|_\infty\le 2(a+s+1).
            \end{equation}
    \item[(3)] $c_n\le p^{-n}(4ar+sr)$ and $\lambda_1:=\prod_{n=1}^{\mathbb N} (1+c_n)<+\infty$. Furthermore, if $c_1<1$, then
            \begin{gather}
                \label{eqlambda}
                \lambda_2:=\prod_{n\ge 1}(1-c_n)^{2/\eta_n}>0.
            \end{gather}
    \end{enumerate}
    \end{lem}

    \begin{proof}
        (1) It follows from the expansion for $\bb(\bk)$ that the $j$th component of $\bu_{n,k}$ is
        \begin{equation}\label{equkj}
        \omega_j(\bu_{n,k})=\dfrac{a_{k}^{(j)}-a_{k-1}^{(j)}}{p^n_j}- \dfrac{a^{(j)}_{r-1}+s_j} {p^n_j(p_j-1)}.
        \end{equation}
        Hence conclusion (1) holds.

        (2)  From \eqref{equkj} we know that
            $$|\omega_j(\bu_{n,k})| \le\dfrac{2a}{p^{n}}+\dfrac{a+s}{p^n(p-1)}
            \le\dfrac{3a+s}{p^n}.$$
        The definition of $\bb(\bk)$ implies that for $\bk\in\Sigma_r^*\cup \Sigma_r^{\mathbb N}$,
            \begin{eqnarray*}
            \|\bb(\bk)\|_\infty&\le&\|\ba(\bk)\|_\infty+s\cdot\|\bs(\bk)\|_\infty
            \\ &\le& a\sum_{n\ge1}\dfrac1{p^n} +s(r-1)\max\Big\{\sum_{n\ge1}s_n^{(j)}:\,1\le j<d\Big\}
            \\&\le&\dfrac{a+s}{p-1}\le a+s.
            \end{eqnarray*}
    Consequently, $\|\bw(\bk)\|_\infty\le \|2^{-1}\bone_{d-1}\|_\infty+\|\bb(\bk)\|_\infty+|\varphi_r(\bk)|
        \le 3/2+a+s$. Using similar arguments, we get
        \begin{gather*}
        \|\bt_n\|_\infty=\|\bb(0^{n-1}\overline{r-1})\|_\infty\le \dfrac{a+s}{p^{n-1}}.
        \end{gather*}
    Thus \eqref{eqtnwk} holds.

    (3) By \eqref{equkj} again, we know that
         $$c_n=\sum_{k=1}^{r-1}\|\bu_{n,k}\|_\infty\le (r-1)\cdot\dfrac{4a+s}{p^n}\le\dfrac{4ar+sr}{p^n}.$$
By the estimation
        $$\sum_{n\ge1}\ln(1+c_n)\le \sum_{n\ge1}c_n\le \dfrac {pc_1}{p-1}
        \le 2c_1,$$
we get $\lambda_1<+\infty$. Next, we assume $c_1<1$. Note that $\eta_n=(n+1)^{-1}$ for all $n$ sufficiently large. This implies the convergence of the series
        $$\sum_{n\ge 1}\dfrac1{\eta_n}|\ln (1-c_n)|\le
        \sum_{n\ge 1}\dfrac1{\eta_n}\Big|\ln \Big(1-\dfrac{4ar+sr}{p^n}\Big)\Big|.$$
Thus $\lambda_2>0$.
    \end{proof}

    In the rest of this paper, we assume $\lambda_1$ and $\lambda_2$ are defined as in Lemma \ref{lemlambda}(3) and let
        \begin{equation}\label{eqlambda3}
            \lambda_3:=2(a+s+2)\sum_{n\ge1}\big(c_n+\|\bt_n\|_\infty+\ep_n\big),
        \end{equation}
    which converges by Lemma \ref{lemlambda}(2,3) and the expression for $\ep_n$ (see \eqref{eqcn}). The following lemma will be used to prove the continuity of the desired homeomorphism.

    \begin{lem}\label{lemcountinuity}
        Let $h_n$ be defined as in \eqref{eqhn}. Then, for $P\in\R^d$, we have
            \begin{equation}\label{eqhndv}
            d_v\big(h_n(P),P\big)\le \ep_n,
            \end{equation}
        and
            \begin{equation}\label{eqhnP-P}
                \|h_n(P)-P\|_\infty\le c_n\|P\|_\infty+2(a+s+2)c_n+\|\bt_n\|_\infty+\ep_n.
            \end{equation}
    \end{lem}
    \begin{proof}
        Since $h_n$ is an identity on $\R^{d-1}\times\big(\R\backslash(0,1)\big)$, we may assume  $P=(\bx,y)\in\R^d\times(0,1)$. Let  $h_n(P):=(\bx',y')$, $\bv_{n}:=\sum_{k=1}^{r-1}\bu_{n,k}$ and $\bt_n':=\bt_n-\bv_n$. Furthermore, we assume $y\in[0,r^{1-n}]+\varphi_r(\bk)$ with $\bk\in\Sigma_r^{n-1}$. Hence $h_n$ can be written as
            \begin{equation}\label{eqhn2}
            h_n(P)=f_{t_{n}'} \circ f_F \circ f_T\circ f_S\circ f_{-\bt_n} \big(P-\bw(\bk)\big)+\bw(\bk).
            \end{equation}
        Among the functions on the right side, $f_F$ is the only one that can change the height of a point. From \eqref{eqvarphi} and \eqref{eqhn2}, we get
            $$|y-y'|=\big|y-\varphi(P')\big|\le \ep_n,$$
        where $P'=f_T\circ f_S\circ f_{-\bt_n}\big(P-\bw(\bk)\big).$
        This is just \eqref{eqhndv}.

        Next, we show \eqref{eqhnP-P}. Let $z=y-\varphi_r(\bk)$ and $z'=y'-\varphi_r(\bk)$. The horizontal coordinate of $h_n(P)$
        is
            $$\rho(z)\Big(\bx-\dfrac zb\cdot \bt_n-\bw(\bk)\Big)+\bx(z) +\dfrac{z'}b(\bt_n-\bv_n)+\bw(k),$$
where $b=r^{1-n}$. This implies that \setcolsep
            \begin{eqnarray*}
            \big\|\cP_h\big(h_n(P)\big)-\bx\big\|_\infty
            &\le&\Big\|\big(\brho(z)-1\big)\Big(\bx-\dfrac zb\bt_n-\bw(\bk)\Big)\Big\|_\infty\\
            &&\quad +\dfrac{|z'-z|}b\|\bt_n\|_\infty+ \Big\|\bx(z)- \frac{z'}b\bv_{n}\Big\|_\infty
            \\ &\le &c_n\big(\|\bx\|_\infty+\|\bt_n\|_\infty+\|\bw(\bk)\|_\infty+2\big) +\|\bt_n\|_\infty\\
            &\le &c_n\|\bx\|_\infty+2(a+s+2)c_n+\|\bt_n\|_\infty.
            \end{eqnarray*}
        The last line follows from \eqref{eqtnwk}.
       Hence \eqref{eqhnP-P} holds, and this completes the proof.
%
    \end{proof}

\section{Proof of Theorem \ref{theo1}}\label{S:sec5}
    This section is devoted to the proof of Theorem \ref{theo1}. To this end, we first present the definition of the homeomorphism $h_\infty$ by using the homeomorphisms $h_n$ defined in \eqref{eqhn}. The map $h_\infty$ sends a $d$-dimensional slant prism to the self-affine tile in $\R^d$. Then we show that $h_\infty$ is a continuous surjection. After that, we use three lemmas to prove the injectivity of $h_\infty$. Finally, we present the proof of Theorem \ref{theo1}.

    Let $g_0$ be the identity on $\R^d$ and let $g_n:=h_n\circ h_{n-1}\circ\cdots\circ h_1$ for $n\ge1$, where $h_n$ is defined as in \eqref{eqhn}. Now, define the desired homeomorphism  as
        \begin{equation}\label{eqh-}
            h_\infty:=\lim_{n\to+\infty}g_n=\lim_{n\to+\infty}h_n\circ\cdots\circ h_1.
        \end{equation}
    For $P\in\R^d$, let $B(P,R)$ be the closed ball with radius $R$ centered at $P$, namely,
        $$B(P,R):=\{P':\,\|P-P'\|_\infty\le R\}.$$
    \begin{lem}\label{lemhinftycontinuity}
      Let $h_\infty$ be defined as above. Then $h_\infty$ is a continuous surjection on $\R^d$.
    \end{lem}
    \begin{proof}
        We first show $h_\infty$ is a well-defined. For $P\in\R^d$, by \eqref{eqhnP-P} we see
            \begin{eqnarray*}
            \|h_n(P)\|_\infty&\le& \|P\|_\infty+\|h_n(P)-P\|_\infty \\ &\le& (1+c_n)\|P\|_\infty+2(a+s+2)c_n+\|\bt_n\|_\infty+\ep_n.
            \end{eqnarray*}
        By induction, for $n\ge1$,
        \begin{eqnarray}
             \nonumber \lefteqn{\|g_n(P)\|_\infty\le(1+c_n)\|g_{n-1}(P)\|_\infty +2(a+s+2)c_n+\|\bt_n\|_\infty+\ep_n}\\
             \nonumber &\le&(1+c_1)\cdots (1+c_n)\bigg(\|P\|_\infty+ \sum_{m=1}^n\big(2(a+s+2)c_m+\|\bt_m\|_\infty+\ep_m\big)\bigg)\\
              \label{eqgnP}&\le&\lambda_1\|P\|_\infty+\lambda_1\lambda_3.
            \end{eqnarray}
        So, for any two integers $n,m\ge0$ with $n>m$,
            \begin{eqnarray}\nonumber
             \lefteqn{\|g_n(P)-g_m(P)\|_\infty\le\sum_{i=m+1}^{n} \|g_{i}(P)-g_{i-1}(P)\|_\infty}\\
              \nonumber&\le&\sum_{i=m+1}^{n} \Big(c_i\|g_{i-1}(P)\|_\infty+2(a+s+2)c_i +\|\bt_i\|_\infty+\ep_i\Big)\\
              &\le&\sum_{i=m+1}^{n}\Big( \big(\lambda_1\|P\|_\infty+\lambda_1\lambda_3+ 2(a+s+2)\big)c_i+\|\bt_i\|_\infty+\ep_i\Big) \nonumber\\
              &\le&\lambda_1\|P\|_\infty\sum_{i>m}c_i+C\sum_{i> m} (c_i+\|\bt_i\|_\infty+\ep_i),\label{eqgngN}
            \end{eqnarray}
        where the constant $C$ is $\lambda_1\lambda_3+ 2(a+s+2)$.
        Thus $\{g_n(P)\}$ is a Cauchy's sequence and $h_\infty(P)$ exists. Hence $h_\infty$ is well-defined on $\R^d$. 

        Next, we show the continuity of $h_\infty$.
        Let $P_0\in\R^d$ and $\ep>0$ be fixed. Since $\sum_{n\ge1}(c_n+\|\bt_n\|_\infty+\ep_n)$ converges, there exists $N>0$ such that
            $$\sum_{n\ge N}(c_n+\|\bt_n\|_\infty+\ep_n)<\ep.$$
        From the continuity of each $h_n$, we obtain the continuity of  $g_N=h_N\circ \cdots \circ h_1$. Hence there exists $\delta>0$, such that $\|P-P_0\|_\infty<\delta$ implies that $\|g_N(P)-g_N(P_0)\|_\infty<\ep$.
        By \eqref{eqgngN}, if $P\in\R^d$ satisfies $\|P-P_0\|_\infty<\delta$, we have \
            \begin{eqnarray*}
            \lefteqn{\|h_\infty(P)-h_\infty(P_0)\|_\infty= \Big\|\lim_{n\to+\infty}\big(g_n(P)-g_n(P_0)\big)\Big\|_\infty} \quad\\
            &=&\lim_{n\to+\infty}\big\|g_n(P) -g_m(P)+g_m(P)-g_m(P_0)+g_m(P_0)-g_n(P_0)\big\|_\infty\\
            &\le&\ep+\lim_{n\to+\infty}\big(\|g_n(P)-g_m(P)\|_\infty +\|g_{
            m}(P_0)-g_n(P_0)\|_\infty\big)\\
            &\le&\ep+\big(\lambda_1\|P\|_\infty+\lambda_1\|P_0\|_\infty\big)\sum_{n> m}c_n+2C\sum_{n> m}\big(c_n+\|\bt_n\|_\infty+\ep_n\big)\\
            &\le&(2\lambda_1\|P_0\|_\infty+\lambda_1\delta+2C+1)\ep.
            \end{eqnarray*}
      So, $h_\infty$ is continuous on $\R^d$.

        Finally, we show $h_\infty$ is surjective. Let $N$ be sufficiently large so that $\lambda_1\sum_{n\ge N}c_n<1$. Let $P_\infty\in \R^d$ be fixed and let
            $$\tilde g_{n}:=g_{n+N}\circ g_N^{-1}=h_{n+N}\circ \cdots\circ h_{N+1},$$
        where $g_N^{-1}$ is the inverse of $g_N$. For convenience, we let $\tilde g_0$ be the identity.
        Since each $\tilde g_n$ is a homeomorphism  on $\R^d$, there exists a sequence of points $\{P_n\}$ such that $\tilde g_n(P_n)=P_\infty$. We claim that $\{P_n\}$ is bounded. Otherwise, there would exist a subsequence $\{P_{n_k}\}$ satisfying $\lim_{k\to+\infty}\|P_{n_k}\|_\infty=+\infty$. Noting that $P_\infty=\tilde g_{n_k}(P_{n_k})$, $k\ge1$, is fixed, we have
            $$\lim_{k\to+\infty}\dfrac{\|\tilde g_{n_k}(P_{n_k})-P_{n_k}\|_\infty} {\|P_{n_k}\|_\infty} \ge\lim_{k\to+\infty}\dfrac{\|P_{n_k}\|_\infty-\|P_\infty\|_\infty}{\|P_{n_k}\|_\infty} =1.$$
        On the other hands, similar to the proof of \eqref{eqgngN} with  $n,m$ replaced by $n_k,0$, we get
            $$\|\tilde g_{n_k}(P_{n_k})-P_{n_k}\|_\infty\le \lambda_1\|P_{n_k}\|_\infty\sum_{n\ge N}c_n+C\sum_{n\ge N}(c_m+\|\bt_n\|_\infty+\ep_n),$$
        which implies, by the choice of $N$, that
            $$\lim_{k\to+\infty}\!\!\dfrac{\|\tilde g_{n_k}(P_{n_k})\!-\!P_{n_k}\|_\infty}{\|P_{n_k}\|_\infty}
            \!\le\!\!\lambda_1\sum_{m\ge1}\!c_m \!+\!\!\lim_{k\to+\infty}\dfrac{C\sum_{m\ge1} \!(c_m\!+\!\|\bt_m\|_\infty\!+\!\ep_m)} {\|P_{n_k}\|_\infty}\!<\!1,$$
 a contradiction. So, $\{P_n\}$ is bounded. Let $\{P_{n_k}\}$ be a subsequence which converges to some $P_0$. For any fixed $\ep>0$, there exists $\delta=\delta(\ep)>0$ such that
            $$P\in B(P_0,\delta)\Rightarrow \|h_\infty\circ g_N^{-1}(P)-h_\infty\circ g_N^{-1}(P_0)\|_\infty<\ep$$
        by the continuity of $h_\infty$ and $g_N^{-1}$.
        By \eqref{eqhnP-P},
        $\{\tilde g_{n}\}$ (and hence $\{\tilde g_{n_k}\}$) converges uniformly (to $h_\infty\circ g_N^{-1}$) on any compact subset of $\R^d$; more precisely, for the above $\ep>0$, there exists $K>0$ such that $\|h_\infty\circ g_N^{-1}(P)-\tilde g_{n_K}(P)\|_\infty\le\ep$ for all $P\in B(P_0,\delta)$. We also assume $K$ satisfies $P_{n_K}\in B(P_0,\delta)$. So, \setcolsep
            \begin{eqnarray*}
            \lefteqn{\|\tilde g_{n_K}(P_{n_K})\!-\!\tilde g_{n_K}(P_0)\|_\infty}\\
            &=&\big\|\tilde g_{n_K}(P_{n_K})\!-\!h_\infty\circ g_N^{-1}(P_{n_K})\! +h_\infty\circ g_N^{-1}(P_{n_K})\!-\! h_\infty\circ g_N^{-1}(P_0)\!\\
            &&\quad +h_\infty\circ g_N^{-1}(P_0)\!-\!\tilde g_{n_K}(P_0)\big\|\\
            &<&3\ep.
            \end{eqnarray*}
        Hence,
            \begin{eqnarray*}
            \|h_\infty\!\circ\!g_N^{-1}(P_0)\!-\!P_\infty\|_\infty &\le&\|h_\infty\!\circ\!g_N^{-1}(P_0)\!-\!\tilde g_{n_K}(P_0)\|_\infty\!+\!\|\tilde g_{n_K}(P_0)\!-\!P_\infty\|_\infty\\
            &\le&\ep+\|\tilde g_{n_K}(P_0)-\tilde g_{n_K}(P_{n_K})\|_\infty<4\ep.
            \end{eqnarray*}
Since $\ep>0$ was arbitrary,  $h_\infty\circ g_N^{-1}(P_0)=P_\infty$. So, $h_\infty\circ g_N^{-1}$ is a surjection.  This proves that $h_\infty$ is surjective.
    \end{proof}

   In the rest of this section, we prove the injectivity of $h_\infty$. For this purpose, we let
        \begin{equation}\label{eqgnPP'}
        (\bx_{n},y_{n}):=g_n(P),\quad (\bx'_{n},y_{n}'):=g_n(P'),\quad n\ge1,
        \end{equation}
    where $P,P'\in\R^{d}$, and let
        \begin{equation}\label{eqPinfty}
       (\bx_\infty,y_\infty):=h_\infty(P),\quad (\bx_\infty',y_\infty'):=h_\infty(P').
        \end{equation}
    We point out that the quantities $y_n,y_n'$ in this section are different from the ones in Section \ref{S:section3}. The following is a key lemma to deal with the horizontal distance between $h_\infty(P)$ and the image of the vertical boundary of $T_{\emptyset,\bi}$ under the map $h_\infty$. Note that $d_h$ is not a real metric on $\R^d$, and hence it may not be valid to interchange $d_h$ and the limit as follows
        $$d_h\big(h_\infty(P),h_\infty(\pt T_{\emptyset,\bi})\big)
        =\lim_{n\to+\infty}d_h\big(g_n(P),g_n(\pt T_{\emptyset,\bi})\big).$$
We will make full use of the construction of $T_{\bk}$ to overcome this difficulty. Recall that for $E\subset\R^d$ and $y\in\R$,
        $$E|_{\pi(y)}=\{P\in E:\,\cP_v(P)=y\}.$$
    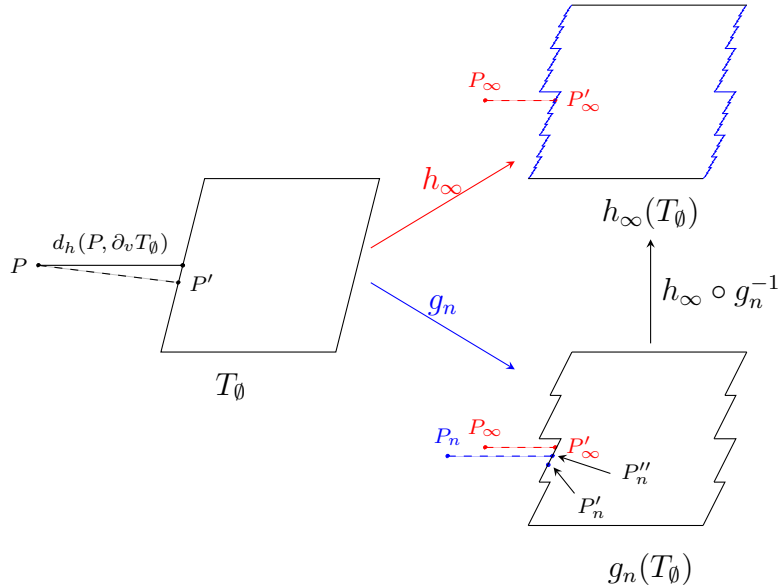
\begin{figure}[htb]
    \begin{tikzpicture}[scale=2.3,>=stealth]
        \pgfmathsetmacro{\h}{2.8}
        \pgfmathsetmacro{\hh}{.3}
        \pgfmathsetmacro{\v}{1}
        \pgfmathsetmacro{\t}{.25}
        \coordinate (p) at (.3,.5);
        \coordinate (q) at (1+.4*\t,.4);
        \coordinate (q0) at (\hh+.155,.45);
        \coordinate (p0) at ($(q0)+(-.4,0)$);
        \coordinate (qn) at (\hh+.115,.35);
        \coordinate (qnn) at (\hh+.14,.4);
        \coordinate (pn) at ($(qnn)+(-.6,0)$);
        \draw[fill=black,dashed](p)circle(.01)node[left,font=\tiny]{$P$} --(q)circle(.01) node[right,font=\tiny]{$P'$};
        \draw(1,0)--++(1,0)--++(\t,1)--++(-1,0)--cycle;
        \draw[fill=black](p)--node[above,font=\tiny] {$d_h(P,\pt_vT_\emptyset)$}++(.7+.5*\t,0)circle(.01);
        \foreach \i/\tex/\col in {1/h_\infty/red,-1/g_n/blue}
        \draw[->,\col](2.2,.5+\i*.1)--node[above]{$\tex$}++(\h/3.4,\i*\v/2);
        \foreach \i/\loc in {1,-1}
            \draw[fill=black,dashed,red]($(p0)+(\h,\i*\v)$)circle(.01) node[above,font=\tiny]{$P_\infty$} --($(q0)+(\h,\i*\v)$)circle(.01)node[right,font=\tiny]{$P'_\infty$};

        \draw[fill=black,dashed,blue]($(pn)+(\h,-\v)$)circle(.01) node[above,font=\tiny]{$P_n$} --($(qnn)+(\h,-\v)$)circle(.01);
        \draw[<-]($(qnn)+(\h+.03,-\v)$)-- ++(.3,-.1)node[right,font=\tiny]{$P_n''$};
        \draw[fill=black,blue]($(qn)+(\h,-\v)$)circle(.01);
        \draw[<-]($(qn)+(\h+.03,-\v-.015)$)--++(.12,-.15) node[below right,xshift=-3pt,yshift=3pt,font=\tiny]{$P'_n$};
        \draw[->] (\h+1,.05)--node[right]{$h_\infty\circ g_n^{-1}$}++(0,\v*.6);

        \foreach \i/\j in {{0., 0.}/{0.0117188, 0.015625}, {0.0117188,
          0.015625}/{0.0195313, 0.03125}, {0.0195313,
          0.03125}/{0.03125, 0.046875}, {0.03125,
          0.046875}/{0.03125, 0.0625}, {0.03125, 0.0625}/{0.0429688,
           0.078125}, {0.0429688, 0.078125}/{0.0507813,
          0.09375}, {0.0507813, 0.09375}/{0.0625, 0.109375}, {0.0625,
          0.109375}/{0.046875, 0.125}, {0.046875, 0.125}/{0.0585938,
           0.140625}, {0.0585938, 0.140625}/{0.0664063,
          0.15625}, {0.0664063, 0.15625}/{0.078125, 0.171875}, {0.078125,
           0.171875}/{0.078125, 0.1875}, {0.078125,
          0.1875}/{0.0898438, 0.203125}, {0.0898438,
          0.203125}/{0.0976563, 0.21875}, {0.0976563,
          0.21875}/{0.109375, 0.234375}, {0.109375,
          0.234375}/{0.0625, 0.25}, {0.0625, 0.25}/{0.0742188,
          0.265625}, {0.0742188, 0.265625}/{0.0820313,
          0.28125}, {0.0820313, 0.28125}/{0.09375, 0.296875}, {0.09375,
          0.296875}/{0.09375, 0.3125}, {0.09375, 0.3125}/{0.105469,
          0.328125}, {0.105469, 0.328125}/{0.113281, 0.34375}, {0.113281,
           0.34375}/{0.125, 0.359375}, {0.125, 0.359375}/{0.109375,
          0.375}, {0.109375, 0.375}/{0.121094, 0.390625}, {0.121094,
          0.390625}/{0.128906, 0.40625}, {0.128906,
          0.40625}/{0.140625, 0.421875}, {0.140625,
          0.421875}/{0.140625, 0.4375}, {0.140625,
          0.4375}/{0.152344, 0.453125}, {0.152344,
          0.453125}/{0.160156, 0.46875}, {0.160156,
          0.46875}/{0.171875, 0.484375}, {0.171875,
          0.484375}/{0.0625, 0.5}, {0.0625, 0.5}/{0.0742188,
          0.515625}, {0.0742188, 0.515625}/{0.0820313,
          0.53125}, {0.0820313, 0.53125}/{0.09375, 0.546875}, {0.09375,
          0.546875}/{0.09375, 0.5625}, {0.09375, 0.5625}/{0.105469,
          0.578125}, {0.105469, 0.578125}/{0.113281, 0.59375}, {0.113281,
           0.59375}/{0.125, 0.609375}, {0.125, 0.609375}/{0.109375,
          0.625}, {0.109375, 0.625}/{0.121094, 0.640625}, {0.121094,
          0.640625}/{0.128906, 0.65625}, {0.128906,
          0.65625}/{0.140625, 0.671875}, {0.140625,
          0.671875}/{0.140625, 0.6875}, {0.140625,
          0.6875}/{0.152344, 0.703125}, {0.152344,
          0.703125}/{0.160156, 0.71875}, {0.160156,
          0.71875}/{0.171875, 0.734375}, {0.171875,
          0.734375}/{0.125, 0.75}, {0.125, 0.75}/{0.136719,
          0.765625}, {0.136719, 0.765625}/{0.144531, 0.78125}, {0.144531,
           0.78125}/{0.15625, 0.796875}, {0.15625,
          0.796875}/{0.15625, 0.8125}, {0.15625, 0.8125}/{0.167969,
          0.828125}, {0.167969, 0.828125}/{0.175781, 0.84375}, {0.175781,
           0.84375}/{0.1875, 0.859375}, {0.1875,
          0.859375}/{0.171875, 0.875}, {0.171875, 0.875}/{0.183594,
          0.890625}, {0.183594, 0.890625}/{0.191406, 0.90625}, {0.191406,
           0.90625}/{0.203125, 0.921875}, {0.203125,
          0.921875}/{0.203125, 0.9375}, {0.203125,
          0.9375}/{0.214844, 0.953125}, {0.214844,
          0.953125}/{0.222656, 0.96875}, {0.222656,
          0.96875}/{0.234375, 0.984375},{0.234375, 0.984375}/{.25,1}}
        \foreach \k in {0,1}
            \draw[blue]($(\i)+(\k+\h+\hh,\v)$)--++({(1-1/32)/64},1/64) --($(\j)+(\k+\h+\hh,\v)$);
        \foreach \i/\k in {0/0,1/0.25}
            \draw[thin](\k+\h+\hh,\i+\v)--++(1,0);
      \foreach \i/\j in {{0., 0.}/{.0625, .25}, {.0625,.25}/{.0625,.5}, {.0625,.5}/{.125,.75}, {.125,.75}/{.25,1}}
        \foreach \k in {0,1}
            \draw[thin]($(\i)+(\k+\h+\hh,-\v)$)--++({(1-1/2)/4},1/4) --($(\j)+(\k+\h+\hh,-\v)$);
        \foreach \i/\k in {0/0,1/0.25}
            \draw(\k+\h+\hh,\i-\v)--++(1,0);
        \foreach \i/\tex in {{1.4,-.2}/T_\emptyset, {\h+1,
        \v*.8}/h_\infty(T_\emptyset), {\h+1,-1.25}/g_n(T_\emptyset)}
            \draw(\i)node{$\tex$};
    \end{tikzpicture}
    \caption{\label{figPQ} An illustration of $d_h\big(h_\infty(P), h_\infty(P')\big)$. $P_n=g_n(P)$, $P'_n=g_n(P')$, $P_\infty=h_\infty(P)$, and $P'_\infty=h_\infty(P')$. $P_n''$ and $P_n$ have the same height.}
 \end{figure}

    \begin{lem}\label{lemboundry>0}
        Let $h_\infty$ be defined as in \eqref{eqh-} and suppose $\bi\in\{-1,0,1\}^{d-1}$ and $P\in(\R^{d-1}\times[0,1])\backslash (\pt_vT_{\emptyset,\bi})$. If $c_1<1/4$, then
            $$d_h\big(h_\infty(P),h_\infty(\pt_vT_{\emptyset,\bi})\big)>0.$$
    \end{lem}
    \begin{proof}
        In view of the equality $h_\infty=\lim_{n\to+\infty}g_n$, we will estimate the horizontal distance between $g_n(P)$ and $g_n(\pt_v T_{\emptyset,\bi})$. Recall that $(\bx_n,y_n)=g_n(P)$. Note that $h_n$ is an identity on each hyperplane $\pi(y)$, $y\in Y_n$, which implies that $h_\infty=g_n$ on such hyperplanes. Since $g_n$ is a homeomorphism  on $\R^d$, for each $\bi\in\{-1,0,1\}$ and for each $P\in(\R^{d-1}\times[0,1])\backslash(\pt_vT_{\emptyset,\bi})$ with $y_n\in Y_n$, we have
            \begin{eqnarray*}
                d_h\big(h_\infty(P),h_\infty(\pt_v T_{\emptyset,\bi})\big)
                &=&d_h\big(h_\infty(P),h_\infty(\pt_v T_{\emptyset,\bi})|_{\pi(y_n)}\big)\\
                &=&d_h\big(g_n(P),g_n(\pt_vT_{\emptyset,\bi})\big|_{\pi(y_n)}\big),
            \end{eqnarray*}
        which is positive by the compactness of $\pt_vT_{\emptyset,\bi}$ and the continuity of $g_n$.  Next, we assume
            \begin{equation}\label{eqPYn}
            P\notin\bigcup_{n\ge1}g_n^{-1}\Big(\bigcup_{y\in Y_n} \pi(y)\Big),
            \end{equation}
        where $g_n^{-1}$ is the inverse of $g_n$. Let $P_\infty$ be given as in \eqref{eqPinfty} and let $\bk_n\in\Sigma_r^{n}$ such that $y_{n}\in(0,r^{-n})+\varphi_r(\bk_n)$. Such a $\bk_n$ does exist since \eqref{eqPYn} holds. By \eqref{eqPYn} again, for $n\ge1$, we see that
            \begin{equation}\label{eqgnP5}
            \begin{array}{l}
            d_h\big(g_{n}(P),g_{n}(\pt_vT_{\emptyset,\bi})\big) =d_h\big(g_n(P),g_n(\pt_vT_{\emptyset,\bi})|_{\pi(y_n)})\\ =d_h\big(g_n(P),\pt_v T_{\bk_n,\bi}|_{\pi(y_n)}\big),
            \end{array}
            \end{equation}
        and $T_{\bk_n}|_{\pi(y_n)}$ (the restriction of $T_{\bk_n}$ on the hyperplane $\pi(y_n)$) is a unit hypercube. Recall that
            \begin{equation}\label{eqTkn|yn}
            T_{\bk_n}\!=\!\!\!\!\!\bigcup_{0\le y-\varphi_r(\bk_n)\le r^{-n}}\!\!\!
            \Big(U+\dfrac{y-\varphi_r(\bk_n)}{r^{-n}}\bt_{n+1}\Big) \times\{y-\varphi_r(\bk_n)\}\Big\}+\bw(\bk_n).
            \end{equation}
        Before the proof, we state some facts which will help to clarify our idea. It follows from the formula $h_n(T_{\bk_{n-1}})=\cup_{k=0}^{r-1} T_{\bk_{n-1} k}$ that with respect to $h_n$, $T_{\bk_{n}}$ is a subset of the image $h_n(T_{\bk_{n-1}})$ and $T_{\bk_{n-1}}$ is the inverse image $h_n^{-1}(\cup_{k=0}^{r-1} T_{\bk_{n-1} k})$. Also, the range of height of $T_{\bk_n}$ is $[0,r^{-n}]+\varphi_r\big(\bk_n).$

    We now divided the proof into three cases: (i) $P\notin\pt_vT_{\emptyset}$, (ii) $P\in T_{\emptyset,\bi}$ and (iii) $P\in\pt_vT_{\emptyset}\backslash T_{\emptyset,\bi}$.

    \noindent \textbf{Case 1. $P\notin\pt_vT_{\emptyset}$.} In this case, Proposition \ref{theoh}(2), together with Lemma \ref{lemlambda}, implies that for $n\ge1$,
            $$d_h\big(g_n(P),g_n(\pt_vT_\emptyset)\big)
            \ge \lambda_2 d_h(P,\pt_vT_\emptyset)>0.$$
        In order to estimate the value
           \begin{equation}\label{eq5lambda2}
            d_h\big(h_\infty(P),h_\infty(\pt_vT_{\emptyset,\bi})\big)
            =d_h\big(h_\infty(P),h_\infty(\pt_vT_{\emptyset,\bi})| _{\pi(y_\infty)}\big),
            \end{equation}
        we let $P'_\infty =(\bx_\infty',y_\infty)\in h_\infty(\pt_vT_\emptyset)$ be fixed. So there exists $P'=(\bx',\by')\in \pt_vT_\emptyset$ such that $P'_\infty=h_\infty(P')$. For $n\ge1$, let
            $$(\bx_{n}',y_{n}'):=g_n(P'),$$
        which belongs to $g_n(\pt_vT_\emptyset)$.
        By Lemma \ref{lemhn}, we see that, for each $n\ge1$, $y_{n}'\in (0,r^{-n})+\varphi_r(\bk_n)$ which might not equal $y_{n}$. In order to apply the quantity
            $$d_h\big(g_n(P),g_n(\pt_vT_{\emptyset})\big) =d_h\big(g_n(P),\pt_vT_{\bk_n}|_{\pi(y_n)}\big),$$
        we will choose (based on $g_n(P')$) another point $P_n''\in\pt_v T_{\bk_n}$, which satisfies $\cP_v(P_n'')=y_n$ and the difference in norm $\|\cdot\|_\infty$ between the horizontal coordinates of $g_n(P)$ and $P_n''$ is $o(1)$ as $n\to\infty$. Since $g_{n}(P')\in\pt_v T_{\bk_{n}}$, it follows from \eqref{eqTkn|yn} that there exists $\by\in\pt U$ such that
            \begin{equation}\label{eqQn}
            \bx_n'=\dfrac{y_n'-\varphi_r(\bk_n)}{r^{-n}}\bt_{n+1} +\by+\bb(\bk_n)+2^{-1}\bone_{d-1}.
            \end{equation}
         Let $P_n''=(\bx_n'',y_n)\in \pt_vT_{\bk_n}$, where
            \begin{equation}\label{eqQn'}
            \bx_n''=\dfrac{y_n-\varphi_r(\bk_n)}{r^{-n}}\bt_{n+1} +\by+\bb(\bk_n)+2^{-1}\bone_{d-1}.
            \end{equation}
       Using \eqref{eqQn}, \eqref{eqQn'} and the fact that $y_n,y'_n\in[0,r^{-n}]+\varphi_r(\bk_n)$, we get
            \begin{eqnarray*}
            \|\bx'_n-\bx''_n\|_\infty=\dfrac{|y_n-y_n'|}{r^{-n}}\|\bt_{n+1}\|_\infty \le \|\bt_{n+1}\|_\infty=o(1).
            \end{eqnarray*}
        The last equality follows from Lemma \ref{lemlambda}(2). Next, we estimate the difference in norm $\|\cdot\|_\infty$ between the horizontal coordinates of $h_\infty(P)$ and $g_n(P)$. From \eqref{eqhnP-P}, we get
            \begin{eqnarray*}
            \|\bx_\infty -\bx_n\|_\infty&\le&\|h_\infty(P)-g_n(P)\|_\infty\\
            &=&\lim_{m\to+\infty}\|h_m\circ \cdots\circ h_{n+1}\big(g_n(P)\big)-g_n(P)\|_\infty\\
            &\le&\sum_{m\ge n}\|h_{m+1}\big(g_m(P)\big)-g_m(P)\|_\infty\\
            &\le&\sum_{m\ge n} \big(c_m\|g_m(P)\|_\infty+2(a+s+2)c_m+\|\bt_m\|_\infty+\ep_m\big)\\
            &=&o(1).
            \end{eqnarray*}
        The last equality holds because $\sum_{n\ge1}(c_n+\|\bt_n\|_\infty+\ep_n)$ converges and $\{g_n(P)\}$ is a Cauchy's sequence. Analogously, we have
            $$\|\bx'_\infty-\bx'_n\|_\infty\le\|h_\infty(P')-g_n(P')\|_\infty=o(1)\quad\text{as}\quad n\to\infty.$$
        Recall that $P_\infty=(\bx_\infty,y_\infty)$, $P'_\infty=(\bx'_\infty, y_\infty)$, and $g_n(P')=(\bx'_n,y_n')$; moreover, $g_n(P)=(\bx_n,y_n)$ and $g_n(P'')=(\bx_n'',y_n)$ have the same height. The discussion above yields, for $n$ sufficiently large,
            \begin{eqnarray*}
            \|\bx_\infty -\bx_n\|_\infty+\|\bx'_\infty-\bx_n'\|_\infty+\|\bx_n'-\bx''_n\|_\infty
            &<&\dfrac12\lambda_2 d_h(P,\pt_vT_\emptyset).
            \end{eqnarray*}
        Now, for such $n$ (see Figure \ref{figPQ}), \setcolsep
            \begin{eqnarray*}
            \lefteqn{d_h\big(h_\infty(P),P'_\infty\big) =\|\bx_\infty-\bx'_\infty\|_\infty}\qquad\\
            &=&\|\bx_\infty-\bx_{n}+\bx_{n}-\bx_n''+\bx_n''-\bx'_n +\bx'_n-\bx'_\infty \|_\infty\\
            &\ge&\|\bx_n-\bx_n''\|_\infty-\|\bx_\infty-\bx_n\|_\infty-\|\bx'_\infty -\bx'_n\|_\infty -\|\bx'_n-\bx_n''\|_\infty\\
            &>&\|g_n(P)-g_n(P'')\|_\infty-\dfrac12\lambda_2 d_h(P,\pt_vT_\emptyset)\\
            &\ge&d_h\big(g_n(P),g_n(\pt_vT_\emptyset)\big)-\dfrac12\lambda_2 d_h(P,\pt_vT_\emptyset)\ge\dfrac12 d_h(P,\pt_vT_\emptyset),
            \end{eqnarray*}
        where the last inequality follows from \eqref{eq5lambda2}. Since $P'_\infty\in h_\infty(\pt_vT_\emptyset)$ was arbitrary, we get
            $$d_h\big(h_\infty(P),h_\infty(\pt_vT_\emptyset)\big)>0.$$
        Now, for each $\bi\in\{-1,0,1\}^{d-1}$, we have $\pt_vT_{\emptyset,\bi}\subset \pt T_\emptyset$, which implies that
            $$d_h\big(h_\infty(P),h_\infty(\pt_vT_{\emptyset,\bi})\big)
            \ge d_h\big(h_\infty(P),h_\infty(\pt_vT_{\emptyset})\big)>0.$$

        \noindent\textbf{Case 2. $P\in T_{\emptyset,\bi}$.} In this case $P\in \myint T_{\emptyset,\bi}$. Proposition \ref{theoh}(3) and Lemma \ref{lemlambda} imply that
        $$d_h\big(g_n(P),g_n(\pt_vT_{\emptyset,\bi})\big)
            \ge \lambda_2 d_h(P,\pt_vT_{\emptyset,\bi}).$$
        For a fixed $P'_\infty \in h_\infty(\pt_vT_{\emptyset,\bi})$ with height $y_\infty$, we let $g_n(P')$, $g_n(P'')$ be as in Case 1, with  $T_{\emptyset,\bi}$ replacing $T_\emptyset$. A proof analogous to that for Case 1 forces $d_h\big(h_\infty(P),h_\infty(\pt_vT_{\emptyset,\bi})\big)>0$.

        \noindent\textbf{Case 3. $P\in \pt_vT_\emptyset\backslash T_{\emptyset,\bi}$.} Since $P=(\bx,y)\in \pt_v T_{\emptyset}$, by using the expression for $T_\emptyset$ in \eqref{eqTbk2}, we see that there exists $\by\in\pt U$ and $y\in[0,1]$ such that
        $$P=(\by+y\bt_1,y)+2^{-1}(\bone_{d-1},0).$$
        Let $\bj=(j_1,\ldots,j_{d-1})\in \{-1,0,1\}^{d-1}$, where the $m$th component of $\bj$ is defined as
        $$j_m=\begin{cases}
          2\omega_m(\by),&2\omega_m(\by)=\pm1;\\
          0,&2\omega_m(\by)\in(-1,1).
          \end{cases}$$
Let $(\bx_n,y_n):=g_n(P)$. In view of \eqref{eqTkn|yn} and the expression for $T_{\bk_n}$, we get
        \begin{equation}\label{eqynUj}
        \bx_n\in U_{\bj}+\bz\quad\text{and}\quad
        T_{\bk_n,\bi}|_{\pi(y_n)}=(U_{\bi}+\bz)\times\{y_n\},
        \end{equation}
    where
        $$\bz=\dfrac{y_n-\varphi_r{\bk_n}}{r^{-n}}\bt_{n+1} +\bb(\bk_n)+2^{-1}\bone_{d-1}.$$

    From $P\in \pt_vT_\emptyset$, we know that $P\in \myint T_{\emptyset,\bj}$ if $\bj$ has a zero component, while $P\in \pt_v T_{\emptyset,\bj}$ if all components of $\bj$ are nonzero. Define
        \begin{gather*}
        \Gamma_1:= \{m:\,i_m=-j_m\ne0\},\\
        \Gamma_2:=\{m:\,i_m=0,j_m\ne0\},\\
        \Gamma_3:= \{m:\,i_m\ne0,j_m=0\}.
        \end{gather*}
    Then the assumption $P\not\in T_{\emptyset,\bi}$ implies that $\Gamma_1\cup\Gamma_2\cup\Gamma_3\neq\emptyset$.
    When $\Gamma_1\ne\emptyset$, \eqref{eqynUj} says that for any $m\in\Gamma_1$ and any $\by'\in U_\bi$,
        $$\|(\bx_n-z)-\by'\|_\infty\ge |\omega_m(\bx_n-\bz)-\omega_m(\by')|
        =\Big|\dfrac12i_m-\dfrac12j_m\Big|=1.$$
    This, together with \eqref{eqgnP5}, implies that
        $$\|d_h\big(g_n(P),g_n(T_{\emptyset,\bi})\big)
        \ge \min\{\|(\bx_n-\bz)-\by'\|_\infty:\,\by'\in U_\bi\}=1.$$
    Similar to Case 1, we get
        $$\|d_h\big(h_\infty(P),h_\infty(T_{\emptyset,\bi})\big)>0.$$

    Next, we assume $\Gamma_1=\emptyset$. It follows that $\Gamma_3\ne\emptyset$, since $T_{\emptyset,\bj}\subset\pt_v T_{\emptyset,\bi}$ when $\Gamma_3=\emptyset$.
    If $\Gamma_2=\emptyset$, we see that $P\in\myint T_{\emptyset,\bj}$ and $T_{\emptyset,\bi}\subset \pt_vT_{\emptyset,\bj}$. From Case 2, we get
        $$d_h\big(h_\infty(P),h_\infty(T_{\emptyset,\bi})\big)\ge d_h\big(h_\infty(P),h_\infty(\pt_vT_{\emptyset,\bj})\big)>0.$$

    Finally, we consider the case $\Gamma_1=\emptyset$ and $\Gamma_2\ne\emptyset$. We delete all components of vectors in $\R^d$ whose indices are in $\Gamma_2$,  and obtain a system with lower dimension. For example, suppose $P=(x_1,\ldots,x_{d-1},y)\in\R^d$ and $\Gamma_2=\{d-1\}$. Then we delete the $(d-1)$th component of $P$ and  get a new point $(x_1,\ldots,x_{d-2},y)$. For convenience, we use the original notation. Then in the new system, $P\in\myint T_{\emptyset,\bj}$ and $T_{\emptyset,\bi}\subset \pt_vT_{\emptyset,\bj}$. Now similar to Case 2, we obtain
        $$d_h\big(h_\infty(P),h_\infty(T_{\emptyset,\bi})\big)\ge d_h\big(h_\infty(P),h_\infty(\pt_v T_{\emptyset,\bj})\big)>0.$$
This completes the proof.
    \end{proof}
    \begin{lem}\label{lemboundryin}
        Let $h_\infty$ be defined as in \eqref{eqh-}. If $c_1<1/4$, then $h_\infty$ is injective on $\pt_vT_\emptyset$.
    \end{lem}
    \begin{proof}
        We suppose on the contrary that $h_\infty$ is not injective on $\pt_v T_\emptyset$, i.e., there exist two distinct points $P=(\bx,y),P'=(\bx',y')$ in $\pt_vT_\emptyset$ satisfying $h_\infty(P)=h_\infty(P')$. Let $P_\infty,P'_\infty$ be given as in \eqref{eqPinfty} and assume that the coordinates of $g_n(P)$ and $g_n(P')$ are given as in \eqref{eqgnPP'}. We first prove the following claim.

      \noindent \textbf{Claim.} \textit{For any $n\ge1$, the following statements hold: $(i)$ $|y_{n}-y_{n}'|<r^{-n}$; $(ii)$ neither $y_{n}$ nor $y_{n}'$ is in $Y_n$; $(iii)$ there does not exist $y_0\in Y_n$ that lies between $y_{n}$ and $y_{n}'$.}

        To prove the claim, we note that if conclusion (i) fails, then there would exist some $N>0$ such that
            $$|y_{N}-y'_{N}|\ge\dfrac1{r^{N}}.$$
        Hence by \eqref{eqhndv} and \eqref{eqcn},
            \begin{eqnarray*}
            |y_\infty -y'_\infty |
            &=&\Big|\Big(y_N+\sum_{n\ge N}(y_{n+1}-y_n)\Big) -\Big(y_N'+\sum_{n\ge N}(y_{n+1}'-y_n')\big)\Big|\\
            &\ge&|y_{N}-y_{N}'|-\sum_{n\ge N}|y_{n+1}-y_{n}| -\sum_{n\ge N}|y_{n+1}'-y_{n}'|\\
            &\ge& \dfrac1{r^N}-2\sum_{n\ge N}\ep_{n}
            \ge\dfrac1{r^{N}}-\dfrac23\sum_{n\ge N}\dfrac1{r^{2n}}>0,
            \end{eqnarray*}
    which is impossible, since $P_\infty =P_\infty'$, $y_\infty=\cP_v(P_\infty)$, and $y'_\infty=\cP_v(P'_\infty)$.

        If conclusion (ii) fails, we assume $y_{N}\in Y_N$ for some $N$. It follows that $P_\infty =g_n(P)=g_N(P)$ for all $n\ge N$. If there exists some $N'>0$ such that $y'_{N'}\in Y_{N'}$, then $h_\infty(P')=g_{n}(P')$ for $n>\max\{N,N'\}$. The injectivity of $g_n$ implies that $P_\infty \ne P'_\infty$.
        It follows from Lemma \ref{lemhn} and conclusion (i) that for each $n>N$, there exists $\bk_n\in\Sigma_r^{n}$ such that
            \begin{equation*}
            g_n(P),g_{n}(P')\in \R^{d-1}\times\big([0,r^{-n}]+\varphi_r(\bk_n)\big)
            \end{equation*}
        and $y_{N}=\varphi_r(\bk_n)$ or $y_{N}=\varphi_r(\bk_n)+r^{1-n}$. So,
            \begin{equation}\label{eqgn}
            d_v\big(g_n(P),g_{n}(P')\big)
            =d_v\Big(g_n(P'),\pi\big(\varphi_r(\bk_n)\big)\cup \pi\big(\varphi_r(\bk_n)+r^{-n}\big)\Big)<\dfrac1{r^n}.
            \end{equation}
        Proposition \ref{theoh}(4-ii) says that
            $$d_v\big(g_{n}(P),g_n(P')\big)
            \ge d_v\big(g_{n-1}(P),g_{n-1}(P')\big)\ge\cdots\ge d_v\big(g_N(P),g_N(P')\big),$$
        which is positive, since $P'_N$ does not lie in the hyperplane $\pi(y_N)$. Letting $n$ tend to infinity, we see that  $P_\infty \ne P_\infty'$, a contradiction.

        If conclusion (iii) fails, there exist $N>0$ and $y_0\in Y_N$ that lies between $y_{N}$ and $y_{N}'$. Without loss of generality, assume $y_{N}>y_{N}'$. By Lemma \ref{lemhn} again, we see that
            $$y'_{n}<y_0<y_{n}\quad\text{ for }\quad n\ge N.$$
        From the proof for conclusion (ii), we get, for $n>N$,
            $$d_v\big(g_n(P),g_n(P')\big)\ge d_v\big(g_n(P),\pi(y_0)\big)
            \ge d_v(g_N(P),\pi(y_0)),$$
        which is impossible, since
            $$\lim_{n\to+\infty}d_v\big(g_n(P),g_n(P')\big) =d_v(P_\infty,P_\infty')=0.$$
    This proves the claim.

        We now use this claim to prove the lemma. Recall $\bone_{d-1}=(1,\ldots,1)\in\R^{d-1}$. Since $P,P'\in\pt_v T_{\emptyset}$, there exist $\by,\by'\in\pt U$ such that
            $$\bx=\by+y\bt_1+2^{-1}\bone_{d-1}\quad\text{and}\quad
            \bx'=\by'+y'\bt_1+2^{-1}\bone_{d-1}.$$
        From $\by,\by'$, we define two words $\bj=j_1\ldots j_{d-1},\bi=i_1\cdots i_{d-1}\in\{-1,0,1\}^{d-1}$ as follows
            $$
            j_m\!:=\!\begin{cases}
                2\omega_m(\by),&2\omega_m(\by)=\pm1;\\
            0,&2\omega_m(\by)\!\in\!(\!-1,1),
            \end{cases}
\quad
            i_m\!:=\!\begin{cases}
            2\omega_m(\by'),&2\omega_m(\by')=\pm1;\\
            0,&2\omega_m(\by')\!\in\!(\!-1,1),
          \end{cases}$$
        and define
        \begin{gather*}
            \Delta_1:= \{m:\,j_m=-i_m\ne0\},\\
            \Delta_2:=\{m:\,j_m=0,i_m\ne0\},\\
            \Delta_3:=\{m:\,j_m\ne0,i_m=0\}.
        \end{gather*}

    Assume $\bi\ne\bj$ so that $\Delta_1\cup\Delta_2\cup\Delta_3\ne\emptyset$. By interchanging the roles of $P$ and $P'$ if necessary, it suffices to consider the following two cases: (i) $\Delta_1\ne\emptyset$ and (ii) $\Delta_1=\emptyset$ but $\Delta_3\ne\emptyset$. However, in both cases we have shown in Case 3 of the proof of Lemma \ref{lemboundry>0} that
            $$d_h(P_\infty,P_\infty')\ge d_h\big(h_\infty(P),h_\infty(T_{\emptyset,\bi})\big)>0,$$
    a contradiction.

We now assume $\bj=\bi$. Note that either $P,P'\in\pt_v T_{\emptyset,\bj}$ or $P,P'\in\myint T_{\emptyset,\bj}$. If $P,P'\in\pt_vT_{\emptyset,\bj}$, which means that all components of $\bj$ are nonzero, we see that $\by=\by'$ and $y\ne y'$, since $P\ne P'$. The above claim, together with Proposition \ref{theoh}(4)(ii), implies that for all $n\ge1$,
            \begin{eqnarray*}
            d_v\big(g_n(P),g_n(P')\big)&\ge& d_v\big(g_{n-1}(P),g_{n-1}(P')\big)\\
            &\ge&\cdots\ge d_v(P,P')\\
            &=&|y-y'|>0,
            \end{eqnarray*}
which contradicts the hypothesis $P_\infty =P_\infty'$. Now, we show the remaining case that $P,P'\in\myint T_{\emptyset,\bj}$. The assumption $c_1<4^{-1}$ implies $c_n<4^{-1}$ for all $n$. Let $\lambda_2$ be  given as in \eqref{eqlambda}. Then
        Proposition \ref{theoh}(3) implies that there exists sufficiently large $N>0$ such that, for $n\ge N$, the horizontal distance between $g_n(P)$ and $g_n(\pt_vT_{\emptyset,\bj})$ is no less than $\lambda_2 d_h(P,\pt_vT_{\emptyset,\bj})$, which is larger than $2c_n$ since $c_n=o(1)$ by Lemma \ref{lemlambda}, i.e.,
            $$d_h\big(g_n(P),g_n(\pt_vT_{\emptyset,\bj})\big)>2c_n\quad\text{for all }n\ge N.$$
        Analogously, we have $d_h\big(g_n(P'),g_n(\pt_vT_{\emptyset,\bj})\big)>2c_n$ for $n\ge N$. Now, Proposition \ref{theoh}(4)(ii), or more precisely \eqref{eqdifhht}, says that
            $$d_v\big(g_n(P),g_n(P')\big)\ge d_v\big(g_N(P),g_N(P')\big).$$
        Letting $n$ tend to infinity, we get
        $$0=d_v\big(h_\infty(P),h_\infty(P')\big)=\lim_{n\to\infty}d_v\big(g_n(P), g_n(P')\big)\ge d_v\big(g_N(P),g_N(P')\big).$$%
       So, for $n\ge N$,  $g_n(P)$ and $g_n(P')$ have the same height. Proposition \ref{theoh}(4)(ii) or \eqref{eqdhPP'2}, together with Lemma \ref{lemdh}, implies that
            $$d_h(P_\infty ,P_\infty')=\lim_{n\to+\infty}d_h\big(g_n(P),g_n(P')\big)
            \ge\lambda_2 d_h\big(g_N(P),g_N(P')\big)>0,$$
        which is a contradiction. Hence $h_\infty$ is injective on $\pt_vT_{\emptyset}$.
    \end{proof}
    \begin{lem}\label{lemh0injective}
        Let $h_\infty$ be defined as in \eqref{eqh-} and assume $c_1<1/4$. Then $h_\infty$ is injective on $\R^d$.
    \end{lem}
    \begin{proof}
       Since $h_\infty$ is an identity on $\R^{d-1}\times(\R\backslash(0,1))$, we only study the behavior of $h_\infty$ on the region $\R^{d-1}\times[0,1]$. From Lemmas \ref{lemboundry>0} and \ref{lemboundryin}, we need to show that $h_\infty$ is injective on
            $$(\R^{d-1}\times[0,1])\backslash \pt_vT_\emptyset.$$
        For this purpose, we let $P,P'\in\R^{d-1}\times[0,1]$ be two distinct points which are not in $\pt_vT_\emptyset$. From Proposition \ref{theoh}(2), we get, for $n\ge1$,
            $$d_h\big(g_n(P),g_n(\pt_v T_\emptyset)\big)\ge \prod_{m=1}^n(1-c_m)^{2/\eta_m}d_h(P,\pt_vT_\emptyset)
            >\lambda_2 d_h(P,\pt_vT_\emptyset).$$
        By the definitions of $\{c_n\}$ and $\{\eta_n\}$,  $(1+2\eta_n+4c_n)^2-1=o(1)$ and $c_n=o(1)$. Hence, there exists an $N_1>0$ such that for $n>N_1$,
            $$d_h\big(g_n(P),g_n(\pt_v T_\emptyset)\big)> \frac12\Big(1+2\eta_n+4c_n\Big)^2-\frac12\ge c_n.$$
        Hence, all the $g_n(P)$ with $n>N_1$ have the same height. Similarly, there exists $N_2>0$ such that all the $g_n(P')$ with $n>N_2$ have the same height. Let $N=N_1+N_2$. If $g_N(P)$ and  $g_N(P')$ have different heights, we get
            $$d_v(P_\infty,P_\infty') =|\cP_v\big(g_N(P)\big)-\cP_v\big(g_N(P')\big)|>0.$$
        If $g_N(P)$ and $g_N(P')$ have the same height, then\setcolsep
            \begin{eqnarray*}
            \lefteqn{d_h(P_\infty,P_\infty')}\\
            &=& d_h\Big(\lim_{n\to+\infty} h_n\circ \cdots\circ h_{N+1}\big(g_N(P)\big), \lim_{n\to+\infty} h_n\circ \cdots\circ h_{N+1}\big(g_N(P')\big)\Big)\\
            &=&\lim_{n\to+\infty}d_h\big(h_n\circ \cdots\circ h_{N+1}\big(g_N(P)\big),h_n\circ \cdots\circ h_{N+1}\big(g_N(P')\big)\big)\\
            &\ge&\lambda_2 d_h\big(g_N(P),g_N(P')\big),
            \end{eqnarray*}
       where the second equality follows from Lemma \ref{lemdh}, and the inequality follows from \eqref{eqdhPP'2} and Lemma \ref{lemlambda}. Since each $g_n$ is injective, the right side of the last inequality above is positive. In both cases, we see that $P_\infty\ne P_\infty'$. This completes the proof that $h_\infty$ is injective on $\R^d$.
    \end{proof}

    We now return to the proof of Theorem \ref{theo1}.
   \begin{proof}[Proof of Theorem \ref{theo1}]
        Since (2) implies (1), and in \cite{Rn} we have shown that (1) is equivalent to (3), we only need to prove that (3) implies (2). By (3), we know that all $\bu_{n,k}$ belong to $2U^\circ$, so the map $h_\infty$ in \eqref{eqh-} is a well-defined continuous surjection on $\R^d$ (Lemma \ref{lemhinftycontinuity}), mapping $T_\emptyset$ onto $T$ (Lemmas \ref{lemTk} and \ref{lemhn}). If $c_1<1/4$, Lemma \ref{lemh0injective} says that $h_\infty$ is injective on $\R^d$. So  by Brouwer's domain invariance theorem (Theorem \ref{theoregion}), $h_\infty$ is a homeomorphism  on $\R^d$, which implies that $T=h_\infty(T_\emptyset)$ is a $d$-dimensional tame ball. If $c_1\ge1/4$, we can choose $N$ sufficiently large so that $c_N<1/4$ (Lemma \ref{lemlambda}(3)). For each $\bk\in\Sigma_r^{N-1}$, we apply Lemma \ref{lemh0injective} on $T_\bk$ to prove that $h_{\bk,\infty}(T_\bk)$ is a $d$-dimensional tame ball, where $h_{\bk,\infty}$ is defined as
            $$h_{\bk,\infty}(P):=
            \begin{cases}
            \dis\lim_{n\to+\infty}h_n\circ\cdots\circ h_N(P),&\cP_v(P)-\varphi_r(\bk)\in [0,r^{N-1}];\\
            P,&\text{otherwise.}
            \end{cases}$$
        Note that if two sets $h_{\bk_1,\infty}(T_{\bk_1})$ and $h_{\bk_2,\infty}(T_{\bk_2})$ intersect, the intersection is a $(d-1)$-dimensional hypercube. So
            \begin{eqnarray*}
            T&=&\Big\{\Big(\bb(\bk)+ \big(\varphi_{p_1}(\bi_1),\ldots,\varphi_{p_{d-1}}(\bi_{d-1})\big), \varphi_r(\bk)\Big):\, \bi_j\in{\Sigma_{p_j}^{\mathbb N}},\ \bk\in\Sigma_r^{\mathbb N}\Big\}\\
            &=&\bigcup_{\bk\in\Sigma_r^{N-1}}
            \Big\{\Big(\bb(\bk\bk')+ \big(\varphi_{p_1}(\bi_1),\ldots,\varphi_{p_{d-1}}(\bi_{d-1})\big), \varphi_r(\bk\bk')\Big):\,\\
             &&\qquad\qquad\ \ \bi_j\in{\Sigma_{p_j}^{\mathbb N}},\ \bk'\in{\Sigma_r^{\mathbb N}}\Big\}\\
            &=&\bigcup_{\bk\in\Sigma_r^{N-1}}h_{\bk,\infty}(T_\bk)
            \end{eqnarray*}
        is a $d$-dimensional tame ball. This completes the proof.
    \end{proof}
\section{Final comments}\label{S:questions-comments}
    In \cite{R3}, the authors of the present paper used a cut-and-paste technique to obtain a necessary and sufficient condition for a class of self-affine tiles $T\subset\R^3$ to be ball-like, under the assumption that $st\ge0$ (see Theorem \ref{theomain2}). If $st<0$, this cut-and-paste technique fails, as $T$ cannot be expressed as a limit of the union of basic blocks. Perhaps by suitably modifying the construction of basic blocks, the restriction $st\ge0$ can be removed.

    Using the techniques in this paper, i.e., re-constructing the self-affine tile by an iterating process and defining a homeomorphism  on $\R^3$ associated to each step of the iteration, we think that the following more general conjecture is true.
    \begin{conj}
        Let $(A,\D)$ be given as follows:
        \begin{equation*}
            A:=
            \begin{pmatrix}
            p&0&0\\[-3pt]
            0&q&0\\[-3pt]
            -t&-s&r
        \end{pmatrix}\end{equation*}
    and
        \begin{equation*}
        \D:=\big\{(i,j,k+a_i+b_j):\,0\le i<|p|,0\le j<|q|,\\
        0\le k<|r|\big\},
    \end{equation*}
    where $p,q,r$ are integers not less than 2 in absolute value and $s,t,a_i,b_j$ are real numbers. If for all $i,j$,
        $$\Big|\dfrac{a_{|p|-1}-a_\infty +t}{r(r-1)}+\frac{a_i-a_{i+1}}r\Big| +\Big|\dfrac{b_{|q|-1}-b_\infty +s}{r(r-1)}+\frac{b_j-b_{j+1}}r\Big| <1,$$
    then $T=T(A,\D)$ is a tame ball in $\R^3$.
    \end{conj}

    \noindent\textbf{Acknowledgements.}~~Part of this work was carried out when the first two authors were visiting the Department of Mathematical Sciences of Georgia Southern University. They thank the Department for its hospitality and support. The authors are indebted to Xianghong Chen for many helpful discussions and for pointing out the relevance of the Brouwer's invariance of domain theorem, which plays a key role in the paper. The authors thank the anonymous referee for some helpful comments and suggestions.


\begin{thebibliography}{99}

\bibitem{Aki2004} S. Akiyama and N. Gjini, On the connectedness of self-affine attractors, \emph{Arch. Math.} \textbf{82} (2004), 153--63.


\bibitem{Aki2000} S. Akiyama  and  J. Thuswaldner, Topological properties of two-dimensional number systems, \emph{J. Th\'{e}or. Nombr. Bordx.} \textbf{12} (2000), 69--79.


\bibitem{AkiThus2005} S. Akiyama  and  J. Thuswaldner, On the topological structure of fractal tilings generated by quadratic number systems, \emph{Comput. Math. Appl.} \textbf{49} (2005), 1439--1485.


\bibitem{Bandtpre}C. Bandt, Combinatorial topology of three-dimensional self-affine tiles, (preprint) http://arxiv.org/pdf/1002.0710.pdf.

\bibitem{BandtMessing}C. Bandt  and M. Mesing, Self-affine fractals of finite type, \emph{Banach Center Publ.}, \textbf{84} (2009), 131--148.

\bibitem{Bandt-Wang_2001} C. Bandt and Y. Wang,
Disk-like self-affine tiles in $\mathbb R\sp 2$,
{\em Discrete Comput. Geom.} {\bf 26} (2001), 591--601.


\bibitem{Barnsley}M.F. Barnsley, \textit{Fractals everywhere, 2nd ed.} Academic Press Professional, Boston, MA (1993). Revised with the assistance of and with a foreword by Hawley Rising, III.

\bibitem{Brouwer1912} L.E.J. Brouwer, Zur invarianz des $n$-dimensionalen gebiets, \emph{Math. Ann.} \textbf{72} (1912), 55--56.

\bibitem{Conner} G.R. Conner  and  J.M. Thuswaldner, Self-affine manifolds, \emph{Adv. Math.} \textbf{289} (2016),  725--783.

\bibitem{DengJiang2012}D.-W. Deng, T. Jiang and S.-M. Ngai, Structure of planar integral self-affine tilings, \emph{Math. Nachr.} \textbf{285}  (2012), 447--475.

\bibitem{R3}G.T. Deng, C.T. Liu and S.-M. Ngai, Topological properties of a class of self-affine tiles in $\R^3$, \emph{Trans. Amer. Math. Soc.} \textbf{370} (2018), 1321--1350.

\bibitem{Rn}G.T. Deng, C.T. Liu and S.-M. Ngai, Topological properties of a class of higher-dimensional self-affine tiles, \emph{Canad. Math. Bull.} \textbf{62} (2019), 727--740.

\bibitem{DengLau2011} Q.R. Deng and K.S. Lau, Connectedness of a class of planar self-affine tiles, \emph{J. Math. Anal. Appl.} \textbf{380}(2) (2011), 493--500.

\bibitem{Falc} K. Falconer, \textit{Fractal geometry: Mathematical foundations and applications}, John Wiley \& Sons, Inc., Hoboken, NJ, 2003.

\bibitem{Gelbrich} G. Gelbrich, Self-affine lattice reptiles with two pieces in $\R^n$, \emph{Math. Nachr.} \textbf{178} (1996), 129--134.

\bibitem{Gilbert} W.J. Gilbert, Radix representations of quadratic number fields, \emph{J. Math. Anal. Appl.} \textbf{83} (1981), 263--274.

\bibitem{Hata}M. Hata, On the structure of self-similar sets, \emph{Japan J. Appl. Math.} \textbf{2} (1985), 381--414.

\bibitem{he}X.G. He, I. Kirat and K.S. Lau, Height reducing property of polynomials and self-affine tiles, \emph{Geom. Dedicata} \textbf{152} (2011), 153--164.

\bibitem{Hutchinson} J.E. Hutchinson, Fractals and self-similarity, \textit{Indiana Univ. Math. J.} \textbf{30} (1981), 713--747.

\bibitem{Kamae}T. Kamae, J. Luo and B. Tan, A gluing lemma for iterated function systems, \emph{Fractals} \textbf{23} (2015), 1550019.

\bibitem{Katai} I. K\'{a}tai and I. K\"{o}rnyei, On number systems in algebraic number fields, \emph{Publ. Math. Debrecen} \textbf{41} (1992), 289--294.

\bibitem{Kenyon}R. Kenyon, Self-replicating tilings,  \emph{symbolic dynamics and its applications} (New Haven, CT, 1991), Contemp. Math., 135, Amer. Math. Soc., Providence, RI, 1992, 239--263.

\bibitem{Kirat2000} I. Kirat and  K.S. Lau, On the connectedness of self-affine tiles, \emph{J. Lond. Math. Soc. (2) \textbf{62} (2000)}, 291--304.

\bibitem{Kirat2004} I. Kirat, K.S. Lau and H. Rao, Expanding polynomials and connectedness of self-affine tiles, \textit{Discrete Comput. Geom.} \textbf{31} (2004), 275--286.

\bibitem{Kovacs}B. Kov\'{a}cs and A. Peth\H{o}, Number systems in integral domains, especially in orders of algebraic number fields, \emph{Acta Sci. Math. Szeged} \textbf{55} (1991), 287--299.

\bibitem{Lagarias-Wang_1996}J.C. Lagarias  and  Y. Wang, Self-affine tiles in $\R^n$, \emph{Adv. Math.} \textbf{121} (1996), 21--49.

\bibitem{Lagarias19962}J.C. Lagarias  and Y. Wang, Integral self-affine tiles in $\R^n$ I. Standard and nonstandard digit sets, \emph{J. Lond. Math. Soc.} \textbf{54} (1996), 161--179.

\bibitem{Lagarias1997}J.C. Lagarias  and Y. Wang, Integral self-affine tiles in $\R^n$~II. Lattice tilings, \emph{J. Fourier Anal. Appl.} \textbf{3} (1997), 83--102.

%
%
\bibitem{Leung-Lau_2007} K.S. Leung and K.S. Lau,
Disklikeness of planar self-affine tiles,  {\em Trans. Amer. Math.
Soc.} {\bf 359} (2007), 3337--3355.

\bibitem{Leung-Luo_2012} K.S. Leung and J.J. Luo, Connectedness of planar self-affine sets associated with non-consecutive collinear digit sets, {\em J. Math. Anal. Appl.} {\bf 395} (2012), 208--217.

\bibitem{Leung2015}K.S. Leung  and J.J. Luo, Connectedness of planar self-affine sets associated with non-collinear digit sets, \emph{Geom. Dedicata} \textbf{175} (2015), 145--157.

\bibitem{LiuJC2014}J.C. Liu, J.J. Luo and H. Xie, On the connectedness of planar self-affine sets, \emph{Chaos Soliton Fract.} \textbf{69} (2014), 107--116.

\bibitem{LiuJC2017}J.C. Liu, J.J. Luo and K. Tang, Connectedness of self-affine sets with product digit sets, \emph{Fractals} \textbf{25} (2017), 1750053.

\bibitem{Liu-Ngai-Tao_2016} J. Liu, S.-M. Ngai and J. Tao, Connectedness of a class of two-dimensional self-affine tiles associated with triangular matrices, \textit{J. Math. Anal. Appl.} \textbf{435} (2016), 1499--1513.

\bibitem{Luo-Rao-Tan_2002}J. Luo, H. Rao and B. Tan, Topological structure of self-similar sets, {\em Fractals} {\bf 10} (2002), 223--227.

\bibitem{Ma-Dong-Deng_2014} Y. Ma, X.-H. Dong and Q.-R. Deng, The connectedness of some two-dimensional self-affine sets,  {\em J. Math. Anal. Appl.} {\bf 420} (2014), 1604--1616.

\bibitem{Thu}J. Thuswaldner and S.Q. Zhang, On self-affine tiles whose boundary is a sphere, \emph{Trans. Amer. Math. Soc.} \textbf{373} (2020), 491--527.

\bibitem{wang}Y. Wang, Self-affine tiles, \emph{Advances in Wavelets} (Hong Kong, 1997), 261--282, Springer, Singapore, 1999.

\end{thebibliography}
\end{document}